\documentclass{article}
\usepackage{graphicx} 
\title{Kernel method for Möbius transformations and hyperbolic representations of groups}
\author{Gonzalo Emiliano Ruiz Stolowicz \thanks{Instituto de Ciencias Matemáticas, ICMAT-Madrid,   gonzalo.ruizstolowicz@gmail.com. The author is supported by the Spanish National Research Council (CSIC) through a Juan de la Cierva Fellowship (grant JDC2023-052425-I).}}
\usepackage[left=1.5in, right=1in, top=1in, bottom=1in]{geometry}               
\usepackage{amsthm}
\usepackage[utf8]{inputenc}
\usepackage{sectsty}   
\usepackage{amsfonts}
\usepackage{calc}
\usepackage{graphicx,import}
\usepackage{amssymb, amsmath}
\usepackage{tikz}
\usepackage{tikz-cd}
\usepackage{amsmath}
\usetikzlibrary{matrix,arrows,decorations.pathmorphing}
\usepackage{enumitem}
\usepackage{mathtools}
\usepackage{dsfont}
\usepackage{amsmath}
\usepackage{titlesec}
\titleformat*{\section}{\normalsize\bfseries}
\titleformat*{\subsection}{\normalsize\bfseries}
\usepackage[capitalise]{cleveref}
\crefname{section}{Section}{Sections}
\Crefname{section}{Section}{Sections}
\crefname{subsection}{Subsection}{Subsections}
\Crefname{subsection}{Subsection}{Subsections}

\usepackage{dirtytalk}

\numberwithin{equation}{section}

\newtheorem{teo}{Theorem}[subsection]
\newtheorem{lem}[teo]{Lemma}
\newtheorem{prop}[teo]{Proposition}
\newtheorem{cor}[teo]{Corollary}
\newtheorem*{teosin}{Theorem}
\theoremstyle{definition}

\newcommand{\N}{\mathbf N}
\newcommand{\He}{\mathbf H}

\newcommand{\R}{\mathbf R}

\newcommand{\C}{\mathbf C}

\newcommand{\F}{\mathbf F}

\newcommand {\cart}{\text{Cart}}

\newcommand{\iso}{\text{Isom}}
\newcommand{\Arg}{{\text{Arg}}}

\newcommand{\dimension}{\mathrm{dim}_\mathrm{H}}

\newcommand{\meddot}{\mathbin{\vcenter{\hbox{\scalebox{1.3}{$\cdot$}}}}}
\newcommand{\findos}{\tag*{\raisebox{-.8\baselineskip}{\qed}}
}
\newcommand{\fintres}{\tag*{\raisebox{-1.6\baselineskip}{\qed}}
}

\begin{document}

	\maketitle
	\begin{abstract}
		A new kernel method is introduced for isometric actions of groups on real or complex hyperbolic spaces that allows us to  classify  such representations for $\mathrm{PSL}_2(\R)$ and give  sufficient and necessary conditions for two  such  representations of any group to be conjugated. 
	\end{abstract}
	
	\tableofcontents
	
		\section*{Introduction}
		\addcontentsline{toc}{section}{Introduction}
		This work is situated within the study of isometric  group actions on infinite-dimensional hyperbolic spaces, see \cite{burger2005equivariant, das2017geometry, duchesne2013infinite,duchesne2021,polishtopology,gromovasymptotic,monod2018notes,monod2014exotic,monodpyselfrepresentations}. The use of \textit{kernels of positive or conditionally negative type} is classical in representation theory (see \textit{e.g.} \cite{propertyTbekka}) and it was  recently introduced in the context of infinite-dimensional hyperbolic spaces in  \cite{monod2018notes,monodpyselfrepresentations} ; we pursue this approach here. 
		 
		The main contributions of this work are the introduction of a kernel method associated to the natural actions by Möbius transformations of the boundary at infinity  induced by the  isometric representations of groups on complex hyperbolic spaces (see \cref{seccionkernelesmobius}), the   classification of isometric representations of $\mathrm{PSL}_2(\R)$ in the aforementioned spaces  (see \cref{horosphcomb}) and the following theorem (see  \cref{dispargclas}).
	\begin{teosin}
		Let $\rho$ and $\tau$ are two irreducible representations of a group $G$ in a complex hyperbolic space. If for every $g$ in $G$, the displacement of $\rho(g)$ and $\tau(g)$ are equal, then for  $g_1,g_2,g_3\in G$ such that $\rho(g_i)$ is hyperbolic and the respective attracting points $\rho(g_i)_+$ are pairwise distinct, then the attracting points $\tau(g_i)_+$ are pairwise distinct. 
		
		Moreover, if  the triples $(\rho(g_1)_+,\rho(g_2)_+,\rho(g_3)_+)$ and $(\tau(g_1)_+,\tau(g_2)_+,\tau(g_3)_+)$  have the same Cartan invariant, then $\rho$ and $\tau$ are conjugated by a holomorphic isometry. 
		\end{teosin}
	
	We will recall in \cref{generalidadeshiperbolicos} the definitions of \textit{displacement}, \textit{Cartan invariant} and \textit{holomorphic isometry}.  This theorem is a generalisation of a result  by Inkang Kim   in the context of  locally compact hyperbolic spaces (see \cite{Kimmarkedlenght}), where only the condition about the displacement is required.   It was noticed by Nicolas Monod that the non-locally compact case exhibits a different behaviour (see Remark 4.5 in \cite{monod2018notes}) and for that reason  it is  required to add  the condition on the Cartan argument.  
	
		Both the classification for $\mathrm{PSL}_2(\R)$ and the previous theorem   are obtained using the kernel method. The \textit{kernels of Möbius type} introduced in this work (see \cref{seccionkernelesmobius}) play a role analogous in relation to isometric representations of groups  in hyperbolic spaces  to that played by the  kernels of real negative type with respect to isometric  affine representations of groups in real Hilbert spaces (see \textit{e.g.} II.C in \cite{propertyTbekka}).
		
		By deforming the \textit{kernels of Möbius type} it is possible to obtain new non-equivalent representations. This   has an immediate consequence for the case of $\mathrm{PSL}_2(\R)$, the description of \say{all} the  representations of  this group  follows the outline introduced in $\cite{stolowicz2022complex}$, but with this new kernel method it is possible to enlarge the family  described in the aforementioned work and to  obtain a description of all irreducible  representations, up to conjugacy.    
		
		This is not the first attempt to develop a kernel method for the  study hyperbolic representations. In \cite{monodpyselfrepresentations} Nicolas Monod $\&$ Pierre Py  and Nicolas Monod  in \cite{monod2018notes}  introduced the notions of kernels of hyperbolic type,   real and complex, respectively. 
		Whilst these kernels \say{model} the orbit maps associated to points in the hyperbolic spaces, the kernels introduced in this work \say{model} the Möbius structure of the limit set of a (non-elementary) representation. This presents advantages that come as  consequences of the "better" dynamical properties of the induced  action by Möbius transformations on the boundary at infinity.  
		
		This work is divided in three sections. The first consists of preliminaries and  is included to fix notation and conventions. The second section studies the complex kernels of negative type and introduces a respective \textit{GNS construction} relaying  on  the Heisenberg groups (see  \cref{gnsnegativocomplejo}). Given the natural relation between the visual boundary of the complex hyperbolic spaces and the Heisenberg groups, the complex kernels of negative type and the\textit{ GNS construction} associated to them are essential for the definition of the \textit{complex kernels of Möbius type} (see \cref{seccionkernelesmobius}). And last, the theorem mentioned above is presented as a natural consequence of the \textit{GNS construction} for the \textit{kernels of Möbius type} (see \cref{GNSMobius}). 
		The third section is devoted to the study of the consequences of the \textit{deformation} of  \textit{kernels of Möbius type}, \textit{e.g.} the classification for the hyperbolic representations of $\mathrm{PSL}_2(\R)$. The study by this kernel method of  representations of non-elementary subgroups of the isometry groups of the hyperbolic spaces is included.  Sufficient conditions are given  for these representations to be infinite-dimensional and  it is shown that some \say{finiteness} properties are preserved by the \say{deformation} of representations introduced in this work. 
		
		Although the \textit{real kernels of Möbius type} are introduced in this work  and the theorem  above can be stated naturally for the real case, I would like to mention a parallel work by Bruno Duchesne $\&$ Christopher-Lloyd Simon in \cite{duchesne2026varietygroupactionsalgebraic}, that approaches the kernel method and that studies in depth the group representations on the real hyperbolic spaces.		I would like also to mention the independent work of Siwei  Liang, Yusen Long $\&$ Federico Viola in \cite{yusensl} that achieves by completely different methods the classification of the irreducible complex hyperbolic representations of $\mathrm{	PSL}_2(\R)$.
		\medskip
		
		\noindent\textbf{Acknowledgements}\\
		I would like to express my deepest gratitude to Nicolas Monod for his invaluable comments and insights, without which this project could never have begun. I would also like to thank  David Xu for his comments  as well as for bringing to my attention Kim's work, and   Pierre Py for his observations, which have undoubtedly helped improve this work.
		I am equally grateful to Yusen Long for his   valuable observations and for suggesting  that  the irreducible complex hyperbolic representations of  $\mathrm{PSL}_2(\R)$ described in this work constitute a complete list of representatives  and the possibility of achieving a classification  with the methods introduced here. Finally, I would like to thank Nicola Cavallucci  for the many discussions on the metric geometry aspects of this work, which ultimately proved to be fundamental.
	\section{Preliminaries}
	 The material presented in this section is not original and is included here just to establish the conventions and notation. A  brief introduction is given to various  concepts related to the theory of metric spaces, Heisenberg groups and hyperbolic spaces, both complex and real, that will be used later in this work. 
	\subsection{Metric spaces}\label{generalidadesmetricos}	 
	
	Let $(X,d)$ be a metric space with more than 4 points.  A quadruple $(x_1,x_2,x_3,x_4)\in X^{4}$ is called \textit{admissible} if for every $i=1,2,3$,  $x_i$ appears at most twice and for every $i\neq 4$, $x_i\neq x_4$. Denote $X^{(4)}$ the set of admissible quadruples in $X^4$. 
	Following \cite{metricmobiusgeoemtryfoeschr}, if  $(X,d)$ is a metric space,  the \textit{cross-ratio} of $X$ associated to $d$  is the map 
	$crt_d: X^{(4)}\rightarrow \mathrm{P(\R^3)}$ given by 
	\begin{equation}\label{eqcrossratio} crt_d(x_1,x_2,x_3,x_4)=\big[d(x_1,x_4)d(x_2,x_3):d(x_1,x_3)d(x_2,x_4):d(x_1,x_2)d(x_3,x_4)\big].\end{equation}
	
	Two metrics $d,d'$ on  $X$ 
	are called \textit{Möbius-equivalent} 
	if $crt_d=crt_{d'}$. A map $f:X\rightarrow X$ is called a $\textit{Möbius map}$ if it is injective and 
	for every $(x_1,x_2,x_3,x_4)\in X^{(4)}$,  
	\[crt(f(x_1),f(x_2),f(x_3),f(x_4))=crt(x_1,x_2,x_3,x_4).\]
	In \cite{buyalo2014mobius} Sergei Buyalo $\&$  Viktor Schroeder characterised the  compact metric spaces that are  Möbius-equivalent to the boundary at infinity of a rank-1 symmetric space provided with a \textit{Bourdon metric} (see also \cite{buyalomoebiusstructur,metricmobiusgeoemtryfoeschr}).
		In \cref{seccionkernelesmobius} it is given  a characterisation of the metric spaces that are   Möbius-equivalent to a subset of the visual boundary of a real or complex hyperbolic space.  
		
The following lemma, which shows that a distinguished topology is determined by a given class of  Möbius-equivalent metrics, will be used later.
	\begin{lem}\label{mobcontinua}
		Every Möbius map $(X,d)\xrightarrow{f}(X,\delta)$ is continuous.
	\end{lem}
	\begin{proof}
		Let $(y_i)$ be a sequence converging to $x\in X$ and, without lost of generality, suppose   that for every $i$, $y_i\neq x$.  Up to taking a subsequence, choose $p,q\in X\setminus\{x,y_i\}_{i\in \N}$  distinct. As $f$ is a Möbius map, 
		\[\frac{d(y_i,x)d(p,q)}{d(y_i,p)d(x,q)}=\frac{\delta(f(y_i),f(x))\delta(f(p),f(q))}{\delta(f(y_i),f(p))\delta(f(x),f(q))}.\] 
		Thus there exists $C>0$ such that for every $i$, 
		\[\frac{d(y_i,x)}{d(y_i,p)}=\frac{\delta(f(y_i),f(x))}{\delta(f(y_i),f(p))}C.\] 
		If $\delta(f(y_i),f(p))$ is unbounded, we can assume  up to taking a a subsequence that for every $i$, $\delta(f(y_i),f(x))>\delta(f(x),f(p)).$ Thus  
		\[\frac{d(y_i,x)}{d(y_i,p)}>\frac{\delta(f(y_i),f(x))}{\delta(f(y_i),f(x))-\delta(f(x),f(p))}C>C,\]which is a contradiction because $\lim\limits_{i\to\infty}\frac{d(y_i,x)}{d(y_i,p)}=0.$ This argument shows that $f$ sends bounded sets into bounded sets. 
		Thus, there exists $C'>0$ such that 
		$\tfrac{d(y_i,x)}{d(y_i,p)}>\delta(f(y_i),f(x))C',$ which shows that
		$\lim\limits_{i\to\infty}\delta(f(y_i),f(x))=0.$
	\end{proof}	
	Denote  $\Delta\subset \mathrm{P}(\R^3)$ the collection  of   ${[}x_1:x_2:x_3{]}\in\mathrm{P(\R^3)}$ such that for every $i,j=1,2,3$, $x_ix_j\geq0$  and such that 
	$|x_1|,|x_2|,|x_3|$ satisfy the triangle inequality.  
	A metric space $(X,d)$ is called  a  $\textit{Ptolemy space }$ if for every $x\in X^{(4)}$, $crt(x)\in \Delta$.   
	If $(X,d)$ is a Ptolemy space, for every $p\in X$, 
	\begin{equation}\label{definversion}d^p(y,z)=\frac{d(y,z)}{d(y,p)d(p,z)}.\end{equation}
	defines a metric in $X\setminus \{p\}$ called the $\textit{inversion of the metric d with respect to  p}.$ The metrics $d$ and $d^p$ are Möbius-equivalent in $X\setminus \{p\}$. 
	
	If $X$ is a CAT(-1) space, $x\in X$ and $\xi,\eta\in X \cup\partial X$,  denote $(\eta,\xi)_x$   the \textit{Gromov product} of $\xi$ and $\eta$ with respect to $x$ (see Lemma 3.4.10 in \cite{das2017geometry}). 
	In Theorem 2.5.1 in  \cite{Bourdonstructureconfome}, Marc Bourdon showed that 
	$d_x(\eta,\xi)=e^{-(\eta,\xi)_x}$ defines a metric  in $\partial X$ and that for  every $x, y\in X$, $d_x$ and $d_y$ are Möbius-equivalent (see p. 89 in \cite{Bourdonstructureconfome} and (\ref{dosmetricasdebourdon})). These metrics will be referred as $\textit{Bourdon metrics}$.  Thomas    Foertsch $\&$ Viktor Schroeder showed that the visual boundary of a CAT(-1) space provided with any  Bourdon metric is a Ptolemy space (see Theorem 1.1 in \cite{hyperbolicityfoerschroeder}).  
	
	Let $(X,d)$ be a (separable) metric space. The \textit{Hausdorff dimension} of $X$,  denoted as $\mathrm{dim}_\mathrm{H}(X)$,  is defined as follows.
	For every $d>0$, denote $\mathcal{U}_d=\mathcal{U}_d(X)$ the set of countable  covers of $X$ such that if  $\{V_i\}_{i\in A\subset \N}\in \mathcal{U}_d$, then for every $i\in A$,  $\mathrm{diam}(V_i)<d$.  For every $\delta>0$, denote 
	\[H^\delta_d(X)=\inf\Bigg\{\sum\limits_{i\in A}\mathrm{diam}(V_i)^\delta\mid \{V_i\}_{i\in A\subset\N}\in \mathcal{U}_d \Bigg\}. \]
	For a given cover $\{V_i\}_{i\in A\subset \N}\in \mathcal{U}_d$ of $X$, for every $\gamma>\delta$, 
	\[\sum\limits_{i\in A\subset\N}\mathrm{diam}(V_i)^\gamma\leq d^{\gamma-\delta}\sum\limits_{i\in A\subset\N}\mathrm{diam}(V_i)^\delta,\]
	and therefore, 	$H_d^\gamma(X)\leq d^{\gamma-\delta}H_d^\delta(X).$
		Denote
	$\mathcal{H}^\delta(X)=\sup_{d>0}\left\{H^\delta_d(X)\right\}$ and 
	notice that if $\delta<\gamma$ and $\mathcal{H}^\delta(X)=0$, 
	then $\mathcal{H}^\gamma(X)=0$.
	The \textit{Hausdorff dimension} of  $(X,d)$ is defined as
	\[\mathrm{dim}_\mathrm{H}(X)=\inf\{\delta\geq 0\mid\mathcal{H}^\delta(X)=0\},\]
	if there exists $\delta\geq0$ such that $\mathcal{H}^\delta(X)=0$,  or $\dim_\mathrm{H}(X)=\infty$ otherwise. For further reading on the Hausdorff dimension see \textit{e.g.} \cite{mattila1999geometry} and the references therein. 	
	The following two lemmas are immediate from the definitions.
	\begin{lem}\label{supremo}
		If $(X,d)$ is a metric space and $X=\bigcup_{i\in \N}Y_i$, then 
	$\mathrm{dim}_\mathrm{H}X=\sup_{i\in \N}\{\mathrm{dim}_\mathrm{H}(Y_i)\}.$
	\end{lem}
  
	\begin{lem}\label{semultiplicaladimensiobn}
		If $(X,d)$ is a metric space and $t\in (0,1)$, then 
		$\dimension(X,d^t)=\tfrac{1}{t}\dimension(X,d).$
	\end{lem}
		\begin{lem}\label{mobiusinvariant}
		If ${f}:X\rightarrow Y$ is a Möbius-equivalence, then $\dimension(X)=\dimension(Y).$
	\end{lem}
	\begin{proof}
		Fix two points $x_1\neq x_2$ in $X$. Define 
		$X'=X\setminus\{x_1,x_2\}$ and for every $r>0$, \[X_r=\{x\in X'\mid d(x,x_1),d(x,x_2)>1/r\}.\] For every $r$, the map  $X_r\xrightarrow{f}f(X_r)$ is  injective and   $K$-Lipschitz continuous, for some $K>0$. 
		This implies that for every $d,\delta>0$,  $H_d^\delta(X_r)\geq H_{Kd}^\delta(f(X_r))$. Thus for every $\delta$,  $\mathcal{H}^\delta(X_r)\geq \mathcal{H}^\delta(f(X_r))$, and therefore, 
		$\dimension(X_r)>\dimension(f(X_r))$. As $X'=\bigcup_{r>0}X_r$ and $f^{-1}$ has the same properties, $\dimension{X}=\dimension{Y}.$
	\end{proof}
	
	\subsection{Hyperbolic spaces}\label{generalidadeshiperbolicos}
	In this subsection only few facts are presented   that will be needed later in this work, for further reading about infinite-dimensional hyperbolic spaces see \cite{burger2005equivariant,das2017geometry, duchesne2013infinite, brunoinfinitedimension,polishtopology,duchesne2021} and for the complex hyperbolic finite-dimensional spaces  see \cite{complexhyperbolic}.	
	
	In this work $\F=\R,\C$ will denote  either the real or the complex numbers.
	Given an 
	$\F$-Hilbert space $(N, \langle\,\meddot,\meddot\,\rangle)$,  one can define  the Hermitian form $B$ on $L=\F\oplus N$  given by \[B\big(a\oplus u, b\oplus v\big)=a\overline{b}-\langle u,v\rangle.\]
	Let   $\mathbf{P}(L)$ be the set of $\F$-lines of $L$ and denote  $[v]$ the line containing $v\in L\setminus\{0\}$. 
	The $m$-dimensional $\F$-hyperbolic space, where $m=\dim_\F(N)$ is the Hamel dimension of $N$ over $\F$,  is the set 
	\[\mathbf{H}^m_\F=\{[v]\in \mathbf{P}(L)\mid B(v,v)>0\}, \]	
	provided with the metric $d$ given by 
	\[\cosh(d([v],[w]))=\tfrac{|B(v,w)|}{B(v,v)^{\frac{1}{2}}B(w,w)^{\frac{1}{2}}}.\]
	The hyperbolic spaces are characterised by  $\dim_\F(N)$.  A \say{model-free} definition  can be given considering an $\F$-vector space $L$ and a \textit{strongly non-degenerate form of signature} $(1,m)$ (see \cite{burger2005equivariant}). 
	
	The hyperbolic spaces $\mathbf{H}^m_\F$ are complete CAT(-1) spaces, thus there exists a natural notion of visual boundary, denoted  $\partial\mathbf{H}^m_\F$. The visual boundary admits at least three  models with their respective topologies: the set of equivalence classes of \textit{asymptotic} geodesic rays, the set of equivalence classes of \textit{Gromov sequences} and the set of isotropic lines in $B$. 
	
	Two geodesic rays $\sigma,\tau:[0,\infty)\rightarrow\mathbf{H}^m_\F$ are \textit{ asymptotic} if the set  $\{d(\sigma(t),\tau(t))\}_{t\geq0}$ is bounded. This is an equivalence relation and on the set of equivalence classes the \textit{cone topology} can be defined. Fix a point $x\in\mathbf{H}^m_\F$. Every geodesic ray is asymptotic to one geodesic ray issuing from $x$. The basis of this topology is given by the sets $U(\tau,R,\epsilon)$, where $\tau$ is a geodesic ray issuing from $x\in\mathbf{H}^m_\F$ and $R,\epsilon >0.$ An equivalence class $[\sigma]$, where $\sigma$ is a geodesic ray issuing from $x$, belongs to $U(\tau,R,\epsilon)$ if $d(\tau(R),\sigma(R))<\epsilon.$ The topology is independent of the choice of $x$.
	
	A sequence $(x_i)$ in $\mathbf{H}^m_\F$ is called a \textit{Gromov sequence} if $\lim_{n,m\to\infty}(x_n,x_m)_{z}=\infty$, for any   $z\in\mathbf{H}^m_\F$. Two Gromov sequences $(x_i)$ and $(y_i)$ are called \textit{equivalent} if $\lim_{n,m\to\infty}(x_n,y_m)_{z}=\infty$. In the set of equivalence classes one can define the topology associated to a Bourdon metric (see \cref{mobcontinua}). 
	
	The last model for $\partial \mathbf{H}^m_\F$ is the set of isotropic $\F$-lines.  If $x\in L$ is such that $B(x,x)>0$, and $x^\perp$ is the space of vectors $B$-orthogonal to $x$,   then for every non-trivial $y\in x^\perp$, $B(y,y)<0$, and $L=\F x\oplus x^\perp$ (see Proposition 2.1 in  \cite{burger2005equivariant}). Thus one can provide $L$ with a non-degenerate and positive definite Hermitian form $B_{\pm}$, where $B_{\pm}|_{\F x}=B|_{\F x}$ and $B_{\pm}|_{x^\perp}=-B|_{x^\perp}.$ 
	One can consider the Hilbert topology on $L$ induced by the positive definite form $B_{\pm}$. This topology   is independent of the choice of $x$ (see Lemma 2.4 in \cite{burger2005equivariant}).  If  ${\pi}:L\setminus\{0\}\rightarrow\mathbf{P}(L)$ is the quotient map, one can provide the set of isotropic lines with the subspace topology of the quotient topology induced by $\pi.$ 

	The three models and the respective topologies are  naturally  equivalent (see Theorem 1.4 in \cite{stolowicz2022real}, Proposition 3.5.3 in \cite{das2017geometry} and Proposition 3.8 in \cite{burger2005equivariant}).  
	
 A  real vector subspace $L'<L$  is called \textit{totally real} if for every $a,b\in L'$, $B(a,b)\in\R.$
	Suppose that $L$ is an $\F$-vector space and consider a  closed (for the Hilbert topology associated to a form $B_{\pm}$) $\F'$-vector subspace $L'<L$ such that $B$ restricts to a (strongly continuous) Hermitian form of signature $(1,m')$ of $L'$. If $\F=\C$,  $\F'=\C$ (resp. $\F'=\R$ and $L'$ is totally real)  and  $m'\geq1$ (resp. $m'\geq2$), and  if  $\F=\R$ and $m'\geq 2$,  then    \[\mathbf{H}^{m'}_{\F'}=\{ [v]\in\mathbf{P}(L')\mid B(v,v)>0\}\] is  isometrically   embedded in $\mathbf{H}^m_\F$. The spaces obtained in this manner  are called   \textit{hyperbolic subspaces}.
	
	If  $x,y,z\in \mathbf{H}^m_\C$ and  $\tilde{x},\tilde{y},\tilde{z}\in L$ are respective lifts, then 
	$B(\tilde{x},\tilde{y})B(\tilde{y},\tilde{z})B(\tilde{z},\tilde{x})$  has strictly positive real part. 
	The \textit{Cartan invariant} of $(x,y,z)\in(\mathbf{H}^m_\C)^3$ is defined as \[\mathrm{Cart}(x,y,z)=\Arg\big(B(\tilde{x},\tilde{y})B(\tilde{y},\tilde{z})B(\tilde{z},\tilde{x})\big).\]
	The definition of the Cartan invariant is independent of the choice of the lifts. 
	
	For a set $X$ with more than 3 elements, let  $X^{\{3\}}$ be the set of triples with pairwise distinct coordinates. 
	The Cartan invariant can be  defined in the same way for any point in  $\big(\mathbf{H}^m_\C\cup\partial\mathbf{H}^m_\C\big)^{\{3\}}.$  
	The function    ${\mathrm{Cart}}:(\partial{\mathbf{H}_\C^m})^{\{3\}}\rightarrow\left[ -\tfrac{\pi}{2},\tfrac{\pi}{2}\right]$ is an alternating 2-cocycle (see Chapter 7 in \cite{complexhyperbolic}).
	
	 A subset  $X$ of $\mathbf{H}^m_\C$ (resp. $\partial\mathbf{H}^m_\C$) is contained in a real hyperbolic subspace  (resp. the boundary of a real hyperbolic subspace) if, and only if, for every  $x\in X^3$ (resp. $X^{\{3\}}$), 
	$\mathrm{Cart}(x)=0$ (see Lemma 2.1 in  \cite{burgeriozziboundedcohomology}).  
	
	Every bijective  $\F$-linear map $L\rightarrow L$ that preserves $B$  induces an isometry $\mathbf{H}^m_\F\rightarrow\mathbf{H}^m_\F.$ Denote $\mathrm{Isom}_\F(\mathbf{H}^m_\F)$ the group of isometries obtained in this way. If $\F=\R$, $\mathrm{Isom}_\R(\mathbf{H}^m_\R)$ is the full isometry group of $\mathbf{H}^m_\R$. If $\F=\C$,  $\mathrm{Isom}_\C(\mathbf{H}^m_\C)$, the group of   \textit{holomorphic isometries},   is an index 2 subgroup of the isometry group of $\mathbf{H}^m_\C.$ The isometries missing, the \textit{anti-holomorphic  isometries},  are those  induced by invertible anti-linear  maps (see Theorem 2.3.3 in \cite{das2017geometry}). If $n<\infty$, denote  $\mathrm{PO}(1,n)=\mathrm{Isom}(\mathbf{H}^n_\R)$ and $\mathrm{PU}(1,n)=\mathrm{Isom}_\C(\mathbf{H}^n_\C)$. 
	
	The holomorphic isometries of $\mathbf{H}^m_\C$ can be characterised as the isometries that preserve the Cartan invariant, either of triples of points in $\mathbf{H}^m_\C$ or of triples of pairwise distinct points in $\partial\mathbf{H}^m_\C$. The anti-holomorphic isometries are those that change   the sign of the Cartan invariant.
	
	\begin{prop}\label{cartancontinuous} 
		The Cartan invariant ${\partial\mathbf{H}_\C^m}^{\{3\}}\xrightarrow{\mathrm{Cart}}\left[-\tfrac{\pi}{2},\tfrac{\pi}{2}\right]$ is continuous. 	  	
	\end{prop} 
	\begin{proof}
		The map $B:L^2\rightarrow\C$ is continuous, thus if $C:L^3\rightarrow\C$ is  given by \[C(a,b,c)=B(a,b)B(b,c)B(c,d),\]then  ${\Arg\circ C}:L_0^{\{3\}}\rightarrow\R$ is a  continuous function.  
		
		Let ${\pi}:L\setminus{\{0\}}\rightarrow\mathbf{P}(L)$ be the natural quotient and suppose for now that  $L$ is finite-dimensional. 
		The set of isotropic vectors $L_0\subset L$  is closed and saturated with respect to the quotient $\pi$, thus $L_0$ is locally compact and  $\pi:L_0\rightarrow\partial\mathbf{H}^m_\C$ is a quotient map. This implies that    $\pi^3:L_0^3\rightarrow(\partial\mathbf{H}^m_\C)^3$ is a quotient map (see \textit{e.g.} Theorem 3.3.17 in \cite{engelking}). The space $(L_0)^{\{3\}}\subset (L_0)^3$ is open and saturated with respect to the quotient $\pi^3$,  thus   $\pi^3:L_0^{\{3\}}\rightarrow{\partial\mathbf{H}^m_\C}^{\{3\}}$ is a quotient map.  Hence one can conclude that if $L$ is finite-dimensional,  the  map ${\mathrm{Cart}}:{\partial\mathbf{H}_\C^m}^{\{3\}}\rightarrow\left[-\tfrac{\pi}{2},\tfrac{\pi}{2}\right]$ is continuous. 
		
		Suppose that $L$ is infinite-dimensional. Recall that given $\xi\in\partial\mathbf{H}_\C^m$ and $x\in\mathbf{H}_\C^m$, a neighbourhood basis for $\xi$ for the topology of  $\partial\mathbf{H}_\C^m$ is given by the sets \[N_t(\xi)=\{\eta\in\partial\mathbf{H}^m_\C\mid (\xi,\eta)_x>t\}.\] 
		Let $(\xi_1,\xi_2,\xi_3)\in {\partial\mathbf{H}_\C^m}^{\{3\}}$ and let $\mathbf{H}_\C^n\subset\mathbf{H}_\C^{n+3}\subset \mathbf{H}_\C^m$ be  two  finite-dimensional hyperbolic subspaces such that $x\in\mathbf{H}_\C^n$ and $\xi_i\in\partial\mathbf{H}_\C^n.$  If   $T\in\mathrm{Isom}_\C(\mathbf{H}^m_\C)$ is such that $T\vert_{\mathbf{H}_\C^n}=\mathrm{Id}$, then for  every $\eta\in\partial\mathbf{H}_\C^m$, $(\xi_i,T(\eta))_x= (\xi_i,\eta)_x.$ For every \[(\eta_1,\eta_2,\eta_3)\in (\partial\mathbf{H}_\C^m)^{\{3\}}\setminus (\partial\mathbf{H}_\C^n)^{\{3\}},\] there exists $T\in\mathrm{Isom}_\C(\mathbf{H}^m_\C)$ such that $T|_{\mathbf{H}_\C^n}=\mathrm{Id}$ and such that $T(\eta_i)\in \partial\mathbf{H}_\C^{n+3}.$  This finishes the proof because the Cartan invariant is continuous on finite-dimensional hyperbolic spaces and is invariant under the diagonal action of  holomorphic isometries. 
	\end{proof}
	Using the  ideas in the proof of the previous proposition one can show that \begin{equation}\label{cartexte}\big(\mathbf{H}^m_\C\cup\partial\mathbf{H}^m_\C\big)^{\{3\}}\xrightarrow{\mathrm{Cart}}\left[-\tfrac{\pi}{2},\tfrac{\pi}{2}\right]\end{equation}
	is a continuous map. 
		

	If  $A,B:L\rightarrow L$ are two invertible $\F$-linear maps  preserving $B$ and inducing the same isometry of the hyperbolic space, then $A=\theta B $ for some unit $\theta\in \F.$ Indeed, recall that  if $x\in L$ is such that $B(x,x)>0$,   then for every non-trivial $y\in x^\perp$, $B(y,y)<0$, and $L=\F x\oplus x^\perp$. This implies that the  set of $B$-positive vectors spans $L$. Notice that if for $i=1,2$,  $z_i=r_ix\oplus v_i\in L$ is such that $B(z_i,z_i)>0$, then  $B(z_1,z_2)\neq 0.$ If  $A(z_i)=\theta_{z_i} B(z_i)$,  for some unit  $\theta_{z_i}\in\F,$   then 	
	$\theta_{z_1}=\theta_{z_2}$. 
		
	Given $\xi\in \partial \mathrm{H}^m_\F$ and a geodesic ray $\sigma$ pointing towards $\xi$,  denote   ${b_\sigma}:\mathbf{H}^m_\F\rightarrow\R$ the \textit{Busemann function} (centred at $\xi$ and normalised in $\sigma(0)$)  given by, 
	$b_\sigma(y)=\lim_{t\to\infty}d(\sigma(t),y)-t.$ If $\sigma$ and $\tau$ are  asymptotic geodesic rays, then there exists $c\in \R$ such that $b_\sigma=b_\tau+c$ (see Lemma 3.4.10 in \cite{das2017geometry}). The level (resp. sublevel) sets of the Busemann functions centred at $\xi$ are called the $\textit{horospheres}$ (resp. \textit{horoballs}) centred at $\xi$.

	Let $p\in\mathbf{H}^m_\F$ and $\tilde{p}$ be a lift of $p$ such that $B(\tilde{p},\tilde{p})=1$. For   $v\in \tilde{p}^\perp$  such that $B(v,v,)=-1$, let  $\sigma$ be the geodesic for which  the  map 
	$t\mapsto \cosh(t)\tilde{p}+\sinh(t)v  $ is a lift.
	From the definitions follows that  if $\tilde{q}$ is a lift of $q\in\mathbf{H}^m_\F$ such that $B(\tilde{q},\tilde{q})=1$,  then
	\begin{equation}\label{buseman}b_\sigma(q)=\ln(|B(\tilde{p}+v,\tilde{q})|).\end{equation}
	
	
	Any two  Bourdon metrics in $\partial\mathbf{H}^m_\F$ are Möbius-equivalent . Indeed,  let  $p,q\in\mathrm{H}^m_\F$ and $\eta,\xi\in\partial\mathbf{H}^m_\F$. Denote $B_\xi(x,y)=b_\sigma(x)-b_\sigma(y)$ and $B_\eta (x,y)=b_\tau(x)-b_\tau(y)$, where $\tau$ and $\sigma$ are any geodesic rays pointing towards  $\xi$ and $\eta$ respectively. Thus one gets from the definitions that  
	\begin{equation}\label{dosmetricasdebourdon}
		d_p(\eta,\xi)=d_q(\eta,\xi)\exp\big(\tfrac{1}{2}(B_\xi(q,p)+B_\eta(q,p)) \big).
	\end{equation}
	
	If $A$ is an isometry of $\mathbf{H}^m_\F$, then    with a slight  abuse of notation, $A$ induces a Möbius map $A:\partial\mathbf{H}^m_\F\rightarrow\partial\mathbf{H}^m_\F$. Indeed,  for every $x\in\mathbf{H}^m_\F$ and $\xi,\eta\in\partial\mathbf{H}^m_\F$,  $d_x(A(\xi),A(\eta))=d_{A^{-1}(x)}(\xi,\eta).$
		In Theorem 0.1 in \cite{surlebirapportBourdon}, Marc Bourdon showed that if $X$ is a CAT(-1) space, then  for every $0<n<\infty$,  if $f:\partial \mathbf{H}_\F^n\rightarrow\partial X$ is an injective Möbius map, then $f$ is induced by an isometric map 
	$\mathbf{H}^n_\F\rightarrow X.$ 
	
	\begin{lem}\label{extensionbourdoninfinito}
		If  $\partial\mathbf{H}^m_\F\xrightarrow{T}\partial\mathbf{H}^m_\F$ is a bijective Möbius map, then $T$ is induced by a unique isometry.
	\end{lem} 
	\begin{proof}
		Let  $\mathbf{H}$ be a finite-dimensional hyperbolic subspace.    There exists   
		$ R_\mathbf{H}:\mathbf{H}\rightarrow\mathbf{H}^m_\F$ an isometric map inducing $T:\partial \mathbf{H}\rightarrow\partial \mathbf{H}_\F^m$. The claim is  
		that $R_\mathbf{H}$ is the unique isometric map  with that property.  Indeed, by contradiction suppose that there exist $S\neq R_\mathbf{H}$, an isometric  map with the same  property. Let  $x\in\mathbf{H}^m_\F$ be such that
		$R_\mathbf{H}(x)\neq S(x).$ Fix $\xi\in\partial\mathbf{H}$ such that $T(\xi)$ is not a limit of the geodesic containing $R_\mathbf{H}(x)$ and $S(x)$. Let $\tau$ be the geodesic in $\mathbf{H}$ containing $x$ and  with $\xi$ as one of its extremes. Thus $R_\mathbf{H}\circ \tau$ and $S\circ\tau$ are two   geodesics with disjoint images   sharing the extreme points, but this is a contradiction. 
		
	The space $\mathbf{H}^m_\C$ is the union of all the finite-dimensional hyperbolic subspaces and observe that for any   finite-dimensional hyperbolic subspace $\mathbf{K}$, there exists a finite-dimensional hyperbolic space $\mathbf{L}$ such that $\mathbf{H}\cup\mathbf{K}\subset\mathbf{L}$, and  $R_\mathbf{L}|_\mathbf{H}=R_\mathbf{H}$ and $R_\mathbf{L}|_\mathbf{K}=R_\mathbf{K}$.  \end{proof}
	
	For $g\in\mathrm{Isom}_\F(\mathbf{H}^m_\F)$, define the \textit{displacement}  as $\ell(g)=\inf\{d(gx,x)\}_{x\in \mathbf{H}^m_\F}.$ If $\ell(g)>0$, $g$ is called \textit{hyperbolic},   $g$ fixes two points in $\partial \mathbf{H}^m_\F$, one \textit{attracting point} $g_+$ and one \textit{repelling point} $g_-$,  and $\ell(g)$ is achieved in the geodesic connecting these two points.  If $\ell(g)=0$, either $\ell(g)$ is achieved and $g$ is called \textit{elliptic}, or  $g$ is called \textit{parabolic} otherwise. The parabolic isometries fix a unique point in $\partial\mathbf{H}^m_\F$ and  preserve the {horospheres} centred at this point (see  \textit{e.g.}  Theorem 6.1.4 in \cite{das2017geometry}).
	
	The following lemma, that will be used later in this work, establishes a Cauchy-Schwarz-type inequality and  shows that the $B$-unitary maps are bounded operators. 
	\begin{lem}\label{controldelanorma}
		Fix $x\in L$ such that $B(x,x)=1$ and consider $\|\,\meddot\,,\meddot\,\|^2$ the norm associated to the Hermitian positive definite form $B_{\pm}$ associated to the  decomposition $L=\F x\oplus x^\perp.$ If  $v,w\in L$ and $A$ is an  $\F$-linear map preserving $B$, then    the following hold.
		\begin{enumerate}
			\item $\|v\|^2=2|B(v,x)|^2-B(v,v)$.
			\item $|B(v,w)|\leq \|v\|\|w\|.$
			\item $\|A\|^2\leq 2(\|A^{-1}x\|^2+\|x\|)^2+1.$
		\end{enumerate}
	\end{lem}
	\begin{proof}
	For	1. and 2., write  $v=B(v,x)x+u$ and $w=B(w,x)x+u'$, for some $u$ and $u'$ in $ x^\perp$.		
		For 3. observe that   for every $v\in L$, 
\[			\|A(v)\|^2=2|B(A(v),x)|^2-B(A(v),A(v))
					\leq\|v\|^2\big(2\|A^{-1}x\|^2+2\|x\|^2+1\big). 
		\qedhere\]
			\end{proof}

Suppose  $\F=\C$. Given $\xi_1,\xi_2\in L$ such that $B(\xi_i,\xi_i)=0$ and $B(\xi_1,\xi_2)=1$, there exists a \textit{horospherical parametrisation} $\sigma:\R^2\times \xi_1^\perp\cap\xi_2^\perp\rightarrow\mathbf{H}_\C^n$
	given by \begin{equation}\label{horosphericalcoordinates}\sigma(s,b,v)=\left[\left(\tfrac{1}{2}(e^s-e^{-s}B(v,v))+ie^{-s}2b\right)\xi_1+e^{-s}\xi_2+e^{-s}v\right].\end{equation}
	Every $q\in \partial\mathbf{H}^m_\F\setminus\{[\xi_1]\}$ admits a unique  representative  of the shape 	\begin{equation}\label{notacionfrontera}\xi(v,b)=K(v,b)\xi_1+\xi_2+{v},\end{equation} where $v\in \xi_1^\perp\cap\xi_2^\perp$, $b\in\R$,  and  $K(v,b)=-\tfrac{{B(v,v)}}{2}+ib$. Observe that the following hold.
	\begin{enumerate}	\item The map $s\mapsto\sigma (s,b,v)$  is a geodesic connecting $[\xi(v,b)]$  with $[\xi_1].$ 
		\item The set $\sigma\big(\{s\}\times\R\times\xi_1^\perp\cap\xi_2^\perp\big)$  is the horosphere centred at $[\xi_1]$ that contains $\sigma(s,0,0)$ .\end{enumerate}
	Notice that  that 
	\begin{equation}\label{distanciahoroesferica}\begin{array}{rcl}
			\cosh\left(d\big(\sigma(s,v,b),\sigma(s,w,d)\big)\right)&=&\\
			e^{-2s}	\left|B\Big(\left(\tfrac{1}{2}e^{2s}+K(v,b)\right)\xi_1+\xi_2+v,\left(\tfrac{1}{2}e^{2s}+K(w,d)\right)\xi_1+\xi_2+w\Big)\right|&=&\\
			e^{-2s}	\left|e^{2s}+B(\xi(v,b),\xi(w,d))\right|.
	\end{array}\end{equation}
	When $\F=\R$ the same formulas are valid if one declares $b=0$. 
	\subsection{Heisenberg groups and the complex hyperbolic spaces}\label{heisengroup}
	
	Let $(H, \langle\,\meddot\,,\meddot\,\rangle)$ be a complex Hilbert space with $m=\mathrm{dim}_\C(H)$, the Hamel dimension of $H$ over $\C$.  Define the   $m$-Heisenberg group  as $\mathbf{H}_m=H\times \R$ provided with the group operation 
	\[(v,b)\cdot(w,c)=\left(v+w,b+c-2\mathrm{Im}\big(\langle v,w\rangle\big)\right).\]
	 Observe that $e=(0,0)$ and $(v,b)^{-1}=(-v,-b)$. 
	For $(u,b),(w,c)\in\mathbf{H}_m$, define 
	\begin{equation}\label{Bk}B_k\big((v,b),(w,c)
		\big)= \|v-w\|^2+i\big(b-c -2\mathrm{Im}(\langle v, w\rangle)  \big).\end{equation}	
	The  \textit{Kor\'{a}nyi metric} on $\mathbf{H}_m$ is defined  as   
	$d_k=|B_k|^\frac{1}{2}$ (see
	   \cref{isometriakoran}).	The following lemma follows from the definitions and  shows that $d_k$ is a left-invariant metric. 
	\begin{lem}\label{esinvariante}
		The map $B_k$ is invariant under the left action of $\He_m$ on itself.
	\end{lem}
	Let $m\geq2$ and consider two distinct   $p,q\in \partial\mathbf{H}^m_\C$  and  respective representatives $\xi_1,\xi_2\in L$ such that $B(\xi_1,\xi_2)=1$. 
	There is a  bijection ${\Gamma}:\partial\mathbf{H}^m_\C\setminus\{p\}\rightarrow\mathbf{H}_{m-1}, $ given by
		$\Gamma(\left[ \xi(v,b) \right])=\big(\tfrac{v}{\sqrt{2}},b\big)$, 
	where $\left(\xi_1^\perp\cap\xi_2^\perp,-B\right)$ is identified with 
	$\C^{m-1}\times\{0\}\subset \mathbf{H}_{m-1}$  (see (\ref{notacionfrontera})). 
	
	Denote $(\mathbf{H}_0,d_k)=(i\R,d^{\frac{1}{2}}_{Euc})
	$ and $\xi(b)=ib\xi_1+\xi_2$. There exists a  bijection ${\Gamma}:\partial\mathbf{H}^1_\C\setminus\{p\}\rightarrow\mathbf{H}_{0}$
	given by $\left[\xi(b)\right]\mapsto ib.$
	In this case define  $B_k(ib,id)=i(b-d).$ 
	
	For $v\in\xi_1^\perp\cap\xi_2^\perp$ and $b\in\R$, let  $g(v,b)\in\mathrm{Isom}_\C(\mathbf{H}^m_\C)$ be   such that $g(v,b)[\xi_1]=[\xi_1]$ and \begin{equation}\label{isometiagvb}g(v,b)[\xi(w,d)]=\big[\xi\big(v+w,b+c +\mathrm{Im}(B( v,w) )\big) \big].\end{equation}
	Denote 
	$x=\left[\xi_1+\xi_2\right]\in\mathbf{H}^m_\C$ and consider  $d^p_x$ the inversion with respect to $p=[\xi_1]$ of the Bourdon metric centred at $x$ (see \cref{generalidadesmetricos}). 
	\begin{prop}\label{isometriakoran}
		The  bijection $\partial\mathbf{H}^m_\C\setminus\{p\}\xrightarrow{\Gamma}\mathbf{H}_{m-1}$ is such that,
		\begin{enumerate}\item	$\Gamma \big(g(v,b)[\xi(w,c)]\big) =\Gamma(\xi(v,b))\cdot\Gamma(\xi(w,c))\big).$
			\item $B\big(\xi(v,b),\xi(w,c)\big)={B_k\big(\Gamma(\xi(v,b)),\Gamma(\xi(w,c))\big)}.$	
			\item	$d^{p}_x\big([\xi(v,b)],[\xi(w,c)]\big)=d_k\left(\Gamma\left(\xi(v,b)\right),\Gamma\left(\xi(w,c)\right) \right).$\end{enumerate}
		When $m=1$, the convention 	 $\xi(b)=\xi(0,b)$ and $(0,b)=ib$ is used. 
	\end{prop}
	\begin{proof}
		The points 1. and 2. follow from the definitions. Denote $\tilde{x}=\frac{1}{\sqrt{2}}(\xi_1+\xi_2)$ and notice that for every $v\in \xi_1^\perp\cap\xi_2^\perp$ and $b\in\R$,  there exist a unique  $\lambda=\lambda(v,b)\in \C$ and  a unique $u=u(v,b)\in \tilde{x}^\perp$, with  $B(u,u)=-1$, such that $\lambda(\tilde{x}+u)=\xi(v,b)$. Indeed, 
		if $u\in \tilde{x}^\perp$ and $B(u,u)=-1$, then the map ${\sigma}:\R_{\geq0}\rightarrow L$, given by  $\sigma(t)= \cosh(t)\tilde{x} +\sinh(t)u,$ is a lift of a generic geodesic ray issuing from $x$ and  pointing towards $[\tilde{x}+u].$ Let $\lambda_0\in\C$ and $u_0\in \tilde{x}^\perp$ be such that $\lambda_0(\tilde{x}+u_0)=\xi_1.$ 
		
		The condition $B(\tilde{x},u)=0$ implies that  $u=r\xi_1-r\xi_2+w$, for some $r\in\C$ and $w\in \xi_1^\perp\cap\xi_2^\perp$.  	In particular, observe that $u_0=r_0\xi_1-r_0\xi_2$, for some $r_0\in \C.$
		Notice that as,  \[\lambda_0\left(\tfrac{1}{\sqrt{2}} +r_0\right)\xi_1+\lambda_0\left(\tfrac{1}{\sqrt{2}} -r_0\right)\xi_2=\xi_1,\]
		then $\lambda_0=\tfrac{1}{\sqrt{2}}.$

		If $u,w\in \tilde{x}^\perp$ are such that $B(u,u)=-1=B(v,v)$, then  
		\begin{equation*}\label{prodgromov}
			{-2([\tilde{x}+u],[\tilde{x}+w])_{x}}=
			\lim\limits_{t\to\infty}\left(\mathrm{arcosh}\left(\left|\cosh(t)^2+\sinh(t)^2B(u,w)\right|\right) -2t\right).\end{equation*}
		Thus, 			
		\begin{equation}\label{representantes}
			\begin{array}{rcl}
				\exp\left({-2\big([\tilde{x}+u],[\tilde{x}+w]\big)_{x}}\right)
				&=&\\	
				\lim\limits_{t\to\infty}e^{-2t}\Big(\left|\cosh(t)^2+\sinh(t)^2B(u,w)\right|+ \sqrt{\left|\cosh(t)^2+\sinh(t)^2B(u,w)\right|^2-1}\Big)&=&\\
				\frac{1}{2}|1+B(u,w)|&=&\\
				\frac{1}{2}\left|B(\tilde{x}+u,\tilde{x}+v) \right|
								\end{array} \end{equation}
		Therefore, \begin{equation}\label{primerade}\begin{array}{rcl}
				\exp\Big(-2\big([\xi(v_1,b_1)],[\xi(v_2,b_2)]\big)_{x}\Big)&=&
				\frac{1}{2}\left|B\Big(\tilde{x}+u(v_1,b_1),\tilde{x}+u(v_2,b_2)\Big)\right|	\\&=&
				\frac{1}{2}|\lambda(v_1,b_1)|^{-1}|\lambda(v_2,b_2)|^{-1}\left|B\big(\xi(v_1,b_1),\xi(v_2,b_2)\big)\right|
			\end{array}
		\end{equation}
		and 
		\begin{equation}\label{prim}
				\exp\big(-2\big([\xi(v,b)],p\big)_{x}\big)= \tfrac{1}{2}|\lambda_0|^{-1}|\lambda(u,b)|^{-1}|B\big(\xi(v,b),\xi_1\big)|=
				\tfrac{1}{\sqrt{2}}|\lambda(u,b)|^{-1}.		\end{equation}
		The identities (\ref{primerade}), (\ref{prim}), the point 2. and the fact that $ \lambda_0=\tfrac{1}{\sqrt{2}}$  imply  that
		\begin{equation}\label{equa}
				d^p_x\big([\xi(v,b)],[\xi(w,c)]\big)^2=
				{\left|B\big(\xi(v,b),\xi(w,c)\big)\right|}=	d_k\big(\Gamma(\xi\left({v},b\right)),\Gamma(\xi\left({w},c\right))\big)^2.
		\end{equation}
	\end{proof}
	With the convention $\xi(v)=\xi(v,0)$, one can prove an  analogous theorem for the real hyperbolic space. In this case one considers the identification  $[\xi(v)]\mapsto \frac{v}{\sqrt{2}},$
	where the metric on $\xi_1^\perp\cap\xi_2^\perp$ is   given by  the restriction of the form $-B$. Hence, one can show in this case that  
	$d^p_x(\xi(v),\xi(w))^2=d\big(\tfrac{v}{\sqrt{2}},\tfrac{w}{\sqrt{2}}\big).$
	The following lemma is an immediate consequence of the definition of $B_k$. 
	\begin{lem}\label{linearrep}	Let $\He_m=H\times\R$ be the $m$-Heisenberg group. If   $\He_m\xrightarrow{T}\He_m$ is  a map   such that  $T(e)=e$ and such that  for every $x,y\in \He_m$, $B_k(T(x),T(y))=B_k(x,y),$ then there exists a unitary map $H\xrightarrow{\pi}H$ such that 
		$T=\pi\times Id.$
	\end{lem}
	\section{Kernels and hyperbolic representations}\label{seccionthekerneles}
		
	In \cref{kernelesnegativos} 
	the real and  complex \textit{kernels of negative type} will be discussed in the context of isometric  representations on hyperbolic spaces. Although these kernels  have been extensively studied, the novelty presented here regarding complex \textit{kernels of negative type} is perhaps the introduction of a \textit{GNS construction}  involving the  Heisenberg groups.
	These kernels and the aforementioned \textit{GNS construction} are essential for  the  method introduced in this work. 
	  For instance, using \say{deformations} and the \textit{Lévy-Khintchine representation} of the complex \textit{functions of negative type}, a classification for the hyperbolic representations of $\mathrm{PU}(1,1)$ is obtained in \cref{algunosejemplos}. 
	
	In \cref{seccionkernelesmobius}  the problem of characterising the embeddings of a set $X$ into the spaces $\partial \mathbf{H}^m_\F$  is addressed and by doing so the  new kernel method is introduced. A family of maps, the \textit{kernels of Möbius type}, is put in bijection with such embeddings, a method to deform this functions is described and a \textit{GNS construction} is given. In other words,  the necessary and sufficient conditions will be described  for a cross-ratio defined on a set $X$ to be the cross-ratio associated to an embedding of such set into the visual boundary of a real or complex hyperbolic space. 
	
	The idea of  studying hyperbolic representations with the use of  kernels can be traced to the work of Pierre Py and Nicolas Monod (see   \cite{monod2018notes,monodpyselfrepresentations}). 
	In the first work    a notion of \say{exponentiation} of complex kernels of hyperbolic type is introduced and with it  Nicolas Monod described a family of  non-trivial infinite-dimensional complex hyperbolic representations of $\mathrm{PU}(1,n)$.
	
	The kernels of hyperbolic type \say{model} the maps of a set into a real or complex hyperbolic space. In their equivariant version these kernels  \say{model}  the  orbit maps of points contained in the hyperbolic space. With  the \textit{kernels of Möbius type} the study of hyperbolic representations is done by looking at the induced action by Möbius transformations on the visual boundary of the hyperbolic spaces. 
	
	In \cref{seccioncomplete}, with the language of   \textit{kernels of Möbius type} at hand,  we  give  necessary and sufficient conditions for two irreducible representations of a group $G$ on a  real or complex hyperbolic space to be equivalent. In  \cite{Kimmarkedlenght} Inkang Kim, developing on techniques of Jean-Pierre Otal (see \cite{Otal}), showed that in the locally compact case the \textit{displacement} of an irreducible  representation is a complete invariant, which is not the case in the non-locally compact case (see \cite{monod2018notes,stolowicz2022complex}).
	Using the \textit{kernels of Möbius type} it  is straightforward to generalize that result  for any hyperbolic representation, possibly on an infinite-dimensional hyperbolic space. 
	
	\subsection{Kernels of negative type, GNS constructions and deformations}\label{kernelesnegativos}
	If $X$ is a set (topological space), a (continuous)  map ${\Psi}:X^2\rightarrow\F$, where $\F=\R,\C$, is called an $\F$-\textit{kernel of  negative type} if the following hold. 
	\begin{enumerate}
		\item For every $c_1,\dots, c_n\in \F$ such that $\sum_{i}c_i=0$ and every $x_1,\dots,x_n\in X$, $\sum_{i,j}c_i\overline{c_j}\Psi(x_i,x_j)\leq 0.$
		\item For every $x,y\in X$, $ \Psi(x,y)=\overline{\Psi(y,x)}.$ 
		\item For every $x\in X$, $\Psi(x,x)=0.$
	\end{enumerate}
	A map ${\Phi}:X^2\rightarrow\F$, where $\F=\R,\C$, is  called an $\F$-\textit{kernel of positive  type} if  	for every $c_1,\dots, c_n\in \F$ and every  $x_1,\dots,x_n\in X$, $\sum_{i,j}c_i\overline{c_j}\Psi(x_i,x_j)\geq 0.$
		The following lemma relates both definitions and can be found in Lemma 3.2.1 in \cite{harmonicemigroups}.
	\begin{prop}\label{positivenegative}
		If  $X^2\xrightarrow{\Psi}\F$  is a (continuous) map, then $\Psi$ is a $\F$-kernel of negative type if, and only if, for any (every) $x_0\in X$, the map $X^2\rightarrow\F$ given by 
		\[(x,y)\mapsto \frac{1}{2}\Big(\Psi(x,x_0)+\Psi(x_0,y)-\Psi(x,y)\Big),\]
		is a $\F$-kernel of positive type. 
	\end{prop}
	The \textit{GNS constructions} will show that the maps considered in the following  lemmas are  {generic}. The first lemma is classical and follows from the definitions.  
	\begin{lem}\label{funcionesykernels}
		If $X$ is a set (topological space), $(H,\langle\,\meddot\,,\meddot\,\rangle)$ is an  $\F$-Hilbert space and $X\xrightarrow{f}\F$ is a (continuous) map, then the following hold, 
		\begin{enumerate}\label{exemplos1}
			\item The map given by $(x,y)\mapsto \langle f(x),f(y)\rangle$ is a (continuous) $\F$-kernel of positive type.
			\item If $\F=\R$, the map given by $(x,y)\mapsto \|f(x)-f(y)\|^2$ is a (continuous) real kernel of negative type. 
		\end{enumerate}
	\end{lem}
	\begin{lem}\label{ejemplos2}If $X$ is a set (topological space) and $X\xrightarrow{f}\mathbf{H}_m$ is a (continuous) map, then  the map  $(x,y)\mapsto B_k(f(x),f(y))$ is a (continuous) complex kernel of negative type. 
	\end{lem}
	\begin{proof}
		Let $c_1,\dots,c_n\in\C$ be such that
		$\sum_{i}c_i=0$ and let $x_1,\dots,x_n\in X$. 
		If $m=0$, recall that $\He_0=i\R$ 
		and $B_k(\xi(b),\xi(d))=i(b-d)$ (see \cref{heisengroup}), thus in this case  the claim  holds. 
			If $m>0$, denote $f(x_i)=(u_i,b_i)$ and notice  that by \cref{funcionesykernels}, 
		\begin{equation*}\begin{array}{rcl}
			\sum\limits_{i,j}c_i\overline{c_j}B_k(f(x_i),f(x_j))&=&\\
			\sum\limits_{i,j}c_i\overline{c_j}\left(\|u_i\|^2+\|u_j\|^2\right)-
			2\sum\limits_{i,j}c_i\overline{c_j}\Big( \mathrm{Re}(\langle u_i,u_j\rangle) + i \mathrm{Im}(\langle u_i,u_j\rangle)    \Big)&=&\\
			-2\sum\limits_{i,j}c_i\overline{c_j}\langle u_i,u_j\rangle&\leq&0.
		\end{array}\fintres\end{equation*}\renewcommand{\qedsymbol}{}
	\end{proof}
	\vspace{-\baselineskip}
	The $\F$-kernels of positive or negative type associated to subsets of $\F$-Hilbert spaces or of  Heisenberg groups are called \textit{tautological} kernels. The constructions in the following theorems, known as the \textit{GNS constructions},   show that the tautological kernels are generic. 
	
	A proof for the following theorem can be found in Theorem C.2.3 in \cite{propertyTbekka}. 
	\begin{teo}
		If $X$ is a topological space and $X^2\xrightarrow{\Psi}\R$ is a (continuous) real kernel of negative type, then for any (every) $x_0\in X$, there exists a real Hilbert space $H$ and a (continuous) map 
		$X\xrightarrow{f}H$ such that \begin{enumerate}
			\item $\Psi(x,y)=\|f(x)-f(y)\|^2$. 
			\item The linear spam of $\{f(x)-f(x_0)\}_{x\in X}$ is dense. 
		\end{enumerate} Moreover, the space $H$ and the map $f$ are   initial with this property. That is to say, if there exist $L$ a real Hilbert space and a (continuous) map  $X\xrightarrow{g}L$ satisfying 1) and 2), then there exists a unique affine isometric map $H\xrightarrow{T}L$ such that for every $x\in X$, $T(f(x))=g(x).$
	\end{teo}
	 Below we sketch  a proof of the following theorem  because the arguments in it will be used later in this work (for a proof see Theorem C.1.4 in \cite{propertyTbekka}.)
	\begin{teo}\label{GNSpos}
		If $X$ is a topological space and $X^2\xrightarrow{\Phi}\F$ is a (continuous) $\F$-kernel of positive type, then there exists $H$ an $\F$-Hilbert space and a (continuous) map $X\xrightarrow{f}H$ such that 
		\begin{enumerate}
			\item For every $x,y\in X$, $\langle f(x),f(y)\rangle=\Phi(x,y)$.
			\item The linear spam of $\{f(x)\}_{x\in X}$ is dense. 
		\end{enumerate}
		Moreover, the space  $H$ and the map $f$ are initial with that property. That is to say, if there exist $L$ an $\F$-Hilbert space and a (continuous) map $X\xrightarrow{g}L$ that satisfy 1) and 2), then there exists a unique unitary  $\F$-linear map $H\xrightarrow{T}L$ such that for every $x\in X$, $T(f(x))=g(x).$
	\end{teo}
	\begin{proof}
		For every $x\in X$, Denote $\Phi_x$ the map given by $y\mapsto \Phi(x,y)$ and denote $N$  the $\F$-vector space generated by $\{\Phi_x\}_{x\in X}$ inside the space of continuous functions $\mathcal{C}(X,\F)$. Define  	$\langle\Phi_x,\Phi_y\rangle_\Phi=\Phi(x,y)$ and observe that if $h=\sum\limits_i\lambda_i \Phi_{x_i}=0$, then for every $y$, 
		\[\begin{array}{rcl} \sum\limits_i\lambda_i \langle \Phi_{x_i},\Phi_{y}\rangle_\Phi&=& h(y)\\&=&0\\&=&\overline{h(y)}\\&=&
			\sum\limits_i\overline{\lambda_i} \langle \Phi_y,\Phi_{x_i}\rangle_\Phi.		
		\end{array}\]
		The previous computation shows that $\langle\,\meddot\,,\meddot\,\rangle_\Phi$ can be extended by linearity to a Hermitian form defined on  $N$ and that for every $y\in X,$ $\langle h, \Phi_y\rangle_\Phi=h(y)$. The fact that the Hermitian product $\langle\,\meddot\,,\meddot\,\rangle_\Phi$ is positive definite  is a consequence of $\Phi$ being of positive type. The space $L$ is the completion of $N$ with respect to $\langle\,\meddot\,,\meddot\,\rangle_\Phi$ and the map $f$ is given by $f(x)= \Phi_x.$
	\end{proof} 
	
	The  GNS construction for complex kernels of negative type  introduced here uses ideas that can be found  in  \textit{e.g.}  Proposition 3.2 in Chapter 3 in \cite{harmonicemigroups}.
	\begin{teo}\label{gnsnegativocomplejo}
		If $X$ is a topological space, $x_0\in X$ and $X^2 \xrightarrow{\Psi}\C$ is a (continuous) complex kernel of negative type, then there exist a Heisenberg group $\He_m$ and  a (continuous) map $X\xrightarrow{f}\He_m$, such that $f(x_0)=e$ and  such that for every $x,y\in X$, 
		$\Psi(x,y)=B_k(f(x),f(y)).$
		Moreover, $(\He_m,f)$ is initial with that property. That is to say, if $X\xrightarrow{g}\He_n$ is another map such that $g(x_0)=e$ and such that  for every $x,y\in X$, 
		$\Psi(x,y)=B_k(g(x),g(y)),$ then there exists a unique map $\He_m\xrightarrow{T}\He_n$ such that
		for every $p,q\in\He_m$, $B_k\big(T(p),T(q)\big)=B_k\big(p,q\big)$ and such that $T\circ f=g.$
		
	\end{teo}		
	\begin{proof}
		If $x_0\in X$,	by \cref{positivenegative} the map $\Phi$, given by \[(x,y)\mapsto \frac{1}{2}\Big(\Psi(x,x_0)+\Psi(x_0,y)-\Psi(x,y)\Big),\] is a complex kernel of positive type. By \cref{GNSpos} there exist $(N, \langle\,\meddot,\,\meddot\,\rangle)$ a complex Hilbert space and a (continuous) map $i:X\rightarrow N$ such that for every $x,y\in X$, $\Phi(x,y)=\langle i(x),i(y)\rangle.$ By definition, for every $x\in X$, $\Phi(x,x_0)=0,$ and therefore,  $i(x_0)=0$ (see \cref{GNSpos}).
		
		Let $\He_m$ be the Heisenberg group associated to $N\times\R$. Denote $j:X\rightarrow\R$ the function  given by $j(x)= \mathrm{Im}(\Psi(x,x_0))$ and denote   $f:X\rightarrow\He_m$ the map  given by 
		$x\mapsto(i(x),j(x))$. It is immediate that $f(x_0)=e$ and observe that  for every $x,y\in X$, 
		$\Psi(x,y)=B_k(f(x),f(y))$.  Indeed, 
		by definition \begin{equation}
			\langle i(x)-i(y),i(x)-i(y)\rangle=
			\mathrm{Re}(\Psi(x,y)),
		\end{equation} 
		and 	therefore,  
		\begin{equation}\label{ecu3}\begin{array}{rcl}B_k(f(x),f(y))&=&\\
				\langle h(x)-h(y),h(x)-h(y)\rangle +i\Big(j(x)-j(y) -2\mathrm{Im}\big(\langle h(x),h(y)\rangle\big)\Big)&=&\\
				\mathrm{Re}(\Psi(x,y))+i\Big(\mathrm{Im}(\Psi(x,x_0))-\mathrm{Im}(\Psi(y,x_0))-\mathrm{Im}(\Psi(x,x_0))-\mathrm{Im}(\Psi(x_0,y))+\mathrm{Im}(\Psi(x,y))\Big)&=&\\
				\Psi(x,y).	
		\end{array}\end{equation}	
		
		Suppose  there exist $\He_{m'}$ a Heisenberg group and  a map $X\xrightarrow{g}\He_{m'}$ such that $g(x_0)=e$ and such that for every $x,y\in X$, $\Psi(x,y)=B_k(g(x),g(y)).$
				 Let $\Psi'$ be the tautological complex kernel of negative type associated to $\He_{m'}$ and consider  
		$\Phi'$, the complex kernel of positive type associated to $\Psi'$ with respect to $g(x_0)$ (see \cref{positivenegative}). As for every
		$x,y\in X$, $\Psi'(g(x),g(y))=\Psi(x,y)$, then 
		$\Phi'(g(x),g(y))=\Phi(x,y).$
		
		If $\He_{m'}$ is identified with $K\times \R$, where $(K,\langle\,\meddot\,,\meddot\,\rangle)$ is a complex Hilbert space and for every $x\in X$,   $g(x)=(u_x,e_x)$, then for every $x,y\in X$,  
		\[\begin{array}{rcl}
			\mathrm{Re}(\Phi'(g(x),g(y)))&=&
			\frac{1}{2}\mathrm{Re}\Big(\Psi'(g(x),g(x_0))+\Psi'(g(x_0),g(y))-\Psi'(g(x),g(y)) \Big)	\\&=&
			\frac{1}{2}\Big(\|u_x\|^2+\|u_y\|^2-\|u_x-u_{y}\|^2\Big)\\&=&
			\mathrm{Re}(\langle u_x,u_y\rangle)
		\end{array}
		\]	
		and
		\[\begin{array}{rcl}
			\mathrm{Im}\big(\Phi'(g(x),g(y))\big)&=&
			\frac{1}{2}\mathrm{Im}\Big(\Psi'(g(x),g(x_0))+\Psi'(g(x_0),g(y))-\Psi'(g(x),g(y)) \Big)\\	&=&
			\frac{1}{2}\Big(e_x-e_y-\big(e_x-e_y-2\mathrm{Im}(\langle w_x,w_y\rangle)\big)\Big)\\&=&
			\mathrm{Im}(\langle w_x,w_{y}\rangle).
		\end{array}\]
		Thus for every $x,y\in X$,
		\[\langle h(x),h(y)\rangle=\Phi(x,y)
			=\Phi'(g(x),g(y))=\langle w_x,w_{y}\rangle,\]
		and  by \cref{GNSpos}, there exists a unique unitary  map $T:L\rightarrow K$ such that for every $x\in X,$ $T(h(x))=u_x$.
		Moreover, 
		\[
			j(x)=\mathrm{Im}\left(\Psi(x,x_0)\big) \right)
			=\mathrm{Im}\left(\Psi'(g(x),g(x_0))\right)=e_x,
		\]
		showing that $(T\times \mathrm{Id})\circ f=g.$
		
		Now suppose that there exists a map $S:\He_m\rightarrow\He_n$ such that for every $p,q\in\He_m$, $B_k(p,q)=B_k(S(p),S(q))$ and such that for every $x\in X$, $S(f(x))=g(x).$ The claim is that $S=T\times \mathrm{Id}.$ By \cref{linearrep}, $S=Q\times \mathrm{Id}$ for some unitary map $Q$. By the uniqueness (coming from the GNS construction) of the map $T$, then $Q=T.$
	\end{proof}
	The following theorems, that  will allow us to \say{deform} representations,  are  Corollaries  2.10 and  2.11 in Chapter 2 of  \cite{harmonicemigroups}.
	\begin{teo}\label{powert}
		If $X^2 \xrightarrow{\Psi}\C$ is a kernel of complex negative type and $t\in(0,1)$, then $\Psi^t$ is a kernel of complex  negative type. 
	\end{teo}
	\begin{teo}\label{power2}
		If $X\xrightarrow{f}\C$ is a map such that $\mathrm{Re}(f)\geq 0$ and $t\in [1,2)$,
		then the map $X^2 \rightarrow\C$ given by
		\[(x,y)\mapsto - (f(x)+\overline{f(y)})^t\] is a complex kernel of negative type. 
	\end{teo}
	
	For $z=re^{i\theta}$, with  $\theta\in (-{\pi},\pi)$,   $t\in\R$ and $s\in(-1,1)$, define $z^{t,s}=|z|^te^{is\theta}$ and  $0^{t,s}=0$.  For $z=ia\in\C$, with $a>0$ and $t\in(1,2)$, \begin{equation}\label{tmayor2}-z^t=|z|^te^{\tfrac{i(t-2)\pi}{2}}={z}^{t,t-2}.\end{equation}
	If $z=-ia$  the same equation holds. 
	
	A (continuous) function $\varphi:\R\rightarrow\C$ is called a (continuous) \textit{complex function of  negative type} if the map $\R^2\rightarrow\C$ given by $(x,y)\mapsto \varphi(x-y)$ is a (continuous) complex kernel of negative type.
	For  further reading on functions of negative type see \textit{e.g.} \cite{harmonicemigroups,bernstein,potential}.   
	
	A function $\varphi:\R\rightarrow\C$  is a function of complex negative type if, and only if, it admits a \textit{Lévy-Khintchine representation}. That is to say,  there exist a Radon measure $\nu$ on $\R\setminus\{0\}$ such that \[\int\limits_{\R\setminus\{0\}} \min(1,y^2)d\nu(y)<\infty,\] $\alpha,\gamma\geq0$ and  $\beta\in\R$  such that 
	\[\varphi(x)= i\beta x+\tfrac{1}{2}\gamma x^2+\int\limits_{\R\setminus\{0\}}\left(1-e^{ixy}+\tfrac{xy}{1+|y|^2}\right)d\nu(y),\]
	where the triple $(\beta,\gamma,\nu)$ is unique (see Theorem 4.15 in \cite{bernstein}). 
	
	The following theorem is  crucial  for the classification of irreducible complex hyperbolic representations of $\mathrm{PU}(1,1)$ (see \cref{horosphcomb}). I would like to  express my gratitude to Yusen Long and David Xu for pointing out an error in both the statement and the proof of an earlier version of this theorem and also to Yusen Long for suggesting the correct statement.   
	\begin{teo}\label{productotorcido}
		Let $\varphi:\R\rightarrow\C$ be a continuous  function of complex  negative type such that there exists $t\in(0,1)\cup(1,2)$ such that for every $\lambda>0$ and every $x\in\R$, $\varphi(\lambda x)=\lambda^t\varphi(x)$. Then if $t\in(0,1)$ (resp. $t\in(1,2)$),     for every  $x\in \R\setminus\{0\}$, $|\Arg(\varphi(x))|\leq \tfrac{t\pi}{2}$  (resp. $|\Arg(\varphi(x))|\leq \tfrac{(2-t)\pi}{2}$). 
	\end{teo}
	\begin{proof}
The function  $\varphi$ admits a Lévy-Khintchine representation associated to a triple $(\beta,\gamma,\nu)$. Let $\lambda>0$ and observe that  on the one hand
		for every $x\in\R$, 
		$\tfrac{\varphi(\lambda x)}{\lambda^t}=\varphi(x)$, and on the other hand,  
		\begin{equation}\label{Levyeq}\tfrac{\varphi(\lambda x)}{\lambda^t}=i\lambda^{1-t}\beta x+\tfrac{1}{2}\lambda^{2-t}\gamma x^2+\lambda^{-t}\int\limits_{\R\setminus\{0\}}\left(1-e^{i\lambda xy}+\tfrac{i\lambda xy}{1+y^2}\right)d\nu(y).\end{equation}
		Denote $\nu_\lambda$  the measure given by $\nu_\lambda(A)=\nu(\lambda^{-1}A)$ and observe that $\int_{\R\setminus\{0\}}\min(1,y^2)d\nu_\lambda<\infty$.
		In order to use the uniqueness of the triple $(\beta,\gamma,\nu)$ associated to $\varphi$ by the Lévy-Khintchine representation, one would like to write \[ \tfrac{\lambda xy}{1+y^2}=xf_\lambda(y) + \tfrac{\lambda xy}{1+\lambda^2y^2},\]
		for some $f_\lambda:\R\setminus\{0\}\rightarrow\R$ such that $\left\lvert\int_{\R\setminus\{0\}}f_\lambda(y)d\nu(y)\right\lvert<\infty.$
		Indeed,  $f_\lambda(y)=\tfrac{\lambda(1-\lambda^2) y^3}{(1+y^2)(1+\lambda^2y^2)}$ and observe that $\tfrac{ |y|^3}{(1+y^2)(1+\lambda^2y^2)}\leq \min(1,y^2)$.
		 Let 
		$R_\lambda=\int_{\R\setminus\{0\}}f_\lambda(y)d\nu(y)\in\R$  and  rewrite (\ref{Levyeq}) as 
		\[\begin{array}{rcl}\tfrac{\varphi(\lambda x)}{\lambda^t}&=&ix(\lambda^{1-t}\beta+\lambda^{-t}R_\lambda) x+\tfrac{1}{2}\lambda^{2-t}\gamma x^2+\lambda^{-t}\int\limits_{\R\setminus\{0\}}\left(1-e^{i\lambda xy}+\tfrac{i\lambda xy}{1+\lambda^2y^2}\right)d\nu(y)\\
			&=&
			ix(\lambda^{1-t}\beta+\lambda^{-t}R_\lambda) +\tfrac{1}{2}\lambda^{2-t}\gamma x^2+\lambda^{-t}\int\limits_{\R\setminus\{0\}}\left(1-e^{i xy}+\tfrac{i xy}{1+y^2}\right)d\nu_\lambda(y).
		\end{array}\]
		Thus for every $\lambda>0$,
		\[(\beta,\gamma,\nu)=(\lambda^{-t}(\lambda\beta+R_\lambda),\lambda^{2-t}\gamma,\lambda^{-t}\nu_\lambda),\] 
		and therefore,  $\gamma=0$ and   for every $\lambda>0$, $\lambda^{-t}\nu_\lambda=\nu$. 
			Up to a positive scalar the Haar measure on $\R_{>0}$    is $\tfrac{dx}{x}$, thus one can deduce that there exist $a,b\geq0$ such that 
		$\nu$ restricted to $\R_{<0}$ (resp. $\R_{>0}$) is equal to $a\tfrac{dy}{|y|^{t+1}}$ (resp. $b\tfrac{dy}{y^{t+1}}$).  
		Hence the Lévy-Khintchine representation of $\varphi$ is given by, 
		\begin{equation}\label{Levy}\varphi(x)=i\beta x + a\int\limits_{\R_{<0}}\left(1-e^{iyx}+\tfrac{ixy}{1+y^2}\right)\tfrac{dy}{|y|^{t+1}}+b\int\limits_{\R_{>0}}\left(1-e^{iyx}+\tfrac{ixy}{1+y^2}\right)\tfrac{dy}{y^{t+1}}.\end{equation}
		When $t<1$, by \textit{e.g.} Example 6.2.14 in \cite{Whittaker}, $S=\int_{\R_{>0}}\tfrac{1}{(1+y^2)y^t}dy<\infty$ and therefore, (\ref{Levy}) can be written as  
		\begin{equation}\label{Levy2}\varphi(x)=ix(\beta +(b-a)S) + a\int\limits_{\R_{<0}}\left(1-e^{iyx}\right)\tfrac{dy}{|y|^{t+1}}+b\int\limits_{\R_{>0}}\left(1-e^{iyx}\right)\tfrac{dy}{y^{t+1}}.\end{equation}
		For $z\in\C$ with $\mathrm{Re}(z)\geq 0$ and for  $0<t<1$, from the integral representation (see Corollary 2.10 in Chapter 2 of \cite{harmonicemigroups})  
		\begin{equation*}\label{integralpotencia}z^t=\tfrac{t}{\Gamma(1-t)}\int\limits_{\R_{>0}}(1-e^{-zy})\tfrac{dy}{y^{t+1}}, \end{equation*}
		one gets that, 
		\[(ix)^t=|x|^t\cos(\tfrac{t\pi}{2})+i|x|^t\mathrm{sign}(x)\sin(\tfrac{t\pi}{2})=\tfrac{t}{\Gamma(1-t)}\int\limits_{\R_{>0}}(1-\cos(xy))\tfrac{dy}{y^{t+1}}-i\tfrac{t}{\Gamma(1-t)}\int\limits_{\R_{>0}}\sin(xy)\tfrac{dy}{y^{t+1}}.\]
		This implies that (\ref{Levy2}) can be written as 
		\[
		\varphi(x)=ix(\beta +(b-a)S)+(a+b)\tfrac{\Gamma(1-t)}{t}|x|^t\cos(\tfrac{t\pi}{2})+i(b-a)\tfrac{\Gamma(1-t)}{t}\mathrm{sign}(x)|x|^t\sin(\tfrac{t\pi}{2})
		.\]
		As $\varphi$ is homogenous of degree $t$, $\beta +(b-a)S=0$, and therefore,  for every $x\in\R\setminus\{0\}$, 
		\[\tan(|\Arg(\varphi(x))|)=\tfrac{|b-a|}{a+b}\tan(\tfrac{t\pi}{2})\leq \tan(\tfrac{t\pi}{2}). \]
		Suppose now that $t=1+t'$, with $0<t'<1$, and consider for $z\in\C$ such that $\mathrm{Re}(z)\geq0$ the integral expression (see Corollary 2.22 in Chapter 2 in \cite{harmonicemigroups}) 
		\begin{equation*}
			-z^{1+t'}=\tfrac{t'(1+t')}{\Gamma(1-t')}\int\limits_{\R_{>0}}(1-e^{-y z}-y z)\tfrac{dy}{y^{2+t'}}.
		\end{equation*}	
		With a change of variable one gets that 
		\[	-z^{1+t'}=c\int\limits_{\R_{<0}}(1-e^{y z}+y z)\tfrac{dy}{|y|^{2+t'}}, \]
		for some constant $c>0$.  Declare $f(y)=\tfrac{y^3}{1+y^2}$ and notice that 
		\[	-z^{1+t'}=c\int\limits_{\R_{<0}}(1-e^{y z}+\tfrac{zy}{1+y^2} +zf(y))\tfrac{dy}{|y|^{2+t'}}. \]
		By \textit{e.g.} Example 6.2.14 in \cite{Whittaker}, 
		$S=c\int\limits_{\R_{<0}}f(y)\tfrac{dy}{|y|^{2+t'}}\in\R,$ and therefore, 
		for every $x\in\R$, 
		\[	-(ix)^{1+t'}=c\int\limits_{\R_{<0}}(1-e^{i xy}+\tfrac{ixy}{1+y^2} +ixf(y))\tfrac{dy}{|y|^{2+t'}}=iSx+c\int\limits_{\R_{<0}}(1-e^{i xy}+\tfrac{ixy}{1+y^2} )\tfrac{dy}{|y|^{2+t'}}. \]
		Using this identity one can rewrite (\ref{Levy})
		as 
		\[\varphi(x)=ix\big(\beta +\tfrac{S}{c}(b-a)\big)-(a+b)\mathrm{Re}((ix)^{1+t'})-(a-b)\mathrm{Im}((ix)^{1+t'}).\]
		As  $\varphi$ is homogeneous of degree $t$, $\beta +\tfrac{S}{c}(b-a)=0$,  and therefore for every $x\in R\setminus\{0\}$,
		\[\tan(|\Arg(\varphi(x))|)=\left|\tfrac{|a-b|\mathrm{Im}((ix)^{t})}{(a+b)\mathrm{Re}((ix)^{t})}
		\right|)\leq\tan(|\Arg(-(ix)^t)|=\tfrac{(2-t)\pi}{2}.\qedhere\]
	\end{proof}
	
	The author believes that this approach could lead  to a classification for the irreducible complex hyperbolic representations of $\mathrm{PU}(1,n)$ with  $n>1$. 
	
	Given a set $X$ and an $\F$-kernel of negative type ${\Psi}:X^2\rightarrow \F$, for $t\in(0,1)$ (resp.  $t\in(1,2)$ if $\Psi$ is as in \cref{power2}),   one can ask weather the Hilbert space or the Heisenberg space associated to $\Psi^t$ by the GNS construction   is infinite-dimensional.  This question and its implications regarding the existence of infinite-dimensional representations for a given group will be addressed in \cref{restrictionofrepresentations}.

	\subsection{Kernels of Möbius type and the  GNS construction}\label{seccionkernelesmobius} 
	
	The cross-ratio of a metric space holds relevant  information, up to a Möbius a transformation, for quadruples of points, but for our purposes with regard to complex hyperbolic spaces this is not sufficient. It is  necessary to consider a map that relates to  the \textit{Korányi-Riemann cross-ratio}, that for   four  distinct points  $x_1,\dots,x_4$ in the boundary of a complex hyperbolic space is defined as \[\tfrac{B(\tilde{x}_1,\tilde{x}_3)B(\tilde{x}_2,\tilde{x}_4)}{B(\tilde{x}_2,\tilde{x}_3)B(\tilde{x}_1,\tilde{x}_4)},\]
	where $\tilde{x}_i$ is any lift of $x_i$ (see 7.2 in \cite{complexhyperbolic}). With this in mind, it is natural to believe that for complex representations the adequate analogous for a cross-ratio would be a map $X^{(4)}\rightarrow \mathbf{P}(\C^3)$ that has   properties analogous to those  of the Korányi-Riemann cross-ratio.  
	
	If a set $X$ is provided with a topology, it will be supposed Hausdorff and 1\textsuperscript{st} countable.	
	Let $\mathrm{D}_8<\mathrm{Sym}(4)$ be the dihedral group  whose elements are the permutations $e$,  $(13)$, $ (24)$, $(12)(34)$, $(13)(24)$, $(14)(23)$, $(1234)$ and $(1432)$. Denote $h=(13)$, $i=(13)(24)$,  $j=(14)(23)$ and $k=(12)(34)$. Observe that $h$ and $j$  can be chosen as generators with the relations 
	$(jh)^4=h^2=e$ and $hj=h^{-1}j^{-1}$. 
	
	Denote \[\Sigma=\{[a:b:c]\in \mathbf{P(F^3)}\mid \{a,b,c\}\subset\F^*\}\]
	and let  ${F_h,F_i,F_j},F_k:\Sigma\rightarrow\Sigma$ be the maps given by,   
	\begin{enumerate}
		\item $F_h\left(\left[a:b:c\right]\right)=\left[\overline{c}:\overline{b}:\overline{a}\right]$.
		\item $F_i\left(\left[a:b:c\right]\right)=\left[\frac{1}{\overline{c}}:\frac{\overline{b}}{\overline{a}\overline{c}}:\frac{1}{\overline{a}}\right]$.
		\item $F_j\left(\left[a:b:c\right]\right)=\left[\frac{\overline{a}}{\overline{b}c}:\frac{1}{c}:\frac{1}{b}\right]$.
		\item $F_k\left(\left[a:b:c\right]\right)=\left[\frac{a\overline{a}}{b}:\overline{a}:\frac{\overline{ac}}{\overline{b}}\right]$.
	\end{enumerate}	
	The notation suggests an isomorphism between the  group generated by $\{F_h,F_i,F_j,F_k\}$ and $\mathrm{D}_8$. 
	
	Recall that for a set $X$,  $X^{\{4\}}\subset X^4$ is the set of 4-tuples  of pairwise distinct elements.  Denote  $X^{(4)}$ the set of $(x_1,x_2,x_3,x_4)\in X^4$ such that for every $i\neq 4$, $x_i\neq x_4$ and such that for $i=1,2,3$, $x_i$ appears at most twice. 
	Observe that $\mathrm{D}_8$ has a natural action in $X^{\{4\}}$.
	
	An $\F$-kernel of negative type $\Psi:X^2\rightarrow\F$  is called \textit{injective} if for every two distinct $x,y\in X$, $\Psi(x,y)\neq0.$ Given a  (continuous) map ${\Gamma}:X^{(4)}\rightarrow\mathbf{P}(\F^3)$,  denote for $x\in X^{(4)}$, \[\Gamma(x)=\left[\Gamma_1(x):\Gamma_2(x):\Gamma_3(x)\right]\] and for  $x\in X^{(4)}$ such that   $\Gamma(x)\in \Sigma$, define $\Gamma_{i,j}(x)=\tfrac{\Gamma_i(x)}{\Gamma_j(x)}$. 
	A (continuous) map $\Gamma:X^{(4)}\rightarrow\mathbf{P}(\F^3)$  is called an $\mathbf{F}$-\textit{kernel of Möbius type} if it has the following properties.	
	\begin{enumerate}[label=A\arabic*.]
		\item If $x\in X^{\{4\}}$, then  
		$\Gamma(x)\in\Sigma$.
		\item 
		The restriction $\Gamma:X^{\{4\}}\rightarrow\Sigma $ is $\mathrm{D}_8$-equivariant. 
		
		
		\item For every $d$ and every two distinct elements $a,b$ of $X\setminus \{d\}$, 
		\begin{enumerate}
			\item $\Gamma(a,b,b,d)=[0:1:1]$.
			\item $\Gamma(b,a,b,d)=[z:0:\overline{z}]$, for some $z\in\F^*$.
			\item $\Gamma(b,b,a,d)=[1:1:0]$.
		\end{enumerate}  
		\item For every $d\in X$ and    $(x_1,\dots,x_4)\in(X\setminus \{d\})^{\{4\}}$, 
		\begin{enumerate}
			\item	 $\Gamma_{2,1}(x_2,x_1,x_3,d)\Gamma_{2,1}(x_1,x_2,x_4,d)=\Gamma_{2,3}(x_1,x_3,x_4,d)\Gamma_{3,2}(x_2,x_3,x_4,d).$ 
			\item $\Gamma_{2,1}(x_1,x_3,x_4,d)=\Gamma_{2,1}(x_2,x_3,x_4,d)\Gamma_{2,1}(x_1,x_2,x_4,d).$
			\item $\Gamma_{3,2}(x_1,x_2,x_4,d)=\Gamma_{3,2}(x_1,x_3,x_4,d)\Gamma_{2,3}(x_1,x_3,x_2,d).$
		\end{enumerate}
		\item If $x_1,\cdots,x_5$ are five distinct elements of $X$, 
		\begin{enumerate}\item 	$\Gamma_{2,1}(x_1,x_2,x_3,x_4)\Gamma_{1,2}(x_1,x_2,x_5,x_4) =
			\Gamma_{2,1}(x_1,x_2,x_3,x_5)  $.
			\item $\Gamma_{2,3}(x_1,x_2,x_3,x_4)\Gamma_{3,2}(x_5,x_2,x_3,x_4)=\Gamma_{2,3}(x_1,x_2,x_3,x_5).$
		\end{enumerate}	
		\item 		For every $d\in X$ and every  two distinct $x_0,x_1\in X\setminus \{d\}$, the map ${\Psi_{x_0,x_1}'}:(X\setminus \{d\})^2\rightarrow\F$,  
		given by  $$
		\Psi_{x_0,x_1}'(a,b)= \begin{cases}\Gamma_{3,2}(a,b,x_0,d)\Gamma_{2,1}(a,x_1,x_0,d)&\text{if}\,\, a\neq x_0\\ 
			\Gamma_{1,3}({x}_1,{x}_0,b,{d})& \text{if}\,\,  a=x_0,\end{cases}$$
		is such that there exists $z_0\in \F^*$ such that $\Psi_{x_0,x_1}=z_0\Psi_{x_0,x_1}'$ is  an injective $\F$-kernel of negative type.\end{enumerate}
	The definition of a  real kernel of Möbius type can be simplified in the following sense. Let  $F_{(12)}$, $F_{(23)}$, $F_{(34)}$ and  $F_{(41)}$ be the maps $\Sigma\rightarrow\Sigma$ given by,   
	\begin{enumerate}
		\item $F_{(12)}\left(\left[a:b:c\right]\right)=\left[b:a:c\right]$.
		\item $F_{(23)}\left(\left[a:b:c\right]\right)=\left[a:c:b\right]$.
		\item $F_{(34)}\left(\left[a:b:c\right]\right)=\left[b:{a}:{c}\right]$.
		\item $F_{(41)}\left(\left[a:b:c\right]\right)=\left[\frac{a}{b}:1:\frac{b}{c}\right]$.
	\end{enumerate}	
	
	The axioms A1, A4, A5 remain the same and the axiom A2 can be stated changing $\mathrm{D}_8$-equivariant for $S_4$-equivariant. The axiom A3 can be stated as follows. For every $d$ and every two distinct $a,b\in X\setminus \{d\}$ the following hold. 
	\begin{enumerate}
		\item $\Gamma(a,b,b,d)=[0:1:1]$.
		\item  $\Gamma(b,a,b,d)=[1:0:1]$.
		\item $\Gamma(b,b,a,d)=[1:1:0]$. 
	\end{enumerate}
	The axiom  A6 can be stated declaring that  $\Psi'_{x_0,x_1}$ is an injective kernel of real hyperbolic type. 
	

	\begin{prop}\label{tautologico}
		If $X\xhookrightarrow{\iota}\partial\mathbf{H}^m_\F$, then the map $X^{(4)}\xrightarrow{\Gamma_\iota}\mathbf{P}(\mathbf{F}^3)$ given by 
		\[\Gamma_\iota(x_1,x_2,x_3,x_4)=\left[\tfrac{B(\tilde{x}_2,\tilde{x}_3)}{B(\tilde{x}_2,\tilde{x}_4)B(\tilde{x}_4,\tilde{x}_3)}:\tfrac{B(\tilde{x}_1,\tilde{x}_3)}{B(\tilde{x}_1,\tilde{x}_4)B(\tilde{x}_4,\tilde{x}_3)}:\tfrac{B(\tilde{x}_1,\tilde{x}_2)}{B(\tilde{x}_1,\tilde{x}_4)B(\tilde{x}_4,\tilde{x}_2)}\right],\]
		where $\tilde{x}_i$ are respective lifts of $x_i$, is a (continuous) $\F$-kernel of Möbius type.  
	\end{prop}
	\begin{proof}
		The map $\Gamma$ is independent of the choice of the lifts of the elements of $X$ and the properties from A1 to  A5 hold trivially. If $x_0,x_1,d$ are three distinct elements of $X$, and $a,b\in X\setminus\{d\}$, where $a\neq x_0$, 
		then 
		\[\Gamma_{3,2}(a,b,x_0,d)\Gamma_{2,1}(a,x_1,x_0,d)=\tfrac{B(\tilde{a},\tilde{b})}{B(\tilde{a},\tilde{d})B(\tilde{d},\tilde{b})}\tfrac{B(\tilde{x}_1,\tilde{d})B(\tilde{d},\tilde{x}_0)}{B(\tilde{x}_1,\tilde{x}_0)}.\]
		When $a=x_0$,  \[\Gamma_{1,3}({x}_1,{x}_0,b,{d})=\tfrac{B(\tilde{x}_0,\tilde{b})}{B(\tilde{x}_0,\tilde{d})B(\tilde{d},\tilde{b})}\tfrac{B(\tilde{x}_1,\tilde{d})B(\tilde{d},\tilde{x}_0)}{B(\tilde{x}_1,\tilde{x}_0)}.\]  
		Hence  the map $\Psi'_{x_0,x_1}$ defined as in A6 is such that, for every $a,b\in X\setminus\{d\}$, 
		\[\Psi'_{x_0,x_1}(a,b)=\tfrac{B(\tilde{a},\tilde{b})}{B(\tilde{a},\tilde{d})B(\tilde{d},\tilde{b})}\tfrac{B(\tilde{x}_1,\tilde{d})B(\tilde{d},\tilde{x}_0)}{B(\tilde{x}_1,\tilde{x}_0)}.\]
		
		Given $\tilde{d}$ a lift of $d$, for every $a\in\partial\mathbf{H}_\C^m\setminus\{d\},$ there exists a lift $\tilde{a}$ such that, $B(\tilde{a},\tilde{d})=1,$ thus 
		\begin{equation}\label{metrica}\Psi'_{x_0,x_1}(a,b)=\tfrac{B(\tilde{a},\tilde{b})}{B(\tilde{x}_1,\tilde{x}_0)}.\end{equation}
		
		If $\F=\C$, after fixing $\tilde{x}_0$ and $\tilde{d}$, there is an identification  $f:\partial\mathbf{H}^m_\C\setminus \{d\}\rightarrow \He_{m-1}$ (see \cref{isometriakoran}) such that 
		$B_k(f(a),f(b))=B(\tilde{a},\tilde{b}).$ By \cref{ejemplos2}, the map 
		$B(\tilde{x}_1,\tilde{x}_0)\Psi'_{x_0,x_1}$ is an injective  complex kernel of negative type. 
			If $\F=\R$  with the same arguments    one can show that  $\Psi'_{x_0,x_1}$ is a real kernel of negative type. 
		
		Recall that if $\pi:L\setminus\{0\}\rightarrow\mathbf{P}(L)$ is the natural quotient, then the topology in $\partial\mathbf{H}^m_\F$ is given by the quotient topology induced by $\pi:L_0\rightarrow\partial\mathbf{H}^m_\F$, where $L_0$ is the set of isotropic vectors of $L$ (see the proof of \cref{cartancontinuous}).  This implies that $\pi:L_0\setminus\F \tilde{d}\rightarrow\partial\mathbf{H}^m_\F\setminus\{d\}$ is also a quotient. The map ${\tilde{\iota}}:\partial\mathbf{H}^m_\F\setminus\{d\}\rightarrow L_0\setminus\F\tilde{d}$ given by
		$a\mapsto \tilde{a}$ is such that for every $w\in L_0\setminus\F \tilde{d}$, 
		$\tilde{\iota}\circ\pi(w)=\tfrac{w}{B(w,\tilde{d})}.$
		This shows that $\tilde{\iota}$ is continuous, and therefore, that the map  $\Psi_{x_0,x_1}$  is  continuous. 
		
		The continuity of $\Gamma_{\iota}$ can be deduced in the same way. Indeed, the previous arguments  show that for every $d\in\partial\mathbf{H}^m_\F$, the map $\Gamma_{\iota}$ is continuous on the set 
		\[\{(x_1,x_2,x_3,x_4)\in X^{(4)}\mid x_4\neq d\}.\qedhere\]
	\end{proof}
	
	Observe that if $\Gamma$ is the $\F$-kernel of Möbius type associated to the inclusion $X\hookrightarrow \partial\mathbf{H}^m_\F,$
	then the map ${|\Gamma|^{\frac{1}{2}}}:X^{\{4\}}\rightarrow\mathbf{P}(\R^3)$ given by, \begin{equation}\label{gamavalorabsoluto}|\Gamma|^{\frac{1}{2}}(x)=\left[|\Gamma_{1}(x)|^{\frac{1}{2}}:|\Gamma_{2}(x)|^{\frac{1}{2}}:|\Gamma_{3}(x)|^{\frac{1}{2}}\right],\end{equation} is such that $|\Gamma|^{\frac{1}{2}}=crt$, where $crt$ is the cross-ratio associated to any Bourdon metric restricted to $X$ (see (\ref{eqcrossratio}) and  \cref{isometriakoran}).

	\begin{lem}\label{cuentasgamma}
		Let $X^{(4)}\xrightarrow{\Gamma}\mathbf{P}(\F^3)$ be  a (continuous) $\F$-kernel of Möbius type. If  $x_0,x_1,y_0\in X$ are pairwise distinct  and $(a,b,c,y_0)\in X^{\{4\}}$, where $a\neq x_0\neq b$, then the following hold. 
		\begin{enumerate}
			\item $\Gamma_{2,1}(a,x_1,x_0,y_0)	\Gamma_{1,2}(a,b,c,y_0)	\Gamma_{3,2}(a,c,x_0,y_0)=
			\Gamma_{2,1}(b,x_1,x_0,y_0)\Gamma_{3,2}(b,c,x_0,y_0).	$
			\item $\Gamma_{3,2}(a,b,c,y_0)	 \Gamma_{3,2}(a,c,x_0,y_0)=\Gamma_{3,2}(a,b,x_0,y_0).$
		\end{enumerate}
	\end{lem}
	\begin{proof}
		Observe that by A4, \[ 
		\Gamma_{1,2}(a,b,c,y_0)\Gamma_{3,2}(a,c,x_0,y_0)=\Gamma_{2,1}(b,a,x_0,y_0)\Gamma_{3,2}(b,c,x_0,y_0),\]
		and 
		\[\Gamma_{2,1}(b,x_1,x_0,y_0)=\Gamma_{2,1}(a,x_1,x_0,y_0)\Gamma_{2,1}(b,a,x_0,y_0).\]
	One can show 	using the two previous identities that 1. holds. Indeed,  
		\[\begin{array}{rcl}
			\Gamma_{2,1}(a,x_1,x_0,y_0)\Gamma_{1,2}(a,b,c,y_0)\Gamma_{3,2}(a,c,x_0,y_0)&=&\\
			\Gamma_{2,1}(a,x_1,x_0,y_0)	\Gamma_{2,1}(b,a,x_0,y_0)\Gamma_{3,2}(b,c,x_0,y_0)&=&\\
			\Gamma_{2,1}(b,x_1,x_0,y_0)\Gamma_{3,2}(b,c,x_0,y_0). \end{array}\]
		
		For 2. observe that 
		by A2, for every $u,v,w,y_0\in X$ pairwise distinct, 
		\[\Gamma_{3,2}(u,v,w,y_0)=\overline{\Gamma_{1,2}(w,v,u,y_0)}, \]
	and therefore, 
		\begin{equation*}\begin{array}{rcl}
			\Gamma_{3,2}(a,b,c,y_0)	 \Gamma_{3,2}(a,c,x_0,y_0)&=&
			\overline{\Gamma_{1,2}(c,b,a,y_0)}\overline{\Gamma_{1,2}(x_0,c,a,y_0)}\\
			&=&\overline{\Gamma_{1,2}(x_0,b,a,y_0)}\\
			&=&\Gamma_{3,2}(a,b,x_0,y_0).
		\end{array}\fintres\end{equation*} \renewcommand{\qedsymbol}{}
	\end{proof}
	\vspace{-\baselineskip}
	The following theorem is the GNS construction for the  kernels of Möbius type. 	
	\begin{teo}\label{GNSMobius}
		If $X^{(4)}\xrightarrow{\Gamma}\mathbf{P}(\F^3)$ is a (continuous) $\F$-kernel of Möbius type, then there exists $\mathbf{H}^m_\F$,   an injective (continuous) map $f:X\rightarrow\partial\mathbf{H}_\F^m$ such that  $\Gamma=\Gamma_f$, where $\Gamma_f$ is the $\F$-kernel of Möbius type induced by $f$,  and such that $(\mathbf{H}^m_\F,f)$ is  initial with those properties. That is to say, if there exists $\mathbf{H}^{m'}_\F$ and    an injective (continuous) map $g:X\rightarrow\partial\mathbf{H}_\F^{m'}$ such that  $\Gamma=\Gamma_g$, then there exists a unique map $T$, induced by a holomorphic  isometry $\mathbf{H}^m_\F\rightarrow\mathbf{H}^{m'}_\F$,  such that $T\circ f=g.$
	\end{teo}	
	\begin{proof}
		The proof  will be done for $\F=\C$; the case $\F=\R$ follows \textit{mutatis mutandis}. Fix  $y_0\in X$ and 
		choose  two distinct  $x_0,x_1\in X\setminus\{y_0\}$. Let ${\Psi_{x_0,x_1}}:(X\setminus\{y_0\})^2\rightarrow\C
		$ be the associated complex kernel of negative type (see A6). By  \cref{gnsnegativocomplejo},  there exists an injective (continuous) map ${\varphi}:X\setminus\{y_0\}\rightarrow\mathbf{H}_{m-1}$, for some $m\geq1$, such that $\varphi(x_0)=e$ and such that for every $a,b\in X\setminus \{y_0\}$,  
		$B_k(\varphi(a),\varphi(b))=\Psi_{x_0,x_1}(a,b).$
		
		Fix two distinct $q_1,q_2\in\partial\mathbf{H}^m_\C$ and respective lifts $\xi_1$ and $\xi_2$ such that $B(\xi_1,\xi_2)=1.$ After these choices there is an  identification $\omega:\mathbf{H}_{m-1}\rightarrow\partial\mathbf{H}^m_\C\setminus\{q_1\}$ (see \cref{isometriakoran}). Denote $f:X\rightarrow\partial\mathbf{H}^m_\C$ the map such that $f(y_0)=q_1$ and such that $f|_{X\setminus\{y_0\}}=\omega\circ\varphi$. 		
		By \cref{isometriakoran}, the (continuous) map $f|_{X\setminus\{y_0\}}$ has a lift $\tilde{f}$  such that $\tilde{f}(y_0)=\xi_1$  and   such that for every $a,b\in X\setminus\{y_0\}$, $B(\tilde{f}(a),\tilde{f}(y_0))=1$ and   
		\begin{equation}\label{kernelprima}\Psi_{x_0,x_1}(a,b)=B(\tilde{f}(a),\tilde{f}(b)).\end{equation}
		
		By A3, for every $x\in X^{(4)}\setminus X^{\{4\}}$, $\Gamma(x)=\Gamma_f(x)$.  	
		Let  $(a,b,c,y_0)\in X^{\{4\}}$ and  suppose that  $a\neq x_0\neq b$. By \cref{cuentasgamma}, 
		\[\begin{array}{rcl}	\Gamma(a,b,c,y_0)
			&=&[\Gamma_{1,2}(a,b,c,y_0):1:\Gamma_{3,2}(a,b,c,y_0)]\\
			&=&[\Gamma_{1,2}(a,b,c,y_0)\Gamma_{3,2}(a,c,x_0,y_0):\Gamma_{3,2}(a,c,x_0,y_0):\Gamma_{3,2}(a,b,c,y_0)\Gamma_{3,2}(a,c,x_0,y_0)]\\
			&=&[\Gamma_{1,2}(a,b,c,y_0)\Gamma_{3,2}(a,c,x_0,y_0):\Gamma_{3,2}(a,c,x_0,y_0):\Gamma_{3,2}(a,b,x_0,y_0)]\hfill\,\,(\text{Lemma}\,\,\ref{cuentasgamma})\\
			&=&[\Gamma_{1,2}(a,b,c,y_0)\Gamma_{3,2}(a,c,x_0,y_0)\Gamma_{2,1}(a,x_1,x_0,y_0):\Psi'_{x_0,x_1}(a,c):\Psi'_{x_0,x_1}(a,b)]\\
			&=&[ \Psi'_{x_0,x_1}(b,c) : \Psi'_{x_0,x_1}(a,c) :\Psi'_{x_0,x_1}(a,b)]\hfill(\text{Lemma}\,\,\ref{cuentasgamma})\\
			&=&[B(\tilde{f}(b),\tilde{f}(c)):B(\tilde{f}(a),\tilde{f}(c)):B(\tilde{f}(a),\tilde{f}(b))]\hfill(\ref{kernelprima})\\
			&=&\Gamma_f(a,b,c,y_0).\end{array}	
		\]
		 By A2 and the previous computation, if $a=x_0$, then 
		\[\begin{array}{rcl}
			\Gamma(x_0,b,c,y_0)&=&\Gamma\big(h\cdot (x_0,b,c,y_0)\big)\\
			&=&F_h\cdot\Gamma(c,b,x_0,y_0)\\
			&=&F_h\cdot [B(\tilde{f}(b),\tilde{f}(x_0)):B(\tilde{f}(c),\tilde{f}(x_0)):B(\tilde{f}(c),\tilde{f}(b))]\\
			&=&[B(\tilde{f}(b),\tilde{f}(c)):B(\tilde{f}(x_0),\tilde{f}(c)):B(\tilde{f}(x_0),\tilde{f}(b))]\\&=&\Gamma_f(x_0,b,c,y_0).
		\end{array}\]	
		If $b=x_0$ and  $y\in\setminus\{x_0,b,c,y_0\}$,   by  A2 and the previous calculations, 
		\[\begin{array}{rcl}
			\Gamma_{2,1}(a,x_0,c,y_0)&=&\Gamma_{2,1}(a,y,c,y_0)\Gamma_{1,2}(x_0,y,c,y_0)\\
			&=&\frac{B(\tilde{f}(a),\tilde{f}(c))}{B(\tilde{f}(y),\tilde{f}(c))}\frac{B(\tilde{f}(y),\tilde{f}(c))}{B(\tilde{f}(x_0),\tilde{f}(c))}
		\end{array}\]
		and
		\[\begin{array}{rcl}
			\Gamma_{2,3}(a,x_0,c,y_0)&=&\Gamma_{2,3}(a,y,c,y_0)\Gamma_{3,2}(a,y,x_0,y_0)\\
			&=&\frac{B(\tilde{f}(a),\tilde{f}(c))}{B(\tilde{f}(a),\tilde{f}(y))}\frac{B(\tilde{f}(a),\tilde{f}(y))}{B(\tilde{f}(a),\tilde{f}(x_0))}.
		\end{array}
		\]
		Thus, 
		\[\begin{array}{rcl}
			\Gamma(a,x_0,c,y_0)&=&\left[\frac{B(\tilde{f}(x_0),\tilde{f}(c))}{B(\tilde{f}(a),\tilde{f}(c))}:1:\frac{B(\tilde{f}(a),\tilde{f}(x_0))}{B(\tilde{f}(a),\tilde{f}(c))}\right]\\
			&=&[B(\tilde{f}(x_0),\tilde{f}(c)):B(\tilde{f}(a),\tilde{f}(c)):B(\tilde{f}(a),\tilde{f}(x_0))]\\
			&=&\Gamma_f(a,x_0,c,y_0),
		\end{array}\]
		concluding the proof for the points $(a,b,c,y_0)\in X^{\{4\}}.$
		
		Suppose  $x=(x_1,x_2,x_3,x_4)\in X^{\{4\}}$ is such that $y_0=x_i$, for some $i=1,2,3.$ There exists $g\in \mathrm{D}_8$ such that $4=g(i)$ and by the  previous computations
		$\Gamma_f(gx)=\Gamma(gx)$. As $\Gamma$ and $\Gamma_f$ are $\mathrm{D}_8$-equivariant,  
		$\Gamma(x)=\Gamma_f(x)$.  
			If $x=(x_1,x_2,x_3,x_4)\in X^{\{4\}}$  is such that  $y_0\neq x_i$, by A5 and   the previous computations,   $\Gamma(x)=\Gamma_f(x)$.
		
		Regarding the  continuity of $f$, it is  only missing to show it in $y_0$. Let  ${{\Psi}'_{y_0,x_1}}:(X\setminus\{x_0\})^2\rightarrow\C$ be the continuous map  associated to $\Gamma$ coming from A6. 
		If $\tilde{f}$ is the lift of $f$ defined in (\ref{kernelprima}), then for every $a\in X\setminus\{x_0\}$,
		\[\begin{array}{rcl}
			\Psi'_{y_0,x_1}(y_0,b)&=&\tfrac{B(\tilde{f}(y_0),\tilde{f}(b))}{B(\tilde{f}(y_0),\tilde{f}(x_0))B(\tilde{f}(x_0),\tilde{f}(b))}\tfrac{B(\tilde{f}(x_1),\tilde{f}(x_0))B(\tilde{f}(x_0),\tilde{f}(y_0))}{B(\tilde{f}(x_1),\tilde{f}(y_0))}\\
			&=&\tfrac{B(\tilde{f}(x_1),\tilde{f}(x_0))}{B(\tilde{f}(x_0),\tilde{f}(b))}.
		\end{array}\]
		Suppose  $(b_i)\to y_0$ and without lost of generality suppose $a_i\in X\setminus\{x_0\}.$ Thus as there exists $z_0\in\C$ such that $z_0\Psi'_{x_0,x_1}=\Psi_{x_0,x_1}$ is a continuous complex kernel of negative type, then
		\[
			\lim_{i\to\infty}\Psi'_{y_0,x_1}(y_0,b_i)=\lim_{i\to\infty}\tfrac{B(\tilde{f}(x_1),\tilde{f}(x_0))}{B(\tilde{f}(x_0),\tilde{f}(b_i))}=0. 
		\]
		This implies that \[\begin{array}{rcl}\lim\limits_{i\to\infty}|B(\tilde{f}(x_0),\tilde{f}(b_i))|^{\frac{1}{2}}&=&\lim\limits_{i\to\infty}d_x^{q_1}(q_2,f(b_i))\\
			&=&\lim\limits_{i\to\infty}\frac{d_x(q_2,f(b_i))}{d_x(q_1,f(b_i))}\\
			&=&\infty.
		\end{array}\]
		See  (\ref{definversion}) and the subsequent paragraph  for the definition of the metric $d_x^{q_1}$. As $d_x(q_2,f(b_i))\leq 1$, then  $\lim_{i\to\infty}d_x(q_1,f(a_i))=0$, and therefore,  the map $f$ is continuous in $y_0$. 
		
		Suppose   there exist $\mathbf{H}^{m'}_\C$    an injective (continuous) map $g:X\rightarrow\partial\mathbf{H}_\C^{m'}$ such that $\Gamma_g=\Gamma$. Denote $g(y_0)=p_1$ and $g(x_0)=p_2$. 
		Let $\mathbf{H}_\C^{m'}$ be induced from $(L',B')$, where $L'$ is a  complex Hilbert space and $B'$ is a strongly non-degenerate Hermitian  form of signature $(1,m')$ (see Section 3 in \cite{burger2005equivariant}). 
		
		By the arguments used above, there exists $\tilde{g}:X\setminus\{y_0\}\rightarrow L'$, a lift of $g$, such that for every $a\in X\setminus\{y_0\}$, 
		$B'(\tilde{g}(a),\tilde{g}(y_0))=1$. 
		Given  $\eta_1,\eta_2\in L'$, respective representatives of $p_1$ and $p_2$ such that $B'(\eta_1,\eta_2)=1$,  then there exists a canonical identification $\kappa:\partial\mathbf{H}^{m'}_\C\setminus\{p_1\}\rightarrow\mathbf{H}_{m'-1}$ such that $\kappa(p_2)=e$ and for every $a,b\in X\setminus\{y_0\}$,   $ B'(\tilde{g}(a),\tilde{g}(b))=\Psi_{x_0,x_1}(a,b).$
		By \cref{gnsnegativocomplejo}, there exists a unique isometric map ${\tilde{T}}:\mathbf{H}_{m-1}\rightarrow\mathbf{H}_{m'-1}$, such that for every $a,b\in X\setminus\{y_0\}$, 
		\[B_k\big((\tilde{T}\circ \omega \circ f(a),\tilde{T}\circ \omega \circ f(b)\big)=B_k( \kappa \circ g(a),\kappa \circ g(b)).\]
	Define $T=\kappa^{-1}\circ\tilde{T}\circ\omega$. 	   
		With a slight abuse of notation, the claim is that $\partial\mathbf{H}^m_\C\xrightarrow{{T}}\partial\mathbf{H}^{m'}_\C$, the natural extension of $T$,  is a Möbius map. Indeed, if $(a,b,c,q_1)\in\partial{\mathbf{H}^m_\C}^{\{4\}}$, 
		\[\begin{array}{rcl}
			crt_{d_x}(a,b,c,q_1)&=&\\
			\left[d_x(b,c)d_x(a,q_1):d_x(a,c)d_x(b,q_1):d_x(a,b)d_x(b,c)\right]&=&\\
			\left[\frac{d_x(b,c)}{d_x(b,q_1)d_x(q_1,c)}:\frac{d_x(a,c)}{d_x(a,q_1)d_x(q_1,c)}:\frac{d_x(a,b)}{d_x(a,q_1)d_x(q_1,b)}\right]&=&\\
			\left[d_x^{q_1}(b,c):d_x^{q_1}(a,c):d_x^{q_1}(a,b)\right]&=&\\
			\left[d_r^{p_1}({T}(b),{T}(c)):d_r^{p_1}({T}(a),{T}(c)):d_r^{p_1}({T}(a),{T}(b))\right]&=&\\
			crt_{d_r}({T}(a),{T}(b),{T}(c),{T}(q_1)).
		\end{array}\]
		For any other $x=(x_1,x_2,x_3,x_4)\in \partial{\mathbf{H}^m_\C}^{\{4\}}$, such that $q_0\in \{x_1,x_2,x_3,x_4\}$, the claim can be deduced from the previous computation and the identities in  A2. 
		
	 The map $T$ is induced by a unique isometric map $R:\mathbf{H}^m_\C\rightarrow\mathbf{H}^{m'}_\C$ (see \cref{extensionbourdoninfinito}). By construction,  
		${\tilde{T}}:\mathbf{H}_{m-1}\rightarrow\mathbf{H}_{m'-1}$ preserves $B_k$, hence  $T$ preserves the Cartan invariant and therefore,  $R$ is holomorphic. 	\end{proof}
	In  the previous theorem if $f$ is such that $\Gamma_f=\Gamma$,  the sufficient condition   for the map $f$ to be initial among the maps $g$ such that 	$\Gamma_g=\Gamma$ is that the image of $f$ is not contained in the visual boundary of any proper $\F$-hyperbolic subspace.

	A representation $\rho:G\rightarrow \mathrm{Isom}_\F(\mathbf{H}^m_\F)$ is called  \textit{continuous}  if  it is \textit{orbitally continuous}, that is to say that for every $x\in\mathbf{H}^m_\F$,  the map $g\mapsto\rho(g)x$ is continuous. 	
	
	If $\Gamma:X^{\{4\}}\rightarrow\mathbf{P}(\F^3)$ is a (continuous) kernel of $\F$-Möbius type, denote $\mathrm{Aut}(\Gamma)$ the group of (continuous) symmetries of $X$ that preserve $\Gamma$. 
	\begin{cor}\label{equivariante}
		Let $X^{(4)}\xrightarrow{\Gamma}\mathbf{P}(\F^3)$ be a (continuous) kernel of $\F$-Möbius type. If a pair  $(\mathbf{H}^m_\F,f)$, where $X\xrightarrow{f}\partial\mathbf{H}^m_\F$ is a (continuous) map such that $\Gamma_f =\Gamma$, is an initial pair with this property, then there exists a representation
		$\mathrm{Aut}(\Gamma)\xrightarrow{\rho}\mathrm{Isom}_\F(\mathbf{H}^m_\F)$ 
		such that for every $x\in X$ and every $g\in \mathrm{Aut}(f)$, 
		$\rho(g)f(x)=f(gx).$
			Moreover, if $\mathrm{Aut}(f)$ is provided with a group topology for which the action of  $\mathrm{Aut}(f)$ on $X$ is orbitally continuous, then $\rho$ is  continuous. 
	\end{cor}
	\begin{proof}
				For every $g\in\mathrm{Aut}(\Gamma)$ and  for every $x=(x_1,x_2,x_3,x_4)\in X^{(4)}$, \[\Gamma_{f\circ g}(x)=\Gamma(g(x_1),g(x_2),g(x_3),g(x_4))=\Gamma_f(x_1,x_2,x_3,x_4)=\Gamma(x),\]
		thus by \cref{GNSMobius}, there exists a unique  Möbius map $\rho(g):\partial\mathbf{H}^m_\F\rightarrow\partial\mathbf{H}^m_\F$ such that for every $x\in X$, $f(g(x))=\rho(g)(f(x))$. 	Notice that the map  $g\mapsto \rho(g)$ is a homomorphism because   for every $g_1,g_2\in \mathrm{Aut}(\Gamma)$, \[f\circ g_1\circ g_2=\rho(g_1)\circ\rho(g_2)\circ f,\] and  the map $\rho(g_1\circ g_2)$ is unique with this property, thus  $\rho(g_1)\circ\rho(g_2)=\rho(g_1\circ g_2)$.  
		
		If $\mathrm{Aut}(f)$ is provided with a  group topology coarser than the  pointwise convergence topology, it is clear that for every $x\in X$,  the map $\mathrm{Aut}(f)\rightarrow{\partial \mathbf{H}^m_\C}$, given by  $g\mapsto \rho(g)f(x)=f(gx)$, is continuous. 
		
		Fix $x_0,x_1,y_0\in X$  pairwise distinct and choose  $\xi_1,\xi_2\in L$ respective representatives of $f(y_0)$ and $f(x_0)$  such that $B(\xi_1,\xi_2)=1$. These choices induce a bijection  $j:\partial\mathbf{H}^m_\C\setminus\{f(y_0)\}\rightarrow\mathbf{H}_{m-1}$, where  $\mathbf{H}_{m-1}$ is  identified with $N\times \R$ and where $N$ is the following space. Consider the  complex kernel  of positive type $\Phi:(X\setminus\{y_0\})^2\rightarrow\C$  given by  \[\Phi(a,b)=\frac{1}{2}\big(\Psi(a,x_0)+\Psi(x_0,b)-\Psi(a,b)\big), \] and where $\Psi=\Psi_{x_0,x_1}$ is the complex kernel of negative type associated to $x_0,x_1$. By the GNS construction for positive type functions (see \cref{GNSpos})   there exist a complex Hilbert space $(N,\langle\,\meddot\,,\meddot\,\rangle)$ and a map $\iota:X\setminus\{y_0\}\rightarrow N$ such that for every $a,b\in X\setminus\{y_0\}$,  $\langle \iota(a),\iota(b)\rangle=\Phi(a,b)$ and such that $(N,\iota)$ is initial with these properties. 
		Furthermore, the map $j\circ f:X\setminus\{y_0\}\rightarrow{}N\times\R$ is such that 
		for every $a,b\in X\setminus\{y_0\}$,  	$j\circ f(a)=\big(\iota(x),\mathrm{Im}(\Psi(x,x_0)\big)$ and 
		$B_k(j\circ f(a),j\circ f(b))=\Psi(a,b).$
		
		The space  $L$ can be chosen to be  $\C^2\times N$  where $(1,0,0)$ and $(0,1,0)$ are identified with $\xi_1$ and $\xi_2$ respectively and where $L$ is  provided with the  Hermitian form given by,
		\begin{enumerate} \item $B\big((1,0,0),(0,1,0)\big)=1$. \item $N=(1,0,0)^\perp\cap(1,0,0)^\perp$. \item $B|_N=-\langle\,\meddot\,,\meddot\,\rangle|_N$. 
		\end{enumerate}
		Let ${\tilde{f}}:X\setminus\{y_0\}\rightarrow\C^2\times N$ be the map given by $\tilde{f}(a)=\big(\Psi(a,x_0),1,\sqrt{2}\iota(a)\big)$ and observe that $\tilde{f}$ is a lift of $f|_{X\setminus\{y_0\}}$ such that  for every $a,b\in  X\setminus\{y_0\}$, 
		\[B(\tilde{f}(a),\tilde{f}(b))=
			\Psi(a,x_0) +\Psi(x_0,b)-2 \Phi(a,b)=\Psi(a,b).\]
		Define $\tilde{f}(y_0)=(1,0,0) $
		and let    ${\tilde{g}}:\C^2\times N\rightarrow\C^2\times N$ in $\mathrm{U}(B)$ be the  representative of $\rho(g)$ such that  $\tilde{g}(1,0,0)= \tilde{f}(gy_0)$. Thus 
		$\tilde{g}(0,1,0)= \tfrac{\tilde{f}(gx_0)}{\Psi(gx_0,gy_0)}$
		and  for every $a\in X\setminus\{y_0,x_0\}$,\[ \tilde{g}(0,0,\sqrt{2}\iota(a))=\tfrac{\tilde{f}(ga)}{\Psi(ga,gy_0)}-\Psi(a,x_0)\tilde{f}(gy_0)-\tfrac{\tilde{f}(gx_0)}{\Psi(gx_0,gy_0)}.\] 
		
		The map $g\mapsto \tilde{g}$, that is not necessarily a homomorphism, defines a lift  $\tilde{\rho}:\mathrm{Aut}(f)\rightarrow\mathrm{U}(B)$ of $\rho$ such that if $v\in\C^2\times K$,   where $K$ is the complex vector space generated by $\mathrm{Im}(\iota)$, then the map $g\mapsto\tilde{\rho}(g)v$ is continuous. 
		
		Denote $x=\tfrac{1}{\sqrt{2}}(1,1,0)$ and let $\|\,\meddot\,,\meddot\,\|$ be the norm associated to the decomposition $L=\C x\oplus x^\perp.$ If $W\subset \mathrm{Aut}(f)$ is a symmetric  open neighbourhood of $\mathrm{Id}$ for which   there exists $M>0$ such that for every $g\in  W$, 
		$|B\big(\tilde{\rho}(g^{-1})x,x\big)|<M,$ then by \cref{controldelanorma}, there exists $M'>0$ such that for every  $g\in W$, $\|\tilde{\rho}(g)\|\leq M'$. 
		To conclude, observe that the space $\C^2\times K$ is total in $L$ and that for every $w\in L$, $v\in\C^2\times K$  and $g\in W\cdot W$, 
		\begin{equation*}\label{orbitafrontera}\begin{array}{rcl}
				\|\tilde{\rho}(g)w-w\|&\leq&\|\tilde{\rho}(g)w-\tilde{\rho}(g)v\|+\|\tilde{\rho}(g)v-v\|+\|v-w\|\\
				&\leq&(M'+1)\|v-w\|+\|\tilde{\rho}(g)v-v\|.
		\end{array}\end{equation*}
		This shows that the map $g\mapsto\tilde{\rho}(g)v$ is continuous with respect to the topology associated to the decomposition $L=\C x\oplus x^\perp$.  
	\end{proof}

	If  ${\Gamma}:X^{(4)}\rightarrow\mathbf{P}(\F^3)$ is an $\F$-kernel of Möbius type, then $\Gamma=\Gamma_f$ for some map 
	$f:X\rightarrow\partial\mathbf{H}_\F^m$. 
	For every pairwise distinct $a,b,c\in X$ and any lift $\tilde{f}$ of $f$, $\tfrac{B(\tilde{f}(a),\tilde{f}(b))}{B(\tilde{f}(a),\tilde{f}(c))B(\tilde{f}(c),\tilde{f}(b))}$ has  positive real part, thus for every $t\in(0,2)$, let     (independently of $\tilde{f}$) 
	${\Gamma^t}:X^{\{4\}}\rightarrow\mathbf{P}(\F^3) $ be the map given by, 
	\[\Gamma^t(a,b,c,d)=\left[\left(\tfrac{B\left(\tilde{f}(b),\tilde{f}(c)\right)}{B\left(\tilde{f}(b),\tilde{f}(d)\right)B\left(\tilde{f}(d),\tilde{f}(c)\right)}\right)^t:
	\left(\tfrac{B\left(\tilde{f}(a),\tilde{f}(c)\right)}{B\left(\tilde{f}(a),\tilde{f}(d)\right)B\left(\tilde{f}(d),\tilde{f}(c)\right)}\right)^t:\left(\tfrac{B\left(\tilde{f}(a),\tilde{f}(b)\right)}{B\left(\tilde{f}(a),\tilde{f}(d)\right)B\left(\tilde{f}(d),\tilde{f}(b)\right)}\right)^t\right].\] 
	
	\begin{lem}\label{argesaditivo} If $t\in(0,1)$ and the lift $\tilde{f}$ is such that there exists $d\in X$ such that for every $x\in X\setminus\{d\}$, $B(\tilde{f}(d),\tilde{f}(x))=1$, 
		then for every three distinct $a,b,c\in X\setminus\{d\}$, 
		\[\left(\tfrac{B(\tilde{f}(a),\tilde{f}(b))}{B(\tilde{f}(a),\tilde{f}(c))B(\tilde{f}(c),\tilde{f}(b))} \right)^t=\tfrac{B(\tilde{f}(a),\tilde{f}(b))^t}{B(\tilde{f}(a),\tilde{f}(c))^tB(\tilde{f}(c),\tilde{f}(b))^t} .\]
	\end{lem}
	\begin{proof}
			Let $a,b,c\in X\setminus\{d\}$ be pairwise distinct. The proof will be done by cases.  
		First suppose that \[\{\Arg(B(\tilde{f}(a),\tilde{f}(c))),\Arg(B(\tilde{f}(c),\tilde{f}(b)))\}\subset \{-\tfrac{\pi}{2},\tfrac{\pi}{2}\}.\]
	
		If   $\Arg(B(\tilde{f}(a),\tilde{f}(c)))$ and $\Arg(B(\tilde{f}(c),\tilde{f}(b)))$  have the opposite sign, 
		then 
		\[\left(\tfrac{B(\tilde{f}(a),\tilde{f}(b))}{B(\tilde{f}(a),\tilde{f}(c))B(\tilde{f}(c),\tilde{f}(b))} \right)^t=\tfrac{B(\tilde{f}(a),\tilde{f}(b))^t}{|B(\tilde{f}(a),\tilde{f}(c))|^t|B(\tilde{f}(c),\tilde{f}(b))|^t}=\tfrac{B(\tilde{f}(a),\tilde{f}(b))^t}{B(\tilde{f}(a),\tilde{f}(c))^tB(\tilde{f}(c),\tilde{f}(b))^t}.\]
		
		Suppose  now, without lost of generality, that   $\Arg(B(\tilde{f}(a),\tilde{f}(c)))$ and  $\Arg(B(\tilde{f}(c),\tilde{f}(b)))$  are equal to $\tfrac{\pi}{2}$.
The claim is that $\Arg(\tilde{f}(a),\tilde{f}(b))=\tfrac{\pi}{2}$. 	   
		Indeed, by (\ref{notacionfrontera}) for every $x\in X\setminus\{d\}$, 
		\[\tilde{f}(x)=ie_x\tilde{f}(d)+\xi+{v_x},\]
		where $\xi$ is one isotropic vector such that $B(\tilde{f}(d),\xi)=1$, $e_x\in\R$ and $v_x\in\tilde{f}(d)^\perp\cap \xi^\perp$. In this notation one has that $e_a>e_c>e_b$ and $v_a=v_b=v_c,$ and from this,  it follows that $B(\tilde{f}(a),\tilde{f}(b))B(\tilde{f}(c),\tilde{f}(a))$ is a strictly positive real number. Thus,
		\[\left(\tfrac{B(\tilde{f}(a),\tilde{f}(b))}{B(\tilde{f}(a),\tilde{f}(c))B(\tilde{f}(c),\tilde{f}(b))} \right)^t=\tfrac{B(\tilde{f}(a),\tilde{f}(b))^t}{B(\tilde{f}(a),\tilde{f}(c))^tB(\tilde{f}(c),\tilde{f}(b))^t}. \] 
		
		Notice  that  if  \[\{\Arg(B(\tilde{f}(a),\tilde{f}(c))),\Arg(B(\tilde{f}(c),\tilde{f}(b)))\}\not\subset \{-\tfrac{\pi}{2},\tfrac{\pi}{2}\},\]
		then \[\Arg\big(B(\tilde{f}(a),\tilde{f}(c))B(\tilde{f}(c),\tilde{f}(b))\big)=\Arg(B(\tilde{f}(a),\tilde{f}(c)))+\Arg(B(\tilde{f}(c),\tilde{f}(b))).
		\]
		If the real part of 
		$B(\tilde{f}(a),\tilde{f}(c))B(\tilde{f}(c),\tilde{f}(b))$ is non-positive, 
		as 
		\[\cart(f(a),f(c),f(b))\in\left[-\tfrac{\pi}{2},\tfrac{\pi}{2}\right]\] and 
		$\mathrm{Re}(B(\tilde{f}(b),\tilde{f}(a)))\geq0$, 
		then \[\Arg\Big(B(\tilde{f}(a),\tilde{f}(c))B(\tilde{f}(c),\tilde{f}(b))\Big)\Arg(\tilde{f}(b),\tilde{f}(a))<0.\] This implies that \begin{equation}\begin{array}{rcl}\label{cartanensuma}\cart(f(a),f(c),f(b))&=&
				\Arg\Big(B(\tilde{f}(a),\tilde{f}(c))B(\tilde{f}(c),\tilde{f}(b))\Big)+\Arg(\tilde{f}(b),\tilde{f}(a))\\
				&=&
				\Arg(B(\tilde{f}(a),\tilde{f}(c)))+\Arg(B(\tilde{f}(c),\tilde{f}(b)))+\Arg(B(\tilde{f}(b),\tilde{f}(a))).\end{array}\end{equation}
		If the real part of  $B(\tilde{f}(a),\tilde{f}(c))B(\tilde{f}(c),\tilde{f}(b))$ is  positive, then  (\ref{cartanensuma}) holds. 
		In the last two cases 
		\[t\cart(f(a),f(c),f(b))=
		t\Arg(B(\tilde{f}(a),\tilde{f}(c)))+t\Arg(B(\tilde{f}(c),\tilde{f}(b)))+t\Arg(B(\tilde{f}(b),\tilde{f}(a))).\qedhere\]
	\end{proof}
	
	\begin{teo}\label{kernelpowert}
		If $\Gamma$ is a (continuous) $\F$-kernel of Möbius type and $t\in(0,1)$, then the map $\Gamma^t$ is a (continuous) $\F$-kernel of Möbius type. 
	\end{teo}
	\begin{proof}
		By \cref{argesaditivo}, the map  $\Gamma^t$ has properties  A1 to  A5. 
		For property A6 observe that if $x_0,x_1,y_0$ are three distinct points in $X$, the map associated ${\Psi_{x_0,x_1}'}:(X\setminus\{y_0\})^2\rightarrow\F$ is such that for every $a,b\in X\setminus\{y_0\}$, 
		\[\Psi_{x_0,x_1}'(a,b)=\left(\tfrac{B(\tilde{f}(a),\tilde{f}(b))}{B(\tilde{f}(a),\tilde{f}(y_0))B(\tilde{f}(y_0),\tilde{f}(b))}\right)^t\left(\tfrac{{B(\tilde{f}(x_1),\tilde{f}(y_0))B(\tilde{f}(y_0),\tilde{f}(x_0))}}{B(\tilde{f}(x_1),\tilde{f}(x_0))} \right)^t.\]
		By \cref{powert},  the map 
		\[(a,b)\mapsto\left(\tfrac{B(\tilde{f}(a),\tilde{f}(b))}{B(\tilde{f}(a),\tilde{f}(y_0))B(\tilde{f}(y_0),\tilde{f}(b))}\right)^t \] is an $\F$-kernel of negative type. 
	\end{proof}
	\begin{teo}\label{kernel2}
		Let $X\xrightarrow{f}\partial\mathbf{H}_\C^1$ be an injective map and let $\Gamma$ be the associated tautological complex kernel of Möbius type. If $t\in(1,2)$, then $\Gamma^t$ is a complex kernel of Möbius type. 
	\end{teo}
	\begin{proof}
		If $\tilde{f}$ is a lift of $f$ such that for some $d\in X$ and for every $a\in X\setminus\{d\}$, $B(\tilde{f}(a),\tilde{f}(d))=1$, then for every distinct  $a,b\in X\setminus\{d\}$, $B(\tilde{f}(a),\tilde{f}(b))$ is purely imaginary. Thus for every $x=(a,b,c,d)\in X^{\{4\}}$, $\Gamma$ is such that   \[\Gamma(x)=[\Gamma_1(x):\Gamma_2(x):\Gamma_3(x)],\] where  $\Gamma_i(x)$ can be chosen to be  purely imaginary. By (\ref{tmayor2}) and the same arguments of \cref{argesaditivo},   
		\[\begin{array}{rcl}
			-\left(\tfrac{B(\tilde{f}(a),\tilde{f}(b))}{B(\tilde{f}(a),\tilde{f}(c))B(\tilde{f}(c),\tilde{f}(b))}\right)^{t}&=&\left(\tfrac{B(\tilde{f}(a),\tilde{f}(b))}{B(\tilde{f}(a),\tilde{f}(c))B(\tilde{f}(c),\tilde{f}(b))}\right)^{t,t-2}\\
			&=&
			\tfrac{B(\tilde{f}(a),\tilde{f}(b))^{t,t-2}}{B(\tilde{f}(a),\tilde{f}(c))^{t,t-2}B(\tilde{f}(c),\tilde{f}(b))^{t,t-2}}\\&=&
			
			-	\tfrac{B(\tilde{f}(a),\tilde{f}(b))^{t}}{B(\tilde{f}(a),\tilde{f}(c))^{t}B(\tilde{f}(c),\tilde{f}(b))^{t}}.	
		\end{array}\]	 
		Hence the map  
		$\Gamma^t=\left[\Gamma_1^t:\Gamma_2^t:\Gamma_3^t\right]$ has properties A1 to A5. 
		For  A6 notice that  if     $x_0,x_1\in X\setminus\{d\}$ are distinct, 
		then		for every $a,b\in X\setminus\{d\},$
		$\Psi'_{x_0,x_1}(a,b)=\tfrac{B(\tilde{f}(a),\tilde{f}(b))^t}
		{B(\tilde{f}(x_1),\tilde{f}(x_0))^t}$
		and by \cref{power2}, the map 
		$(a,b)\mapsto -{B(\tilde{f}(a),\tilde{f}(b))}^t$ is a complex kernel of negative type. 
	\end{proof}
	\subsection{Argument and displacement: a complete invariant}\label{seccioncomplete}
	Given a subgroup $G<\mathrm{Isom}_\F(\mathbf{H}_\F^m)$,     the \textit{limit set} of $G$, denoted $\Lambda G$, is defined as the set of  $q\in\partial\mathbf{H}^m_\F$ such that if $x_0\in \mathbf{H}_\F^m$,  there exists a sequence $(g_ix_0)$, with  $g_i\in G$,   converging to $q$. The set $\Lambda G$ is independent of  the choice of $x_0$ and it is initial among the closed $G$-invariant subsets of $\partial\mathbf{H}^m_\F$. The proofs of these claims and a general treatment of limit sets of subgroups of the isometry groups of the hyperbolic spaces can be found in Chapter 7 in \cite{das2017geometry}.
	
	A representation   $\rho:G\rightarrow\mathrm{Isom}_\F(\mathbf{H}_\F^m)$ is called \textit{non-elementary} if the action of $G$ on $\mathbf{H}_\F^m\cup\partial\mathbf{H}_\F^m$ does not have finite orbits. 
	If $\rho$ is non-elementary, then $\Lambda_\rho G$, the limit set of $\rho(G)$, is uncountable (see \textit{e.g.} Proposition 10.5.4 in \cite{das2017geometry}). 
		A representation $\rho$ is called \textit{irreducible} if it is non-elementary and  does not preserve any proper  $\F$-hyperbolic subspace of $\mathbf{H}_\F^m$. If $\rho$ is non-elementary, there exists  an $\F$-hyperbolic subspace,  the \textit{irreducible part} of $\rho$, that is $G$-invariant and initial with this property (see Proposition 4.3 in \cite{burger2005equivariant}).  
		
		If $G<\mathrm{Isom}(\mathbf{H}^n_\F)$ is non-elementary, denote $\Lambda^+G$ the set of attracting points of hyperbolic isometries in $G$. The set $\Lambda^+G$ is non-empty and $G$-invariant, and therefore, dense in $\Lambda G.$
	
	The following lemmas are classical (in the locally compact case) but their proofs in the non-locally compact setting  are presented here  for the sake of completeness. 
	\begin{lem}\label{uniformconvergence}
		If $g$ is an hyperbolic isometry of $\He_\F^m$, then for every open $U,V\subset\partial\He_\F^m$ such that $g_+\in U$ and $g_-\not\in{V}$, then there exists $N>0$ such that for every $n>N$, $g^n(V)\subset U$.  Moreover, there exists $p\in\mathbf{H}^m_\F$ such that   for every $q_1,q_2\not\in \overline{V}$ and every $n>N,$ \[d^{g_-}_p(g^n(q_1),g^n(q_2))\leq 2e^{-{\ell(g)n}}d^{g_-}_p(q_1,q_2).\]
		
	\end{lem}
	\begin{proof}
		Let $\xi_1$ and $\xi_2$ be  respective  lifts of $g_+$ and $g_-$  such that $B(\xi_1,\xi_2)=1$ and let  $ p=\left[(\xi_1+\xi_2)\right]$.
	  Without lost of generality, suppose that $g_+\not\in V$ and that for some $R>0$, $V$ is the complement of the ball of radius $R$ centred at $g_-$ of the metric $d_p^{g_+}$. Suppose also that  $U$ is the ball of radius $S$  centred at $g_+$ of the metric $d_p^{g_-}$.   Denote $\tilde{g}$ the linear  lift of $g$ such that $B(\tilde{g}(\xi_1),\xi_2)=e^{\ell(g)}.$ 
		
		Recall that for every  $\xi(v,b)$ such that $g_-\neq [\xi(v,b)]$ (see (\ref{notacionfrontera}) and  \cref{isometriakoran}),   \[d_p^{g_-}([\xi(v,b)],g_+)^2=d_p^{g_+}([\xi(v,b)],g_-)^{-2}=|B(\xi(v,b),\xi_2)|^{-1}.\]
			To conclude, notice that $\left[\tilde{g}\xi(v,b)\right]=\left[\xi(e^{\ell(g)}v,e^{2\ell(g)}b)  \right]$ and observe  that   
		\[d_p^{g_-} (g^n[\xi(v,b)],g_+)^2
			=
			e^{-2\ell(g)n}|B(\xi(v,b), \xi_2)|^{-1}
			=e^{-2\ell(g)n}d_p^{g_-} ([\xi(v,b)],g_+)^2.\qedhere\]
	\end{proof}
	If $p\in\mathbf{H}^m_\F$ and $d_p$ is the Bourdon metric associated, following \cite{das2017geometry},  for $g\in\mathrm{Isom}(\mathbf{H}^m_\F)$ such that  $\mathrm{Fix}(g)\cap\partial\mathbf{H}^m_\F\neq\emptyset$, define for $q\in \mathrm{Fix}(g)\cap\partial\mathbf{H}^m_\F$,
	${g}'(q)=\lim_{z\to q}\tfrac{d_p(gz,q)}{d_p(z,q)}. $
	The number $g'(q)$ is called the \textit{dynamical derivative} (with respect to $d_p$) of $g$ in $q$. A point $q\in \mathrm{Fix}(g)\cap\partial\mathbf{H}^m_\F$  is called an  \textit{attracting} (resp. \textit{repelling}) point of $g$ if 
	$g'(q)<1$ (resp. $g'(q)>1$). This notation is  coherent with that introduced in \cref{generalidadeshiperbolicos} (see Theorem 6.1.10 in \cite{das2017geometry}). 
	By Lemma 6.1.7 in \cite{das2017geometry}, if $g$ has either an attracting or a  repelling point, then $g$ is hyperbolic. 
	
	If  $G<\mathrm{Isom}(\mathrm{H}^m_\F)$ fixes a point $q\in\partial\mathbf{H}^m_\F$, define the \textit{Busemann quasi-character} of $G$ as the map ${B_{q}}:G\rightarrow\R$ given by,  $B(g)=  b_q(x)-b_q(gx)$, where $x\in \mathbf{H}^m_\F$ and $b_q$ is any Busemann function centred at $q$. The Busemann quasi-character is independent of  the choice of $x$. This map was studied   in  \cite{MonodCaprCornTess} in the context of group actions on Gromov-hyperbolic spaces. In the context of $\F$-hyperbolic spaces the  Busemann quasi-character is a homomorphism and  $|B_{q}(g)|=\ell(g)$ (see \textit{e.g.} p. 3729 in \cite{stolowicz2022complex}).  The relation between the Busemann quasi-character and the dynamical derivative is given by $e^{-B_q(g)}=g'(q)$ (see Remark 4.2.13 in \cite{das2017geometry}).  
	
	\begin{lem}\label{commutatorhyperbolic}
		If  $g$ and $h$ are two hyperbolic isometries  of $\mathbf{H}^m_\F$ such that  $Fix(g)\cap Fix(h)=\emptyset$, then the following hold.
		\begin{enumerate}
			\item There exists $N>0$ such that for every $n>N$, $g^nh^n$ is hyperbolic and $\lim_{n\to \infty}(g^nh^n)_+=g_+$ and $\lim_{n\to \infty}(g^nh^n)_-=h_-$.  
			\item There exist $m>0$ and $N>0$ such that for every $n>N$, $g^nh^mg^{-n}h^{-m}$ is hyperbolic. 
		\end{enumerate} 
	\end{lem}
	\begin{proof}
		For 1. let $U_+,U_-,V_+,V_-\subset \partial\mathbf{H}^m_\F$ be pairwise disjoint closed  neighbourhoods of $g_+,g_-,h_+,h_-$, respectively. Let $p,q\in \mathbf{H}^m_\F$ be  in the axis of $g$ and $h$ respectively and let $\lambda=e^{\ell(g)}$ and $\gamma=e^{\ell(h)}$. 
		
		 By \cref{uniformconvergence}, there exists $N>0$ such that for every $n>N$, $g^{\pm n}(\partial\mathbf{H}^m_\F\setminus U_{\mp})\subset U_\pm$ and $h^{\pm n}(\partial\mathbf{H}^m_\F\setminus V_{\mp})\subset V_\pm.$
		Suppose that for every $\eta\in U_+$, $d_q(\eta,h_-)\geq T$, for some $T>0$. For every $\eta\in U_+$ and every $n>N$, $h^n(\eta)\in V_+$, thus as the metrics $d_p$ and $d_q$ are bounded, there exists $S>0$ such that  for every $n>N$ and every $\eta_1,\eta_2\in U_+$, 
		\[\frac{d_q(h^n(\eta_1),h_-)d_q(h^n(\eta_2),h_-)}{d_p(h^n(\eta_1),g_-)d_p(h^n(\eta_2),g_-)}\leq S.\] 
		Observe that $g^nh^n(U_+)\subset U_+$ and by \cref{uniformconvergence}, for every $\eta_1,\eta_2\in V_+$, 
		\[\begin{array}{rcl}
			d^{g_-}_p\big(g^nh^n(\eta_1),g^nh^n(\eta_2)\big)&\leq&2\lambda^{-{n}}d^{g_-}_p\big(h^n(\eta_1),h^n(\eta_2)\big)\\
			&=&2\lambda^{-{n}}d^{h_-}_q\big(h^n(\eta_1),h^n(\eta_2)\big)\frac{d_q(h^n(\eta_1),h_-)d_q(h^n(\eta_2),h_-)}{d_p(h^n(\eta_1),g_-)d_p(h^n(\eta_2),g_-)}	\\
			&\leq &4\lambda^{-{n}}\gamma^{-{n}}d^{h_-}_q\big(\eta_1,\eta_2\big)S\\
			&=&4\lambda^{-{n}}\gamma^{-{n}}d^{g_-}_p(\eta_1,\eta_2)S\frac{d_p(\eta_1,g_-)d_p(\eta_2,g_-)}{d_q(\eta_1,h_-)d_q(\eta_2,h_-)}\\
			&\leq&4\lambda^{-{n}}\gamma^{-{n}}d^{g_-}_p(\eta_1,\eta_2)\frac{S}{T^2}.
		\end{array}\]
		Thus there exist $N'>0$ and $0<R<1$ such that for every $n>N'$, 
		\begin{equation}\label{contraction}d^{g_-}_p\big(g^nh^n(\eta_1),g^nh^n(\eta_2)  \big)<Rd^{g_-}_p(\eta_1,\eta_2  ).
		\end{equation}
		By the Banach fixed-point theorem,   there exists $\xi\in U_+$ fixed by $g^nh^n.$ Rewriting (\ref{contraction}) one gets, 
		\[d_p\big(g^nh^n\eta_1,g^nh^n\eta_2  \big)<Rd_p(\eta_1,\eta_2  )\tfrac{d_p\big(g^nh^n\eta_1,g_-\big)d_p\big(g^nh^n\eta_2,g_-\big)}{d_p(\eta_1,g_-)d_p(\eta_2,g_-)}.  \]
		
		As $\eta_i\in U_+$, there exists $T'>0$ such that  $d_p(\eta_i,g_-)>T'$. Thus there exists $0<R'<1$  such that  for every $n$ large enough,
		\[d_p\big(g^nh^n\eta_1,g^nh^n\eta_2  \big)<R'd_p(\eta_1,\eta_2  ).  \]
		This implies that $g^nh^n$ is hyperbolic (see the comments before the lemma). 
		
		For 2.  observe that $\mathrm{Fix}(hg^{-1}h^{-1})=h(\mathrm{Fix}(g))$, thus the claim follows from 1. because for $m$ large enough, $\mathrm{Fix}(h^m g h^{-m})\cap \mathrm{Fix}(g)=\emptyset.$
	\end{proof}

	Given a group $G$ and two non-elementary   representations $\rho, \tau:G\rightarrow\iso_\F(\He_\F^m)$, we say that $\rho$ and $\tau$ \textit{have the same displacement} if for every $g\in G$, $\ell(\rho(g)) =\ell(\tau(g)).$
	
	The arguments used in the  following lemmas,  in the locally compact setting,  are contained in Lemma 2 of \cite{Kimmarkedlenght}. 
	\begin{lem}\label{welldefined}
		If $\rho, \tau:G\rightarrow\iso_\F(\He_\F^m)$ are two non-elementary  representations with the same displacement, then the following hold.
		\begin{enumerate}
			\item If $g,h\in G$ are mapped to two hyperbolic isometries then 
			$Fix(\rho(g))\cap Fix(\rho(h))=\emptyset$ if, and only if, $Fix(\tau(g))\cap Fix(\tau(h))=\emptyset$. 
			\item $\rho(g)_+=\rho(h)_+$ if, and only if, $\tau(g)_+=\tau(h)_+$.
		\end{enumerate}
	\end{lem}		
	\begin{proof}
		For 1., let $g,h\in G$ be such that $\rho(g)$ and $\rho(h)$ are hyperbolic isometries.  By \cref{commutatorhyperbolic} if   $\mathrm{Fix}(\rho(g))\cap \mathrm{Fix}(\rho(h))=\emptyset$, then  there exists $m>0$ such that  $[\rho(g)^m,\rho(h)^m]$  is  hyperbolic.  	Every isometry contained in the commutator group of $\mathrm{Stab}(q)<\mathrm{Isom}_\F(\mathbf{H}^m_\F)$, for some  $q\in\partial\mathrm{H}^m_\F$, has to be either  parabolic or elliptic (see the comments before \cref{commutatorhyperbolic}). Thus as both representations have the same displacement,  $\mathrm{Fix}(\tau(g))\cap \mathrm{Fix}(\tau(h))=\emptyset.$  
		
		For 2. suppose that $\rho(g)_+=\rho(h)_+=q$.  
		By the comments before \cref{commutatorhyperbolic}, \[\ell(\rho(gh))=|B_q(\rho(gh))|=|B_q(\rho(g))+B_q(\rho(h))|=\ell(\rho(g))+\ell(\rho(h)).\]
		If   $\tau(g)_+=\tau(h)_-$,. 
		\[\ell(\tau(gh))=|B_\xi(\tau(gh))|=|B_\xi(\tau(g))+B_\xi(\tau(h))|\neq \ell(\tau(g))+\ell(\tau(h)),\] which is a contradiction because $\rho$ and $\tau$ have same displacement. 
	\end{proof}
	The following lemma is contained in  Theorem 1 in \cite{Kimmarkedlenght}, where it  is stated for locally compact rank-1 symmetric spaces of non-compact type. This proof is given here to clarify that in the  non-locally compact setting   the same arguments in  the aforementioned work can be used to prove   \cref{crossratiolimitset}. Then  following lemma covers the only  complication that the non-locally compact case could bring.  
	\begin{lem}\label{complication}
		Let $\tau$ be a geodesic in $\mathbf{H}^m_\F$ with attracting and repelling  points $\tau_+$ and $\tau_-$. Suppose that  for every $i$, $\sigma_i$ is  a geodesic with attracting and repelling points $\sigma_{+,i}$ and $\sigma_{-,i}$ respectively, and suppose  that $\lim\limits_{i\to\infty}\sigma_{+,i}=\tau_+$ and $\lim\limits_{i\to\infty}\sigma_{-,i}=\tau_-$.   Then there exist sequences $\sigma_i(t_i)$ and $\tau(t_i)$, with $\lim\limits_{i\to\infty}t_i=\infty$,  such that 
		$\lim\limits_{i\to\infty}d\big(\sigma_i(t_i),\tau(t_i))=0.$
	\end{lem}
	\begin{proof}
		Fix $p\in \mathbf{H}^m_\F$ contained in $\tau$ and fix  $\tilde{p}$ a lift of $p$ such that $B(\tilde{p},\tilde{p})=1$. Every $\omega\in\partial\mathbf{H}^m_\F$ admits a representative $\tilde{\omega}=\tilde{p}+v$, with $B(v,v)=-1$   and recall (see \ref{representantes}) that 
		\[d_p([\tilde{p}+v_1],[\tilde{p}+v_2])^2=\tfrac{1}{2}|B(\tilde{p}+v_1,\tilde{p}+v_2)|.\]
		Suppose that for every $i$, $\sigma_i(0)=\left[{\tilde{\sigma}_{+,i}+\tilde{\sigma}_{-,i}}\right]$ and that $\tau(0)=\left[{\tilde{\tau}_{+}+\tilde{\tau}_{-}}\right]$.  
		Let $(t_i)$ be such that  $\lim_{i\to\infty}t_i=\infty$ and that  $\lim_{i\to\infty}e^{t_i}d_p(\tau_+,\sigma_{+,i})=0.$
		Define for every $i$, \[f(i)=\tfrac{\left|B\big(e^{t_i}\tilde{\tau}_+ +e^{-t_i}\tilde{\tau}_-,e^{t_i}\tilde{\sigma}_{+,i}+e^{-t_i}\tilde{\sigma}_{-,i} \big)\right|}{2|B(\tilde{\tau}_+,\tilde{\tau}_-
			)|^{\frac{1}{2}}|B(\tilde{\sigma}_{+,i},\tilde{\sigma}_{-,i})|^{\frac{1}{2}}}.\]
	If $\tilde{\sigma}_{\pm,i}=\tilde{p}+u_{\pm,i}$ and $\tilde{\tau}_+=\tilde{p}+w$, then $\lim_{i\to\infty}u_{+,i}= w$ and $\lim_{i\to\infty}u_{-,i}= -w$. By definition, 	
		\[
			\lim\limits_{i\to\infty}	f(i)			=\lim\limits_{i\to\infty}\tfrac{|B(\tilde{\tau}_{-},\tilde{\sigma}_{+,i})+B(\tilde{\tau}_{+},\tilde{\sigma}_{-,i})|}{2|B(\tilde{\tau}_+,\tilde{\tau}_-
				)|}=1. 			 
		\]
		To conclude, observe that 
		\[
		\lim\limits_{i\to\infty}d\big(\tau(e^{t_i}), \sigma_i(e^{t_i}))=\lim\limits_{i\to\infty}\ln\big(f(i)+\sqrt{f(i)^2-1}\big)=0.		
		\qedhere\]
		\end{proof}
	Given pairwise distinct $x_1,\dots,x_4\in\partial\mathbf{H}^m_\F$ and a Bourdon metric $d_p$, for some $p\in\mathbf{H}^m_\F$,  denote  
	\[[x_1,x_2,x_3,x_4]=\frac{d_p(x_1,x_3)d_p(x_2,x_4)}{d_p(x_1,x_4)d_p(x_2,x_3)}.\]
	The proof of the following theorem is, \textit{mutatis mutandis}, the same as that  of   Theorem 1 in \cite{Kimmarkedlenght} (see \cref{complication} and the comment before it).  
	\begin{teo}\label{crossratiolimitset}
		If $g$ and $h$ are two hyperbolic isometries of $\mathbf{H}^m_\F$ such that $Fix(g)\cap Fix(h)=\emptyset$, then \[\lim\limits_{i\to\infty}\exp\big(\ell(g^i)+\ell(h^i)-\ell(h^ig^i)\big)=[g_-,h_-,g_+,h_+].\]
	\end{teo}
	Two non-elementary representations  $\rho,\tau: G\rightarrow\mathrm{Isom}_\F(\mathbf{H}^m_\F)$    \textit{have the same argument} if for every $g_1,g_2,g_3\in G$ such that    $\rho(g_i)$ (or equivalently $\tau(g_i)$) are hyperbolic and the respective attracting points are pairwise distinct (see \cref{welldefined}), then 
	\[\mathrm{Cart}\big(\rho(g_1)_+,\rho(g_2)_+,\rho(g_3)_+\big)=\mathrm{Cart}\big(\tau(g_1)_+,\tau(g_2)_+,\tau(g_3)_+\big).\]
	
	\begin{teo}\label{dispargclas}
		Two irreducible representations $\rho,\tau: G\rightarrow\mathrm{Isom}_\F(\mathbf{H}^m_\F)$ are conjugated by a holomorphic isometry if, and only if, they have the same displacement and argument. 
	\end{teo}  	 
	\begin{proof}
		Denote $\Gamma_\rho$ (resp. $\Gamma_\tau$) the kernel of Möbius type associated to the inclusion $\Lambda_\rho G\hookrightarrow\partial \mathbf{H}^m_\F$
		(resp. $\Lambda_\tau G\hookrightarrow\partial \mathbf{H}^m_\F$). The  claim is that there exists $f:\Lambda_\rho(G)\rightarrow\Lambda_\tau(G)$, a continuous $G$-equivariant  bijection,   
		such that for every $(x_1,x_2,x_3,x_4)\in(\Lambda_\rho G)^{\{4\}}$, \[\Gamma_\rho(x_1,x_2,x_3,,x_4)=\Gamma_\tau(f(x_1),f(x_2),f(x_3),f(x_4)).\] 
				
		Following \cite{Kimmarkedlenght} and by \cref{welldefined},  define $f:\Lambda^+_\rho G\rightarrow\Lambda^+_\tau G$ by declaring  $f(\rho(g)_+)=\tau(g)_+$. By construction the map $f$ is $G$-equivariant. 
		
		Let 
		\[x=(\rho(g)_-,\rho(h)_-,\rho(g)_+,\rho(h)_+)\in(\Lambda_\rho^+ G)^{\{4\}}\] and \[y=(\tau(g)_-,\tau(h)_-,\tau(g)_+,\tau(h)_+)\in(\Lambda_\tau^+ G)^{\{4\}}.\]
		Denote $a_{\pm},b_{\pm},\alpha_{\pm},\beta_{\pm}$ the  respective lifts of $\rho(g)_{\pm},\rho(h)_{\pm},\tau(g)_{\pm},\tau(h)_{\pm}$, such that \[B(a_{\pm}, b_+)=B(b_-,b_+)=B(\alpha_{\pm}, \beta_+)=B(\beta_-,\beta_+)=1.\]
		Thus,   
		\[\begin{array}{rcl}\Gamma_\rho(x)&=&\left[B(b_-,a_+):{B(a_-,a_+)}:{B(a_-,b_-)}\right]\\
			&=&\left[|B(b_-,a_+)|e^{ir_1}:|{B(a_-,a_+)}|e^{ir_2}:|{B(a_-,b_-)}|e^{ir_3}\right]\end{array}\] and \[\begin{array}{rcl}\Gamma_\tau(y)&=&\left[B(\beta_-,\alpha_+):{B(\alpha_-,\alpha_+)}:{B(\alpha_-,\beta_-)}\right]\\&=&\left[|B(\beta_-,\alpha_+)|e^{is_1}:|{B(\alpha_-,\alpha_+)}|e^{is_2}:|{B(\alpha_-,\beta_-)}|e^{is_3}\right],\end{array}\]
		where 
		\[\begin{array}{ccccccc}
			r_1&=&\cart(\rho(h)_-,\rho(g)_+,\rho(h)_+)&=&\cart(\tau(h)_-,\tau(g)+,\tau(h)_+)&=&s_1,\\
			r_2&=&\cart(\rho(g)_-,\rho(g)_+,\rho(h)_+)&=&\cart(\tau(g)_-,\tau(g)_+,\tau(h)_+)&=&s_2,\\
			r_3&=&\cart(\rho(g)_-,\rho(h)_-,\rho(h)_+)&=&\cart(\tau(g)_-,\tau(h)_-,\tau(h)_+)&=&s_3. 
		\end{array}\]					
		By (\ref{equa}) and \cref{crossratiolimitset},  
		\[\left|\tfrac{B(a_-,a_+)}{B(b_-,a_+)}\right|^{\frac{1}{2}}=[\rho(g)_-,\rho(h)_-,\rho(g)_+,\rho(h)_+]=[\tau(g)_-,\rho(h)_-,\tau(g)_+,\tau(h)_+]=\left|\tfrac{B(\alpha_-,\alpha_+)}{B(\beta_-,\alpha_+)}\right|^{\frac{1}{2}}\] and \[\left|\tfrac{B(a_-,b_-)}{B(b_-,a_+)}\right|^{\frac{1}{2}}=[\rho(g)_-,\rho(g)_+,\rho(h)_-,\rho(h)_+]=[\tau(g)_-,\tau(g)_+,\tau(h)_-,\tau(h)_+]=\left|\tfrac{B(\alpha_-,\beta_-)}{B(\beta_-,\alpha_+)}\right|^{\frac{1}{2}},\] 
	hence $\Gamma_\rho(x)=\Gamma_\tau(y)$.	This implies that $f$ is a homeomorphism (see \cref{mobcontinua}).   
	
	  The set 
	\[\{(\rho(g)_+,\rho(g)_-)\mid g\in G, \rho(g)\,\textrm{is hyperbolic}\}\] is dense in $(\Lambda_\rho G)^2$ (see \textit{e.g.} Theorem 7.4.7 in \cite{das2017geometry}). Then, for every $x\in (\Lambda_\rho^+G)^{\{4\}}$, with a slight abuse of notation, $\Gamma_\rho (x)=\Gamma_\tau(f(x))$. 
		
The claim now is that $f$ can be extended continuously  to $\Lambda_\rho G$. 	Suppose that $(\eta_i)$ is a Cauchy sequence in $\Lambda^+_\rho G$.  Let $\xi_1,\xi_2\in\Lambda^+_\rho G$ be the attracting and repelling points of an hyperbolic isometry in $\rho(G)$  such that $(\eta_i)$ does not converge to neither $\xi_1$ nor  $\xi_2$. This implies that $(f(\eta_i))$ does not accumulate in neither $f(\xi_1)$ nor $f(\xi_2)$.  Let $x$ (resp. $y$) be contained in the geodesic connecting $\xi_1$ and $\xi_2$	(resp. $f(\xi_1)$ and $f(\xi_2)$).  Thus, without lost of generality, there exists $T>0$ such that  for every $i,j$, 
\[d_y(f(\eta_i),f(\eta_j))\leq\tfrac{d_y(f(\eta_i),f(\eta_j))}{d_y(f(\xi_1),f(\eta_j))d_y(f(\xi_2),f(\eta_i))}=\tfrac{d_x(\eta_i,\eta_j)}{d_x(\xi_1,\eta_j)d_x(\xi_2,\eta_i)}\leq Td_y(\eta_i,\eta_j).\]
Hence $f$ can be extended to a homeomorphism $f:\Lambda_\rho G\rightarrow \Lambda_\tau G$ such that, with a slight abuse of notation, 
for every $x\in (\Lambda_\rho G)^{\{4\}}$, $\Gamma_\rho(x)=\Gamma_\tau(f(x))$ (see \cref{cartancontinuous}). 
By \cref{GNSMobius} and the comment after it, $f$  has a unique extension to a  bijective Möbius map $F:\partial\mathbf{H}^m_\F\rightarrow\partial\mathbf{H}^m_\F$. 
		
For every $g\in G$ and $\eta\in \Lambda_\rho G$, $\tau(g)^{-1}F(\rho(g)\eta)=f(\eta)$, thus as $F$ is the unique extension of $f$, then $F$ is $G$-equivariant. If $A\in\mathrm{Isom}_\F(\mathbf{H}^m_\F)$ induces $F$, then  $A$ conjugates $\rho$ and $\tau$. 	 
	\end{proof} 
If $\F=\R$ in the previous theorem  the condition on the Cartan arguments is void, thus two such representations are conjugated if, and only if, they have the same displacement.

	\section{Infinite-dimensional representations} 
	In \cref{seccionkernelesmobius} it was introduced   a method to deform kernels of Möbius type  by \say{taking powers}. This method will be used in  \cref{deformations}  to describe new representations. 
In   \cref{algunosejemplos},  the complex   hyperbolic representations of $\mathrm{PSL}_2(\R)$ are classified   using the aforementioned  method, together with the \textit{horospherical combination} of representations (see \cite{stolowicz2022complex}) and the  results obtained in \cref{kernelesnegativos}  for functions of complex negative type. 
	
	By the Karpelevich-Mostow theorem (see  \cite{boubelkarpelevichmost,Mostowmostowkarpelevich}), the representations of $\mathrm{PU}(1,n)$ and $\mathrm{PO}(1,n)$  coming from a \say{power} have to be infinite-dimensional. When the group considered is  a non-elementary subgroup of the  groups above another argument is needed; in \cref{restrictionofrepresentations} this problem is  addressed.	
	It will also be shown  that the properties \textit{strongly discrete}, \textit{convex-cobounded} and \textit{geometrically finite} are preserved when  the representations introduced in this work are considered. 
	\subsection{Deforming  representations}\label{deformations}

	Given a representation $\rho$, let  $\Gamma$ be the tautological kernel of Möbius type associated to  $\Lambda_\rho  G\hookrightarrow\partial\mathbf{H}_\F^m$. If $m>1$ (resp. $\F=\C$ and $n=1$), for every  $t\in(0,1)$ (resp.   $t\in(0,2)$),  	by \cref{equivariante},  the kernel $\Gamma^t$  has associated 
 a representation ${\rho_t}:G\rightarrow\mathrm{Isom}_\F(\mathbf{H}_\F^{m'})$ (resp. ${\rho_{t}}:G\rightarrow\mathrm{Isom}_\C(\mathbf{H}_\C^{m'})$). This subsection will be devoted to study these \say{deformations} of  representations.
	
	The following theorem is an immediate consequence of  \cref{GNSMobius}.
	\begin{teo}\label{respectirreducible}
		If $\rho$ is a continuous non-elementary  representation of a group $G$, then $\rho_t$ is an irreducible representation.
	\end{teo}
	\begin{prop}\label{cartant1}
		Let $G<\mathrm{Isom}_\C(\mathbf{H}^m_\C)$ be  non-elementary and let $\Gamma$ be the tautological kernel of Möbius type associated. If  $t\in(0,1)$ and  $f:\Lambda(G)\rightarrow\Lambda_{\rho_t}G$ is the $G$-equivariant map associated to $\Gamma^t$, then for every  $(x_1,x_2,x_3)\in \Lambda_\rho G^{\{3\}}$, 
		\[\mathrm{Cart}(f(x_1),f(x_2),f(x_3))=t\mathrm{Cart}(x_1,x_2,x_3).\] 
	\end{prop}
	\begin{proof}
		Let $d\in \Lambda_\rho(G)\setminus\{x_i\}$ and fix $\tilde{d}$ one of its representatives.  Take $\tilde{x}_i$  and $\tilde{f}(x_i)$, lifts of $x_i$ and $f(x_i)$ respectively, such that  $B(\tilde{x}_i,\tilde{d})=1=B(\tilde{f}(x_i),\tilde{d}).$ On the one hand,
		\begin{equation*}\label{cocientes1}	\begin{array}{rcl}
				\Gamma_{1,2}(x_1,x_2,x_3,d)\Gamma_{2,1}(x_3,x_2,x_1,d)\Gamma_{1,2}(x_3,x_1,x_2,d)&=&
				\tfrac{B(\tilde{x}_2,\tilde{x}_3)}{B(\tilde{x}_1,\tilde{x}_3)}
				\tfrac{B(\tilde{x}_3,\tilde{x}_1)}{B(\tilde{x}_2,\tilde{x}_1)}	\tfrac{B(\tilde{x}_1,\tilde{x}_2)}{B(\tilde{x}_3,\tilde{x}_2)}\\
				&=&
				\exp\big(i2\mathrm{Cart}(x_1,x_2,x_3)\big).
		\end{array}\end{equation*}
		On the other hand, by \cref{argesaditivo},  \begin{equation*}\label{cocientes2}\begin{array}{rcl}
				\Gamma^t_{1,2}(x_1,x_2,x_3,d)\Gamma^t_{2,1}(x_3,x_2,x_1,d)\Gamma^t_{1,2}(x_3,x_1,x_2,d)&=&
				\tfrac{B(\tilde{x}_2,\tilde{x}_3)^t}{B(\tilde{x}_1,\tilde{x}_3)^t}
				\tfrac{B(\tilde{x}_3,\tilde{x}_1)^t}{B(\tilde{x}_2,\tilde{x}_1)^t}	\tfrac{B(\tilde{x}_1,\tilde{x}_2)^t}{B(\tilde{x}_3,\tilde{x}_2)^t}\\&=&
				\tfrac{B(\tilde{f}(x_2),\tilde{f}(x_3))}{B(\tilde{f}(x_1),\tilde{f}(x_3))}
				\tfrac{B(\tilde{f}(x_3),\tilde{f}(x_1))}{B(\tilde{f}(x_2),\tilde{f}(x_1))}	\tfrac{B(\tilde{f}(x_1),\tilde{f}(x_2))}{B(\tilde{f}(x_3),\tilde{f}(x_2))}

				\\&=&
				{\exp\big(i2\mathrm{Cart}(f(x_1),f(x_2),f(x_3))\big)}.
		\end{array}\end{equation*}
		Again by \cref{argesaditivo},  \[\begin{array}{rcl}
			\tfrac{B(\tilde{x}_2,\tilde{x}_3)^t}{B(\tilde{x}_1,\tilde{x}_3)^t}
			\tfrac{B(\tilde{x}_3,\tilde{x}_1)^t}{B(\tilde{x}_2,\tilde{x}_1)^t}	\tfrac{B(\tilde{x}_1,\tilde{x}_2)^t}{B(\tilde{x}_3,\tilde{x}_2)^t}&=&\left(\tfrac{B(\tilde{x}_2,\tilde{x}_3)}{B(\tilde{x}_1,\tilde{x}_3)B(\tilde{x}_2,\tilde{x}_1)}\right)^t\left(\tfrac{\overline{B(\tilde{x}_2,\tilde{x}_3)}}{\overline{B(\tilde{x}_1,\tilde{x}_3)}\overline{B(\tilde{x}_2,\tilde{x}_1)}}\right)^{-t}\\&=&\exp(2t\cart(x_1,x_2,x_3)).
		\end{array}\]
		Thus, as $|2(1-t)\mathrm{Cart}(x_1,x_2,x_3)|<\pi$,
		 \[\mathrm{Cart}(f(x_1),f(x_2),f(x_3))=t\mathrm{Cart}(x_1,x_2,x_3).\qedhere\] 	\end{proof}
	
	\begin{prop}\label{cartant2}
		Let $G<\mathrm{Isom}_\C(\mathbf{H}^1_\C)$ be non-elementary and let $\Gamma$ be the tautological kernel of Möbius type associated to $G$. If $t\in(1,2)$ and $f:\Lambda (G)\rightarrow\Lambda_{\rho_{t}}G$ is the $G$-equivariant map associated to $\Gamma^{t}$, then for every  $(x_1,x_2,x_3)\in \Lambda_\rho G^{\{3\}}$, 
		\[\mathrm{Cart}(f(x_1),f(x_2),f(x_3))=(t-2)\mathrm{Cart}(x_1,x_2,x_3).\]  
	\end{prop}
	\begin{proof}
		In the proof of \cref{kernel2} one can find  the analogous statement of \cref{argesaditivo} but for the case $t\in(1,2)$. The arguments for the proof follow, \textit{mutatis mutandis}, those    in  \cref{cartant1}.
	\end{proof}
	It remains to compute the \textit{displacement} of the representations coming from \say{deformations} of  kernels. Following Section 1.3 in \cite{surlebirapportBourdon}, for $(a,b,c)\in(\partial\mathbf{H
	}^m_\F)^{\{3\}},$
	denote $p=p(a,b,c)\in   \mathbf{H
	}^m_\F$, the unique point in the geodesic connecting $a$ and $b$ such that $d_{p}(a,d)=d_p(b,d)$, where $d_p$ is the Bourdon metric centred at $p$.  If $(a,b,c,d)\in (\partial\mathbf{H
	}^m_\F
	)^{\{4\}}$ and $x\in\mathbf{H
	}^m_\F$, then  \begin{equation}\label{desplazamiento}
		\left|\ln\left(\tfrac{d_x(a,c)d_x(b,d)}{d_x(a,d)d_x(b,c)} \right)\right|=d\big(p(a,b,c),p(a,b,d)\big).\end{equation}
	
	\begin{prop}\label{rescalamientodeldesplazamiento}
		Let $G<\mathrm{Isom}_\F(\mathbf{H}^m_\F)$ be non-elementary and let $\rho_t$ be the irreducible representation associated to $\Gamma^t$, where $\Gamma$ is the tautological kernel of Möbius type associated to $G$. 	Then for every $g\in G$, $\ell(\rho_t(g))=t\ell(g).$
	\end{prop}
	\begin{proof}
		
		Let $g\in G$ be hyperbolic and let  $f:\Lambda G \rightarrow\Lambda_\rho G$ be the $G$-equivariant map associated to $\rho_t$. Let    $\eta\in\partial \mathbf{H}^m_\F\setminus\{g_+,g_-\}$ and  consider $p=p(g_+,g_-,\eta)$.  By definition $gp=p(g_+,q_-,g\eta)$ and 
		 by  the comment after  (\ref{definversion}),  \cref{isometriakoran} and (\ref{gamavalorabsoluto}), for every $x\in\mathbf{H}^m_\F$,   \begin{equation*}\left|\ln\left(|\Gamma|^{\frac{1}{2}}_{2,1}(g_+,g_-,\eta,g\eta)\right)\right|=\left|\ln\left(\tfrac{d_x(g_+,\eta)d_x(g_-,g\eta)}{d_x(g_+,g\eta)d_x(g_-,\eta)} \right)\right|=d(p,gp),\end{equation*}
		and for the same reasons,  
		\[\left|\ln\left(|\Gamma^t|^{\frac{1}{2}}_{2,1}(g_+,g_-,\eta,g\eta)\right)\right|=d\Big(p\big(f(g_+),f(g_-),f(\eta)\big),gp\big(f(g_+),f(g_-),f(\eta)\big)\Big).\]
		By definition, 
		\[\ln\left|\left(|\Gamma^t|^{\frac{1}{2}}_{2,1}(f(g_+),f(g_-),f(\eta),f(g\eta))\right)\right|=t\left|\ln\left(|\Gamma|^{\frac{1}{2}}_{2,1}(g_+,g_-,\eta,g\eta)\right)\right|,\] 
		and therefore, 
		\[d\Big(p\big(f(g_+),f(g_-),\eta\big),gp\big(f(g_+),f(g_-),\eta\big)\Big)=td(gp,p).\] Hence  $\rho_t(g)$ is hyperbolic  and  $\ell(\rho_t(g))=t\ell(g).$ 
			To conclude,  observe that if   $\rho_t(g)$ is hyperbolic, then  $\rho_t(g)_+$ and $\rho_t(g)_-$ are contained in $\Lambda_{\rho_t}G=f(\Lambda G)$, hence by the arguments above, $g$ is hyperbolic. 
	\end{proof}
	\begin{prop}\label{mapeoequivariante}
		If $G<\mathrm{Isom}_\F(\mathbf{H}^m_\F)$ is non-elementary and $\rho:G\rightarrow \mathrm{Isom}_\F(\mathbf{H}^{m'}_\F)$ is a  representation such that there exists $t>0$ such that for every $g\in G$, $\ell(\rho(g))=t\ell(g)$, then $\rho$ is non-elementary and there exists a unique  continuous  $G$-equivariant map $f:\Lambda G\rightarrow \Lambda_\rho G$ such that for every $x\in\mathbf{H}^m_\F$, $y\in\mathbf{H}^{m'}_\F$ and $\eta\in\Lambda G,$ there exists $\lambda>0$ such that for every $a,b\in\Lambda G\setminus\{\eta\}$, 
		\[d_x^{\eta}(a,b)^t=\lambda d_y^{f(\eta)}(f(a),f(b)).\]
	\end{prop}
	\begin{proof}
		The arguments  in \cref{welldefined} can be adapted to show that one can define  the map $f:\Lambda^+ G\rightarrow \Lambda_\rho G$, given by $f(g_+)=\rho(g)_+.$  
		By \cref{crossratiolimitset}, for every $g,h\in G$ hyperbolic, $x\in\mathbf{H}^m_\F$ and $y\in\mathbf{H}^{m'}_\F$, 
		\[\begin{array}{rcl}\tfrac{d_y(f(g_-),f(g_+))d_y(f(h_-),f(h_+))}{d_y(f(g_-),f(h_+))d_y(f(h_-),f(g_+))}&=&[f(g)_-,f(h)_-,f(g)_+,f(h)_+]\\&=&[g_-,h_-,g_+,h_+]^t\\&=&\tfrac{d_x(g_-,g_+)^td_x(h_-,h_+)^t}{d_x(g_-,h_+)^td_x(h_-,g_+)^t}.\end{array}\]
		
		Suppose  $x$ and  $y$ are contained in the geodesics connecting $h_+$ and $h_-$ and $f(h_+)$ and $f(h_-)$ respectively. Observe that 
		\[d_y(f(g_-),f(g_+))=\tfrac{d_x(g_-,g_+)^td_y(f(g_-),f(h_+))d_y(f(h_-),f(g_+))}{d_x(g_-,h_+)^td_x(h_-,g_+)^t}\leq \tfrac{d_x(g_-,g_+)^t}{d_x(g_-,h_+)^td_x(h_-,g_+)^t}.		\]
		Thus, for every $r>0$, $f$ can be uniquely extended  to \[\Lambda G\setminus {\big(B_{d_x}(h_+,r)\cup B_{d_x}(h_-,r)\big)},
		\]	that  implies that $f$ can be uniquely extended to a continuous map defined on  $\Lambda G\setminus \{h_+,h_-\}$. With the same argument, but choosing  a geodesic preserved by an hyperbolic isometry with extremes distinct  from  $h_+$ and $h_-$,  one can conclude that $f$ can be extended continuously to $\Lambda G$. Using the same arguments one can show  that the map $f^{-1}:f(\Lambda G)\rightarrow \Lambda G$ is defined. By definition the map $f$ is $G$-equivariant on $\Lambda^+G$ and by density and \cref{mobcontinua}, $f$ is $G$-equivariant on $\Lambda G$. 
		
		Let $\xi,\eta\in\Lambda G$ and and let $x$ and $y$ be contained in the geodesics connecting $\xi$ and $\eta$ and $f(\xi)$ and $f(\eta)$ respectively. If $t<1$ (resp. $t\geq1$),  consider the metrics $(d_x^{\eta})^t$ (resp. $d_x^{\eta}$) and the pullback of the metric $d_y^{f(\eta)}$ (resp. $\big(d_y^{f(\eta)}\big)^{{1}/{t}}$).  The arguments above show that both metrics are Möbius-equivalent, thus by Lemma 2.1 in \cite{buyalomoebiusstructur}, there exists $\lambda>0$ such that for every $a,b\in\Lambda G\setminus\{\eta\}$, \[d_x^{\eta}(a,b)^t=\lambda d_y^{f(\eta)}(f(a),f(b)).\] 
		By the  aforementioned lemma the result is not specific of  $x$ and $y$. 
			\end{proof}
	
	\subsection{Examples and classification for PU(1,1)}\label{algunosejemplos}
	Throughout this subsection  $n$ will denote   a natural number.
	Denote $\mathrm{Aut}(\mathrm{T}_d)$ the isometry group of the homogeneous tree of degree $d\geq 3.$ If $G=\mathrm{PO}(1,n),\mathrm{Aut}(\mathrm{T}_d)$,  there exists $\ell(\rho)>0$ such that for every $g\in G$, 
	$\ell(\rho(g))=\ell(\rho)\ell(g).$ The scalar $\ell(\rho)$ is called the \textit{displacement} of $\rho$. 
	Let  ${\rho}:G\rightarrow\mathrm{Isom}(\mathbf{H}^\infty_\R)$ be a continuous irreducible representation. If $G=\mathrm{Aut}(\mathrm{T}_d)$, $\ell(\rho)\in(0,\infty)$ (see Theorem C in \cite{burger2005equivariant}) and  if $G=\mathrm{PO}(1,n)$, $\ell(\rho)\in(0,1)$ (see Propositions 2.1 and 2.3 in \cite{monod2014exotic}).
	Moreover, it was shown in Theorem C in \cite{burger2005equivariant} for $G=\mathrm{Aut}(\mathrm{T}_d)$, and in Theorem B and Remark 3.13 in \cite{monod2014exotic} for $G=\mathrm{PO}(1,n)$, with $n>2$,   that the displacement classifies the continuous irreducible real hyperbolic representations. These are particular cases of \cref{dispargclas}. By the Karpelevich-Mostow theorem, if ${\rho}:\mathrm{PO}(1,n)\rightarrow\mathrm{Isom}(\mathbf{H}_\R^m)$ is irreducible and $\ell(\rho)\in(0,1)$, then $m\not\in\N$. 
	
	For the non-elementary continuous complex hyperbolic representations of the groups $\mathrm{Aut}(\mathrm{T}_d)$ and $\mathrm{PO}(1,n)$, with $n>2$,   the displacement is well-defined, it has the same possible values  and the same classification for irreducible representations holds. In both cases the continuous non-elementary   complex hyperbolic  representations preserve a real hyperbolic subspace (containing all the axis of the hyperbolic isometries in $\rho(G)$). This was shown  for $\mathrm{PO}(1,n)$, with $n>2$, in  Lemmas 2.6 and 4.2 in \cite{monod2018notes}. Similar arguments apply for $\mathrm{Aut}(\mathrm{T}_d)$.

	In \cref{pfijaunpunto} it is   shown that  the displacement is  defined for the continuous non-elementary  representations ${\rho}:\mathrm{PU}(1,m)\rightarrow\mathrm{Isom}_\C(\mathbf{H}^{m'}_\C)$.  
	For $m=1$,  $\ell(\rho)\in(0,2)$ (see Lemmas 1.3.2 and 1.3.4  in \cite{stolowicz2022complex}) and for  $m>1$,  $\ell(\rho)\in(0,1)$ (see Proposition 2.2.1 and Theorem 2.2.3 in \cite{tesisdoc}).  By the Karpelevich-Mostow theorem, if  $n>1$ and  ${\rho}:\mathrm{PU}(1,n)\rightarrow\mathrm{Isom}_\C(\mathbf{H}_\C^m)$ is irreducible and   $\ell(\rho)\in(0,1)$, then $m$ is not finite.  
	
	If ${\rho}:\mathrm{PU}(1,1)\rightarrow \mathrm{Isom}_\C(\mathbf{H}_\C^m)$ is continuous and  non-elementary and $\ell(\rho)\in(0,1)\cup(1,2)$, then $m\not\in\N$. A direct proof   can be found in Lemma 1.3.7 in \cite{stolowicz2022complex}. For  $\ell(\rho)=1$, see \cref{horosphcomb} below. 
	The results  above show that the values of  $t$ in  Theorems \ref{kernelpowert} and \ref{kernel2} are optimal.
	
	It was noted in \cite{monod2018notes} that for  $\mathrm{PU}(1,1)$ the displacement  does not classify the hyperbolic representations. In  \cite{stolowicz2022complex} it is shown that in this case two continuous irreducible representations are equivalent if, and only if, they have same displacement and same \textit{argument}. This result turns out to be a particular case of \cref{dispargclas}. 
	
	If ${\rho}:\mathrm{Isom}_\F(\mathbf{H}^n_\F)\rightarrow\mathrm{Isom}_\C(\mathbf{H}^m_\C)$ is non-elementary and continuous,   
	there exists $s\in[-1,1]$ such that for every $a,b,c\in \mathbf{H}^n_\C,$
	\begin{equation}\label{cart}\mathrm{Cart}(f(a),f(b),f(c))=s\mathrm{Cart}(a,b,c), \end{equation}
	where $f:\partial\mathbf{H}^n_\C\rightarrow \partial\mathbf{H}^{m}_\C$ is the equivariant associated to $\rho$ (see  Remark 2.5 in \cite{monod2018notes},  (\ref{cartexte}) and \cref{mapeoequivariante}).  Up to conjugating $\rho$ with an anti-holomorphic isometry, we will always suppose that $s\in[0,1]$. The argument of $\rho$ is defined as $\Arg(\rho)=s.$
	
	\begin{prop}\label{argument}
		If $n>1$ and  $\mathrm{PU}(1,n)\xrightarrow{\rho}\mathrm{Isom}_\C(\mathbf{H}^\infty_\C)$ is continuous and   irreducible, then $\mathrm{Arg}(\rho)\in(0,1)$. 
	\end{prop}
	\begin{proof}
		Suppose by contradiction that $\Arg(\rho)=1$.	The restriction of $\rho$ to a copy of $\mathrm{PU}(1,1)$ is non-elementary, thus   by Lemma 2.1.2 in \cite{stolowicz2022complex} and by   7. and 9. in Proposition 1.3.1 in \cite{stolowicz2022complex},  $\ell(\rho)=1$, but this is a contradiction (see Proposition 2.2.1 and Theorem 2.2.3 in \cite{tesisdoc}). 
		
		It is shown in \cite{stolowicz2022real}  that if $n>1$, then $\mathrm{PU}(1,n)$ does not admit continuous non-elementary representations in a real hyperbolic space, or in other words,   $\mathrm{Arg}(\rho)>0.$
	\end{proof}

	By deforming the tautological kernel of Möbius  type  associated to  $\mathrm{PO}(1,n)$ (see \cref{respectirreducible}) one obtains another model for the representations constructed in \cite{monod2014exotic}. By \cref{cartant1}, the argument used in the proof of \cref{horosphcomb} and by   \cref{rescalamientodeldesplazamiento}, if one deforms the tautological representations of $\mathrm{PU}(1,n)$, with $n>1$, one obtains for every $t\in(0,1)$  representations ${\rho_t}:\mathrm{PU}(1,n)\rightarrow\mathrm{Isom}_\C(\mathbf{H}^\infty_\C)$ such that $\ell(\rho_t)=t$ and $\Arg(\rho_t)=t.$   After \cref{dispargclas}  we know that  these representations where described   in  \cite{monod2018notes}, that  to the best of the author's knowledge,   is the first reference to  infinite-dimensional continuous irreducible representations  of $\mathrm{PU}(1,n)$, for $n>1$. 
	
	In Section 2.3 in \cite{stolowicz2022complex} the \textit{horospherical combination} of two irreducible representations of $\mathrm{PU}(1,1)$ with the same displacement is defined. Given $\rho$ and $\tau$ such that $\ell(\rho)=\ell(\tau)=t$, then for every $u\in(0,1)$, there exists a representation $\rho\underset{u}{\wedge}\tau$  such that $\ell(\rho\underset{u}{\wedge}\tau)=t$ and \begin{equation}\label{combconv}\Arg(\rho\underset{u}{\wedge}\tau)=(1-u)\Arg(\rho)+u\Arg(\tau).\end{equation}
	Combining in this way  the families of representations constructed in \cite{monod2018notes,monod2014exotic}, it was described, in \cite{stolowicz2022complex},   a 2-parameter family of representations. 
This result is extended	using the methods introduced in this work.  The following theorem, together with \cref{dispargclas},  gives a classification for  the continuous irreducible complex hyperbolic representations of $\mathrm{PU}(1,1)$.  
	
	Let $A\subset \R^2$ be the set such that $(t,s)\in A$ if, and only if,  $t\in(0,1)$ and $0\leq s\leq t$,
	$t=1$ and $0\leq s<1$, or
	$t\in(1,2)$ and $0\leq s\leq 2-t.$ 
	
	\begin{teo}\label{horosphcomb}
		For every $(t,s)\in A$ there exists,  up to conjugation,  a unique continuous irreducible  representation \[\mathrm{PU}(1,1)\xrightarrow{\tau_{t,s}}\mathrm{Isom}_\C(\mathbf{H}^\infty_\C)\] such that $\ell(\tau_{t,s})=t$ and $\Arg(\tau_{t,s})=s$. Moreover, for every $\rho$ irreducible complex hyperbolic infinite-dimensional representation of $\mathrm{PU}(1,1)$, $(\ell(\rho),\Arg(\rho))\in A$.  
	\end{teo}
	\begin{proof}
		Let $\rho:\mathrm{PU}(1,1)\rightarrow \mathrm{Isom}_\C(\mathbf{H}^\infty_\C)$ be a continuous irreducible representation.  Fix $\xi\in\partial\mathbf{H}^1_\C$ and $x\in\mathbf{H}^1_\C$ and  let $f:\partial\mathbf{H}^1_\C\rightarrow\partial\mathbf{H}^\infty_\C$ be the equivariant map associated . Then  there exists $y\in \mathbf{H}^\infty_\C$  such that for every $\eta_1,\eta_2\in\partial\mathbf{H}^1_\C\setminus\{\xi\}$ (see \cref{mapeoequivariante}), \[d_x^\xi(\eta_1,\eta_2)^{\ell(\rho)}= d_y^{f(\xi)}(f(\eta_1),f(\eta_2)).\] Let $\tilde{f}$ be a lift of $f$ such that for every $\eta\in\partial\mathbf{H}^1_\C\setminus\{\xi\}$, $B(\tilde{f}(\xi),\tilde{f}(\eta))=1$. Observe that \[\Arg(\rho)\tfrac{\pi}{2}=|\cart(f(\xi),f(\eta_1),f(\eta_2))|=|\Arg(B(\tilde{f}(\eta_1),\tilde{f}(\eta_2)))|.\]
		By \cref{isometriakoran},  there exist homeomorphisms 
		$h_2:\R\rightarrow \partial\mathbf{H}^1_\C\setminus\{\xi\}$ and $h_1:\partial\mathbf{H}^\infty_\C\setminus\{f(\xi)\}\rightarrow\mathbf{H}_\infty$ such that $h_1\circ f\circ h_2(0)=e$ and such that  for every $x_1,x_2\in\R$,
		\[\begin{array}{rcl}
			B_k\Big(h_1\circ f\circ h_2(x_1),h_1\circ f\circ h_2(x_2)\Big)&=&B\Big(\tilde{f}\circ h_2(x_1),\tilde{f}\circ h_2(x_2)\Big)\\
			&=&e^{\pm i\frac{\Arg(\rho)}{2}}\left|B\Big(\tilde{f}\circ h_2(x_1),\tilde{f}\circ h_2(x_2)\Big)\right|\\
			&=&e^{\pm i\frac{\Arg(\rho)}{2}}d_y^{f(\xi)}(f\circ h_2(x_1),f\circ h_2(x_2))^{2}\\
			&=&e^{\pm i\frac{\Arg(\rho)}{2}}d_x^{\xi}(h_2(x_1),h_2(x_2))^{2t}\\
						&=&e^{\pm i\frac{\Arg(\rho)}{2}}|x_1-x_2|^{t}.
		\end{array}\]
		Thus, the  map given by $x\mapsto B_k(h_1\circ f\circ h_2(x),e)$ is a function of complex negative type defined on $\R$ that is homogenous of degree $\ell(\rho)$. By \cref{productotorcido},  $(\ell(\rho),\Arg(\rho))\in A$.
		 
		 For every $(t,s)\in A$, with $t\leq1$,  the representations $\tau_{t,s}$  were described in Theorem 2.3.1 in \cite{stolowicz2022complex}.
			If $t\in (1,2)$, let $\Gamma$ be the tautological kernel of Möbius type associated to $\mathbf{H}^1_\C$ and let  $\rho_t$ be the irreducible representation associated to the $\Gamma^t$ (see Theorems  \ref{kernel2} and \ref{respectirreducible}).
			Consider the representation $\tau_{t,0}$ (see Theorem 2.3.1 in \cite{stolowicz2022complex}).
			 Observe that, with a slight abuse of notation, $\rho_t=\tau_{t,2-t}$ (see \cref{kernel2} and Propositions  \ref{cartant2} and \ref{rescalamientodeldesplazamiento}) and  that    for every $s\in[0,1]$,
		\[\tau_{t,s(2-t)}=\tau_{t,0}\underset{s}{\wedge}\tau_{t,2-t}.\qedhere\] 
	\end{proof}
	
	Regarding the representations of $\mathrm{PU}(1,m)$ with $m>1$, the author does not know if those described for the first time in $\cite{monod2018notes}$ and  revisited in this work are a complete list, up to conjugation.

	Given a group $G<\mathrm{Isom}_\F(\mathbf{H}^m_\F)$ and a representation ${\rho}:G\rightarrow\mathrm{Isom}_\F(\mathbf{H}^n_\F)$, is natural to  ask weather $\rho$ is \textit{type preserving};   weather $\rho(g)$ is hyperbolic (resp. parabolic, elliptic), when $g$ is hyperbolic (resp. parabolic, elliptic). 
	
	Nicolas Monod $\&$ Pierre Py  showed that the non-elementary representations  of $\mathrm{PO}(1,m)$ are type preserving (see Proposition 2.1 in
	\cite{monod2014exotic} and Proposition 4.4 and Theorem 5.5 in \cite{monodpyselfrepresentations}). 
	The arguments used  in the following proposition are adapted from those in \cite{monod2014exotic} and the statements about fixed points are contained in Proposition 5.1 in \cite{monod2018notes}.  
	\begin{prop}\label{pfijaunpunto}
		If $m>1$  and  $\mathrm{Isom}_\C(\mathbf{H}^m_\C)\xrightarrow{\rho
		}\mathrm{Isom}_\C(\mathbf{H}^{m'}_\C)$  is a continuous  non-elementary representation, then $\rho$ is type preserving and for every $q\in\partial\mathbf{H}^m_\C$,  the stabiliser of $q$ in $\mathrm{Isom}_\C(\mathbf{H}^m_\C)$ fixes a unique point in $\partial\mathbf{H}^{m'}_\C$. Moreover, the displacement and argument of $\rho$ are defined   and $\ell(\rho)\in(0,1]$ and $\Arg(\rho)\in(0,1].$  
	\end{prop}
	\begin{proof}	
		If $g\in\mathrm{Isom}_\C(\mathbf{H}^m_\C)$ is elliptic, by the Cartan fixed-point theorem if $m<\infty$ (see Proposition II.2.7 in \cite{bridson2013metric}) or by Proposition 5.1 in \cite{monod2018notes} if $m\not\in\N$,  $\rho(g)$ is elliptic. By Proposition 5.1 in \cite{monod2018notes}, the hyperbolic type is preserved and  for every  $q\in\partial\mathbf{H}^m_\C$, $\mathrm{Stab}(q)$ fixes a unique point in $\partial\mathbf{H}^{m'}_\C$.
		
	The only type missing is the parabolic. 
		Let  $q_1,q_2\in\partial\mathbf{H}_\C^m$ be  distinct  and  let  $\xi_1,\xi_2$ be   respective lifts such that $B(\xi_1,\xi_2)=1$. 
	 With respect to the decomposition 	$L=\C\xi_1\oplus\C\xi_2\oplus(\xi_1^\perp\cap\xi_2^\perp)$, every $g\in\mathrm{PU}(1,m)$ fixing $q_1$ admits a linear representative 
		\begin{equation}\label{matriz}g(\lambda,v,A,b)=\begin{pmatrix}
				\lambda&-\lambda\tfrac{B(v,v)}{2}+ib&-\lambda B(A(\meddot),v)\\
				0&\lambda^{-1}&0\\
				0&v&A
			\end{pmatrix},   	
		\end{equation}
		where $\lambda>0$, $b\in\R$,  $v\in \xi_1^\perp\cap\xi_2^\perp$ and $A$ is a unitary  map of $\xi_1^\perp\cap\xi_2^\perp$.  
		
		It is immediate that every  map  $g(1,v,\mathrm{Id},b)$ represents a parabolic isometry.    
		By Corollary 5.2 in \cite{monod2018notes}, the restriction of $\rho$ to a copy of  $\mathrm{Isom}_\F(\mathbf{H}^n_\F)$  is non-elementary and every map $g(1,v,\mathrm{Id},0)$ is contained in the isometry group of a copy of $\mathrm{H}^2_\R$. With a slight abuse of notation, for every $b\neq0$ and every $v\neq0$, $\rho(1,0,\mathrm{Id},b)$ and $\rho(1,v,\mathrm{Id},0)$ are parabolic (see Proposition 1.2.2 in \cite{stolowicz2022complex} and  Proposition 2.1 in \cite{monod2014exotic}). 
		
		Let $p\in\mathbf{H}^{m'}_\C$ be the unique point fixed by $\mathrm{Stab}(q_1)$
		and let $F:\mathrm{Stab}(q_1)\rightarrow \R$ be the map  given by 
		$F(g)=B_{p}(\rho(g))$, where $B_{p}$ is the Busemann quasi-character associated to $p$. The map $F$ is a continuous homomorphism (see the comments before  \cref{commutatorhyperbolic}). By (\ref{buseman}) and the comments before \cref{commutatorhyperbolic}, if $g(\lambda,v,A,b)$ represents a non-hyperbolic isometry, then $\lambda=1$. Thus, as  
				 \[g(\lambda,v,A,b)=g(\lambda,0,\mathrm{Id},0)g(1,v,\mathrm{Id},0)g(1,0,\mathrm{Id},\lambda^{-1}b)g(1,0,A,0), \]  with a slight abuse of notation,   $F(g(\lambda,v,A,b))=F(\lambda,0,\mathrm{Id},0)$. Hence, the displacement of $\rho$ is defined and  $\ell(\rho)\in(0,1]$ (see Proposition 2.3 in \cite{monod2014exotic}).
				
	Let $\sigma:\R^2\times\xi_1^\perp\cap\xi_2^\perp\rightarrow\mathbf{H}^m_\C$ be the horospherical parametrisation (see (\ref{horosphericalcoordinates})). If $g\in\ker(F)$, then  $d(g^i\sigma(0,0,0),\sigma(0,0,0))\to\infty$ if, and only if, $(g^iq_2)\to q_1$ (see (\ref{distanciahoroesferica})). Thus, for such a $g$, $g$ is parabolic if, and only if,  $(g^iq_2)\to q_1$.  As  $f:\partial\mathbf{H}^m_\C\rightarrow\mathbf{H}^{m'}_\C$, the equivariant map associated to $\rho$, preserves this property,  then  the parabolic type is preserved (see \cref{mapeoequivariante}). 
				 
				The claim about the argument   follows from  Propositions  \ref{mapeoequivariante} and  \ref{argument} and  the fact that any six elements in $\partial\mathbf{H}^m_\C$ are contained in some $\partial\mathbf{H}^n_\C$.	\end{proof}
	\begin{prop}
		If $m>1$ and  $\mathrm{Isom}_\C(\mathbf{H}^m_\C)\xrightarrow{\rho
		}\mathrm{Isom}_\C(\mathbf{H}^{m'}_\C)
		$ is an irreducible representation such that $\ell(\rho)=1$ or $\Arg(\rho)=1$, then $m=m'$ and $\rho$ is conjugated to the identity. -	\end{prop}
	\begin{proof}		
		The strategy used here is adapted from  \cite{monod2014exotic}.  If $\mathrm{Arg}(\rho)=1$, by Lemma 1.3.1 and  Proposition 2.1.3 in \cite{stolowicz2022complex}, $\ell({\rho})=1$ (see Corollary 5.2 in \cite{monod2018notes}). 
		Hence, in both cases the  proposition follows from  Remark 2.1.7, Proposition 2.1.11 and Theorem 2.3.3 in \cite{tesisdoc}. 
	\end{proof}
	
	\subsection{Restricting  representations}\label{restrictionofrepresentations}
	 	For every  non-elementary subgroup $\Gamma$ of $G=\mathrm{Aut}(\mathrm{T}_d), \mathrm{Isom}_\F(\mathbf{H}^n_\F)$ and  every  irreducible representation ${\rho}:G\rightarrow\mathrm{Isom}_\F(\mathbf{H}^m_\F)$,  the restriction of $\rho$ to $\Gamma$ is non-elementary ($\Lambda_\rho\Gamma$ is infinite.) In this subsection we give sufficient conditions for those restrictions  to be  \say{infinite-dimensional}. 
	This question is addressed   by studying the dimension of the real Hilbert spaces or the Heisenberg groups associated to  the $\F$-kernels of negative type coming from   \say{deformations} of $\F$-kernels of Möbius type. The fact that  these representations are infinite-dimensional comes from metric obstructions.
	
	It is  shown at the end of this subsection that the properties  \textit{strongly discrete}, \textit{convex-cobounded} or \textit{geometrically-finite} are preserved by the representations coming from the \say{deformations} of the tautological kernels of Möbius  type. For further reading on these properties and their consequences see \textit{e.g.} \cite{das2017geometry} and the references therein. 
	
I would like to thank Nicola Cavallucci for bringing \cite{ledonnesnowflakes} to my attention, as well as for encouraging me to seek a metric obstruction in order to prove  that certain representations are infinite-dimensional.

		\begin{teo}\label{noelementalarbol}
		For $t\in(0,\infty)$, denote $\mathrm{Aut}(\mathrm{T}_d)\xrightarrow{\rho_t}\mathrm{Isom}(\mathrm{H^\infty_\R})$ the irreducible representation such that $\ell(\rho_t)=t.$ If $G<\mathrm{Aut}(\mathrm{T}_d)$ is non-elementary, then for every $t$,  the irreducible part of $\rho_t|_G$ is infinite-dimensional.
	\end{teo}
	\begin{proof}
		Denote $\Gamma_{2t}$ the tautological real kernel of Möbius type associated to $\Lambda_{\rho_{2t}}\mathrm{Aut}(\mathrm{T}_d)\subset\partial\mathbf{H}^\infty_\C.$ Observe that the irreducible representation $\tau_t$ associated to $\Gamma_{2t}^{\frac{1}{2}}$ is equivalent to $\rho_t$. Given $y_0,x_0,x_1\in \Lambda_{\rho_{2t}}G$ pairwise distinct, if  
		${\Psi'_{x_0,x_1}}:\big(\Lambda_{\rho_{2t}}G\setminus y_0\big)^2\rightarrow\R$ is the real kernel of negative type associated to the restriction of $\Gamma_{2t}^{\frac{1}{2}}$ to $\Lambda_{\rho_{2t}}G$, then by definition, $\Psi'_{x_0,x_1}=\Psi^{\frac{1}{2}}_{x_0,x_1}$, where ${\Psi_{x_0,x_1}}:\big(\Lambda_{\rho_{2t}}G\setminus y_0\big)^2\rightarrow\R$ is the real kernel of negative type associated to the restriction of $\Gamma_{2t}$ to $\Lambda_{\rho_{2t}}G$. By Theorem 1.2 in \cite{ledonnesnowflakes}, the irreducible part of $\rho_t|_G$ is infinite-dimensional.	
	\end{proof}
With the similar arguments one can show the following theorem. 
	\begin{teo}\label{noelementalreal}
		If $\mathrm{PO}(1,n)\xrightarrow{\rho}\mathrm{Isom}(\mathbf{H}^\infty_\R)$ is an irreducible representation, then for every non-elementary $G<\mathrm{PO}(1,n)$, the irreducible part of $\rho|_G$ is infinite-dimensional. 
	\end{teo}
	
In Theorem 1.2 in \cite{ledonnesnowflakes} the authors showed that if $(X,d)$ is an infinite metric space and $t\in(0,1)$, then $(X,d^t)$ cannot be embedded isometrically in any finite-dimensional Euclidean space. 	A non-embeddability result of this kind in  the finite-dimensional Heisenberg groups for infinite  $t$-snowflaked metric spaces (for $t<\frac{1}{2}$)  is not known by the author. The  following proposition is a partial result in that direction that implies for the complex case an analogous result to \cref{noelementalreal}.
	\begin{teo}\label{infinito}
		Let   $G<\mathrm{PU}(1,n)$ be a non-elementary subgroup and  let  $q_1, q_2\in \Lambda G$ be distinct.  Suppose that   $\Lambda G\setminus\{q_1\}$ is identified with  $ X\subset \mathbf{H}_{n-1}$, using  $\xi_1$ and $\xi_2$, two  respective lifts of $q_1$ and $q_2$ such that $B(\xi_1,\xi_2)=1$.   Denote $\Psi$  the  complex kernel of  negative type associated to $X$. If
			 $t<\frac{1}{2}$ when $n>1$, or if    $t<1$ when $n=1$,  
		then the Heisenberg group associated to $\Psi^t$ is infinite-dimensional.
	\end{teo}
	\begin{proof}
		The map $\phi:X^2\rightarrow\C$ given by  \[\phi(x,y)=\frac{1}{2}\Big(\Psi(x,e)^t+\Psi(e,y)^t -\Psi(x,y)^t\Big)\]
		is a kernel of complex positive type (see \cref{positivenegative}).
		The claim is that the complex Hilbert space associated to this kernel  by the GNS construction (see \cref{GNSpos}) is infinite-dimensional. If $\mathbf{H}_{m-1}=N\times\R$  is the Heisenberg group associated by the GNS construction to $\Psi^t$,   where $(N,\langle\,\meddot\,,\meddot\,\rangle)$ is a complex Hilbert space, then for every  $(u,b),(v,d)\in X$, 	$\phi\big((u,b),(v,d)\big)=\langle u,v\rangle$.

		Let $X^+\subset X$ be the set identified with $\Lambda^+G$. 
		For $x\in X$ denote $\phi_x=\phi(x,\meddot\,)$  and define  $A=\{x\in X^+\mid \phi_x\neq0\}$. The first claim is that $A$ is infinite.
		Indeed, by contradiction suppose that  $A$ is finite. Then there exist  $r>0$  and an injective map such that  \[X^+\xrightarrow{\iota}\bigcup\limits_{j=1}^r\{u_j\}\times \R\]
		and 	such that for every $x,y\in X^+$,
		$d_k(\iota(x),\iota(y))=d_k(x,y)^t$ (see the proof of \cref{GNSpos}).  
		For every $j$, $\big(\{u_j\}\times\R,d_k\big)$ is isometric to $(\R,d_{Euc}^{\frac{1}{2}})$ and as $X^+$ is infinite, there exists $u_j$ such that $(\{u_j\}\times\R) \cap\iota(X^+)$ is infinite, but this is a contradiction (see Theorem 1.2 in \cite{ledonnesnowflakes}).

		Suppose by contradiction that $N$ is finite-dimensional. 
		By the GNS construction (see \cref{GNSpos}) there exist $x_i\in A$ and  $\lambda_i\in\C$, with $i=1,\dots,s$,  such that
		$\sum_{i=1}^s\lambda_i\phi_{x_i}=0.$
		This implies that there exists a constant $c_0$ such that
		\begin{equation}\label{lincomb}
				\tfrac{1}{2}\Big(\sum\limits_{i=1}^s\lambda_{x_i}
				\Big)\Psi(e,\meddot\,)^t - \tfrac{1}{2}\sum\limits_{i=1}^s\lambda_{x_i}\Psi(x_i,\meddot\,)^t=
				c_0.\end{equation}
		The claim  is that this is a contradiction. 
	If $t>0$,   on the one hand, every function ${\Psi^t_x}:X\rightarrow\C$, given by ${\Psi^t_x}(y)=\Psi(x,y)^t$, is the restriction to $X$ of a function defined on  $N\times \R$ that is differentiable away from $x$. 
		On the other hand, the claim is that for $x\in A$,  $\Psi_x^t$ is not \say{differentiable} in  $x$, that is to say, there exists a sequence  $(y_i)$ in $ X$ such that the set \[\left\{\tfrac{\Psi(x,y_i)}{\|y_i-x\|_{N\times\R}}\right\}_i \] is unbounded. 
			Let  $x\in \Lambda^+ G$ be the attracting point of a hyperbolic isometry $\delta$.   Denote $y$ the repelling point of $\delta$ and take representatives $\eta_1$, $\eta_2$ of $y$ and $x$ respectively such $B(\eta_2,\xi_1)=1$ and such  that $B(\eta_1,\eta_2)=1$. Consider  $p=[a\eta_1+\eta_2+u]\in \Lambda^+G$, where  $u\in\eta_1^\perp\cap\eta_2^\perp$ and $a\in\C$. There exist $\lambda>0$ and $T$, a unitary map of $\eta_1^\perp\cap\eta_2^\perp$, such that for every $n\in\N$,
		\[\delta^n(p)=\left[\lambda^{-n}a\eta_1+\lambda^n\eta_2+T^n(u)\right]. \]
		Denote
		$\mu_n=\lambda^{-n}a\eta_1+\lambda^n\eta_2+T^n(u)$ and 
		recall that for every $[\mu],[\nu]\in \Lambda G$
		\[\Psi([\mu],[\nu])=\tfrac{B(\mu,\nu)}{B(\mu,\xi_1)B(\xi_1,\nu)}.\]
		With a slight abuse of notation, the claim is that \[\left\{ \tfrac{\Psi\big(\delta^n(p),x\big)^t}{\|\delta^n(p)-x\|_{N\times\R}}   \right\}_{n\in\N}\] is  unbounded. When $n=1$, declare $\C^0=\{0\}$.
		Here ${\|\delta^n(p)-x\|_{N\times\R}}$ denotes the distance with respect to  the usual norm in $N\times\R$ and  the identification 
		\[(v,b)\mapsto \left[\left({-\|v\|^2}+ib\right)\xi_1+\xi_2 +\sqrt{2}{v}\right].\] 
		Observe that 
		\begin{equation}\label{elkernelt}
				\Psi\big(\delta^n(p),x\big)=
				\tfrac{B(\mu_n,\eta_2)}{B(\mu_n,\xi_1)}=
				\tfrac{a}{\lambda^{2n}\big( \lambda^{-2n}aB(\eta_1,\xi_1)+1+\lambda^{-n}B(T^n(u),\xi_1)  \big)}.
		\end{equation}
			If $\delta^n(p)$ is identified with $(w_n,d_n)\in N\times\R$, then \[\sqrt{2}w_n=\tfrac{\mu_n}{B(\mu_n,\xi_1)}- \tfrac{B(\mu_n,\xi_2)}{B(\mu_n,\xi_1)}\xi_1-\xi_2\] and 
		\[d_n=\mathrm{Im}\left(\tfrac{B(\mu_n,\xi_2)}{B(\mu_n,\xi_1)} \right).\]
		Similarly, if $x$ is identified with $(z,c)$, then 
		\[\sqrt{2}z={\eta_2}- {B(\eta_2,\xi_2)}\xi_1-\xi_2\] and 
		$c=\mathrm{Im}\left({B(\eta_2,\xi_2)} \right).$
		Thus,  
		\begin{equation}\label{limites1}				2\|w_n-z\|^2 = 
				\left\| \tfrac{\mu_n}{B(\mu_n,\xi_1)}+\left(B(\eta_2,\xi_2)-\tfrac{B(\mu_n,\xi_2)}{B(\mu_n,\xi_1)}\right)\xi_1- \eta_2 
				\right\|^2=								
				2\mathrm{Re}\left(\tfrac{B(\mu_n,\eta_2)}{B(\mu_n,\xi_1)}\right).					\end{equation}	
		
		If one denotes   $B(\mu_n,\xi_1)=\lambda^n +A_n$ and $B(\mu_n,\xi_2)=\lambda^nB(\eta_2,\xi_2)+B_n$, then \[\lim\limits_{n\to\infty}\lambda^{-n}A_n=0=\lim\limits_{n\to\infty}\lambda^{-n}B_n.\]
		Observe that  
		\begin{equation}\label{limites2} \begin{array}{rcl}d_n-c&=&
			\mathrm{Im}\left(	\frac{B(\mu_n,\xi_2)}{B(\mu_n,\xi_1)}\right)-\mathrm{Im}(B(\eta_2,\xi_2)) \\&=&
				\frac{\mathrm{Im}\Big(\big(\lambda^nB(\eta_2,\xi_2)+B_n\big)(\lambda^n+\overline{A_n}) -|B(\mu_n,\xi_1)|^2B(\eta_2,\xi_2)\Big)}{|B(\mu_n,\xi_1)|^2}
				\\&=&
				\frac{\mathrm{Im}\Big(\lambda^{2n}B(\eta_2,\xi_2)+\lambda^n\big(B(\eta_2,\xi_2)\overline{A_n}+B_n\big)+B_n\overline{A_n}-\big(\lambda^{2n}+2\lambda^n\mathrm{Re}(A_n)+|A_n|^2\big)B(\eta_2,\xi_2)\Big)}{\lambda^{2n}+2\lambda^n\mathrm{Re}(A_n)+|A_n|^2}\\&=&
				\frac{\mathrm{Im}\Big(\lambda^n\big(B(\eta_2,\xi_2)(\overline{A_n}-2\mathrm{Re}(A_n))+B_n\big)+B_n\overline{A_n}-|A_n|^2B(\eta_2,\xi_2)\Big)}{\lambda^{2n}+2\lambda^n\mathrm{Re}(A_n)+|A_n|^2}.	\end{array}\end{equation}
		Thus if $n>1$ and  $t<\frac{1}{2}$, 	  by 	(\ref{limites1}),
		\[\lim\limits_{n\to\infty}\lambda^{4nt}\|w_n-z\|^2=	\lim\limits_{n\to\infty}\lambda^{n(4t-2)}\mathrm{Re}\left(\tfrac{a}{a\lambda^{-2n}B(\eta_1,\xi_1)+ 1+\lambda^{-n}B(T^n(u),\xi_1)}\right)=0		\]							and by 	(\ref{limites2}),							$\lim_{n\to\infty}\lambda^{2nt}\left(d_n-c\right)=0.									$									These   computations, together with  (\ref{elkernelt}),  show that  the set 									\[\left\{ \tfrac{\Psi^t\big(\delta^n([\xi_2]),[\eta_2]\big)}{\|\delta^n([\xi_2])-x\|_{N\times\R}}   \right\}_{n\in\N}\] 									is unbounded, proving that   (\ref{lincomb}) is a contradiction.
		
		If $n=1$, 
		$A_n=\lambda^{-n}aB(\eta_1,\xi_1)$ and $B_n=\lambda^{-n}aB(\eta_1,\xi_2)$. Thus,  (\ref{limites2}) takes the following shape, 
		\begin{equation*} \begin{array}{rcl}\frac{B(\mu_n,\xi_2)}{B(\mu_n,\xi_1)}-B(\eta_2,\xi_2) &=&						\frac{\lambda^n\big(-B(\eta_2,\xi_2)A_n+B_n\big)+B_n\overline{A_n}-|A_n|^2B(\eta_2,\xi_2)}{\lambda^{2n}+2\lambda^n\mathrm{Re}(A_n)+|A_n|^2}\\&=&
				\frac{a\big(-B(\eta_2,\xi_2)B(\eta_1,\xi_1)+B(\eta_1,\xi_2)\big)+B_n\overline{A_n}-|A_n|^2B(\eta_2,\xi_2)}{\lambda^{2n}+2\lambda^n\mathrm{Re}(A_n)+|A_n|^2}. 
		\end{array}\end{equation*}
		This shows that  
		\[\begin{array}{rcl}										\lim\limits_{n\to\infty}\lambda^{2nt}\left(\frac{B(\mu_n,\xi_2)}{B(\mu_n,\xi_1)}-B(\eta_2,\xi_2)\right)&=&0									\end{array}\]	
		and therefore, the set\[\left\{ \tfrac{\Psi^t\big(\delta^n([\xi_2]),[\eta_2]\big)}{\|\delta^n([\xi_2])-x\|_{\R}}   \right\}_{n\in\N}\] is not bounded, proving that (\ref{lincomb}) is a contradiction. 
	\end{proof}
	
	With the previous proposition at hand and with similar arguments to those used in \cref{noelementalarbol}, one can prove the following theorem. 
	
	\begin{teo}\label{noelementalcomplejo}
		Let $G<\mathrm{PU}(1,n)$ be non-elementary and denote $\mathrm{PU}(1,n)\xrightarrow{\rho_t}\mathrm{Isom}_\C(\mathbf{H}^\infty_\C)$ the irreducible representations constructed by deforming  the tautological complex kernel of hyperbolic type associated to $\partial\mathbf{H}^n_\C$. If
		$t<\frac{1}{2}$ when $n>1$, or if    $t<1$ when $n=1$,
		 then the irreducible part of $\rho_t|_G$ is infinite-dimensional. 
	\end{teo}
	For $S\subset \mathbf{H}^m_\F\cup\partial\mathbf{H}^m_\F$, denote $\mathrm{Hull}(S)$ the smallest closed \textit{convex} subset of $\mathbf{H}^m_\F\cup\partial\mathbf{H}^m_\F$ that contains $S$. A  set  is called \textit{convex} if for every pair of points (possibly at infinity) there exists a geodesic or a segment of a geodesic connecting them. For an explicit construction of $\mathrm{Hull}(S)$ see Proposition 7.5.2 in \cite{das2017geometry}. If $G<\mathrm{Isom}_\F(\mathbf{H}^m_\F)$ is  non-elementary, denote $C_G=\mathrm{Hull}(\Lambda G)\cap \mathbf{H}^m_\F.$ Observe that $C_G$ is a closed, convex and $G$-invariant subset of $\mathbf{H}^m_\F$. If $\Lambda G$ is compact, then $C_G$ is locally compact (see Proposition 7.7.2 in \cite{das2017geometry}). 
	
	Given a  group $G<\mathrm{Isom}_\F(\mathbf{H}^m_\F)$ and a point $x_0\in \mathbf{H}^m_\F$,  for every $s>0$ denote \[\Sigma_s(G)=\sum\limits_{g\in G}e^{-sd(gx_0,x_0)}\in\R_{>0}\cup\{\infty\}.\]
	The \textit{critical  exponent} of the group $G$ is defined as 
	$\delta_G=\inf\{s\mid \Sigma_sG<\infty\},$
	where the convention $\inf\emptyset=\infty$ is used. The critical exponent is independent of  the choice of $x_0$. For further reading on the properties of $\delta_G$ see \textit{e.g.} Chapter 8 in \cite{das2017geometry} and the references therein. 
	
If  $G<\mathrm{Isom}_\F(\mathbf{H}^n_\F)$ is non-elementary, with respect to any Bourdon metric (see \cref{mobiusinvariant}), $0<\mathrm{dim}_{\mathrm{H}}(\Lambda G)<\infty$. Indeed, on the one hand,  \[\dimension(G)\leq\dimension(\partial\mathbf{H}^n_\F)=k(n+1)-2,\] where $k=\mathrm{dim}_\R(\F)$ (see \textit{e.g.} Remark 2.4 in \cite{magnani}). 
	On the other hand, if $G<\mathrm{Isom}_\F(\mathbf{H}^n_\F)$ is non-elementary,  there exists $G'<G$ discrete (for the compact-open topology) and non-elementary (see Lemma 10.2.2 and the proof of  Proposition 10.5.4 in \cite{das2017geometry}). If   $\Lambda G$ is compact, then $C_{G'}$ is locally compact and by   Theorem 1.2.1 and  Propositions 5.2.7 and  10.5.4 in \cite{das2017geometry},  $\dimension(\Lambda G)\geq\delta_{G'}>0.$ 
	
	If  $G<\mathrm{Isom}_\F(\mathbf{H}^n_\F)$ is non-elementary and   $\rho:G\rightarrow\mathrm{Isom}_\F(\mathbf{H}^m_\F)$  is a representation for which $\ell(\rho)$ is defined, by Lemmas  \ref{supremo}, \ref{semultiplicaladimensiobn}, \ref{mobiusinvariant} and \cref{mapeoequivariante}, 
	 $\mathrm{dim}_{\mathrm{H}}(\Lambda_{\rho}G)=\tfrac{1}{\ell(\rho)}\mathrm{dim}_{\mathrm{H}}(\Lambda G).$
	The results at the end of this subsection, together with this fact, will  show that there are infinite-dimensional \textit{strongly discrete}, \textit{convex-cobounded} or \textit{geometrically finite} subgroups of $\mathbf{H}^\infty_\F$ with limit sets with Hausdorff dimension as large as desired. 
	
	Let   $(x_i)$ be a sequence in $\mathbf{H}^m_\F$ converging to a point $\xi\in\partial\mathbf{H}^m_\F$. The sequence   \textit{converges radially} to  $\xi$ if  $\{x_i\}_{i}$ remains within a bounded distance of a geodesic ray pointing at $\xi$. Whilst $(x_i)\to\xi$ if, and only if, $((x_i,\xi)_{x_0})\to\infty$, for the radial convergence the condition changes; the sequence  converges  radially to $\xi$ if, and only if,
	there exists $R>0$ such that for every $i$, $(x_0,\xi)_{x_i}\leq R$ (see Proposition 7.1.1 in \cite{das2017geometry}). 
	The sequence  \textit{converges horospherically} to $\xi$ if for every  horoball $H$ centred at $\xi$ there exists $N>0$ such that for every $n>N$, $x_n\in H.$
	
	Let $G<\mathrm{Isom}_\F(\mathbf{H}^m_\F)$ be a non-elementary group and fix $x_0\in \mathbf{H}^m_\F.$ Define  the \textit{radial} and \textit{horospherical limit sets} respectively as  
	\[\Lambda_rG=\{\xi\in\Lambda G\mid (g_ix_0)_i\,\,\text{converges radially to }\,\xi,\,\,\text{for some}\,\,g_i\in G\}\]and
	\[\Lambda_hG=\{\xi\in\Lambda G\mid (g_ix_0)_i\,\,\text{converges horospherically to }\,\xi,\,\,\text{for some}\,\,g_i\in G\}. \]
	By the triangle inequality,    $\Lambda_rG$ and $\Lambda_hG$ are independent of  the choice of $x_0$, and from the definitions, it is clear that    both sets are $G$-invariant and  $\Lambda_rG\subset\Lambda_hG$ (see Observation 7.1.4 in \cite{das2017geometry}).
	
	If $G$  and $x_0$ are as before, a point $\xi\in\partial\mathbf{H}^m_\F$ is called \textit{a parabolic fixed point of} $G$ if  $G_\xi$  does not contain hyperbolic elements and $G_\xi\cdot x_0$ is unbounded. The definition is independent of the choice of $x_0$. 	If $\xi\in\partial\mathbf{H}^m_\F$ is a parabolic fixed point of $G$, then  $\xi\in \Lambda G$ (see Observation 6.2.11 in \cite{das2017geometry}).
		A
	set $S\subset \mathbf{H}^m_\F\cup \partial\mathbf{H}^m_\F$ is called $\xi$-\textit{bounded}, if $\xi\not\in \overline{S}$. 	A point $\xi$ is called \textit{bounded parabolic}  if it is a parabolic fixed point of $G$
	and if there exists $S\subset \mathbf{H}^m_\F\cup\partial\mathbf{H}^m_\F $  $\xi$-{bounded} such that $G\cdot x_0\subset G_\xi\cdot S.$
	The definition is independent  of the choice of $x_0$. Indeed, suppose $y\in \mathbf{H}^m_\F$ is such that $d(x_0,y)<R$ and denote 
	\[S(R)=\big(\partial\mathbf{H}^m_\F\cap S\big)\cup \big(S\cap \mathbf{H}^m_\F\big)_R, \]
	where $\big(S\cap \mathbf{H}^m_\F\big)_R$ is the $R$-neighbourhood of $S\cap \mathbf{H}^m_\F$. If $S$ is $\xi$ bounded,  then  $S(R)$ is $\xi$-bounded and 
	$G\cdot y\subset
	G_\xi\cdot S(R).$ 
	Denote \[\Lambda_{bp} G=\{\xi\in \Lambda G\mid \xi\,\,\text{is bounded parabolic}\}.\]   
	
	Let   $\rho:G\rightarrow\mathrm{Isom}_\F(\mathbf{H}^m_\F)$ be a continuous non-elementary representation. The representation $\rho$ is called \textit{convex-cobounded} if $\Lambda_\rho G$ is compact and $\Lambda_{\rho}G=\Lambda_r\rho(G)$.  The representation $\rho$ is called \textit{geometrically finite} if $\Lambda_{\rho}G$ is compact and 
	$\Lambda_{\rho} G=\Lambda_r\rho(G)\cup\Lambda_{bp}\rho(G).$
	
	\begin{lem}
		If $G<\mathrm{Isom}_\F(\mathbf{H}^m_\F)$ is non-elementary and  convex-cobounded, then 
		the action of $G$ on $C_{G}$ is cobounded.
	\end{lem}
	\begin{proof}
By contradiction suppose  that there is no bounded subset $B\subset C_G$ such that 
		$G\cdot B=C_G$. Observe that $\Lambda_h G=\Lambda G$ and fix $x_0\in C_G$. For every $n\in\N_{>0}$, $C_G\not\subset G\cdot B(x_0,n),$  thus there exists a sequence $(x_n)_{n>0}$ in $ C_G$ such that $\inf\limits_{g\in G}\{d(x_n,gx_0)\}\geq n$. Choose $g_n\in G$ such that $$d(x_n,g_n^{-1}x_0)<\inf\limits_{g\in G}\{d(x_n,gx_0)\}+1.$$
		As $\Lambda G$ is compact, up to taking a sequence, one can suppose that $(g_nx_n)_n\to\xi\in\Lambda_h G$ (see Proposition 7.7.2 in \cite{das2017geometry}). Let $b_{\xi,x_0}$ be the Busemann function  centred  at $\xi$ and normalised in  $x_0.$ By Lemma 3.4.10 in \cite{das2017geometry} for every $g\in G$, 
	\[			b_{\xi,x_0}(gx_0)
			=\lim\limits_{n\to\infty}d(gx_0,g_nx_n)-d(x_0,g_nx_n)
			>-1.	\]	
		This is a contradiction because $\xi\in\Lambda_hG.$ 
	\end{proof}
	
	A subgroup $G<\mathrm{Isom}_\F(\mathbf{H}^m_\F)$ is called \textit{strongly discrete} if for every bounded $B\subset \mathbf{H}^m_\F$, the set \[\{g\in G\mid g\cdot B\cap B\neq\emptyset\}\]
	is finite. For a discussion on different definitions of discreteness for subgroups of $\mathrm{Isom}_\F(\mathbf{H}^m_\F)$ see Chapter 5 in \cite{das2017geometry}. 
		If $G$ is strongly discrete, then  $G$ is convex-cobounded if, and only if,  the action of $G$ on $C_{G}$ is cobounded (see Theorem 12.2.7 in \cite{das2017geometry}).
	
	\begin{lem}\label{fesquasiisometrico}
		If $\mathrm{Isom}_\F(\mathbf{H}^n_\F)\xrightarrow{\rho}\mathrm{Isom}_\F(\mathbf{H}^\infty_\F)$ is non elementary, then there exists  $\mathbf{H}^n_\F\xrightarrow{f}\mathbf{H}^\infty_\F$  equivariant and quasi-isometric. 
	\end{lem}
	\begin{proof}
By Proposition 5.7 in \cite{monod2018notes}, there exists $f:\mathrm{H}^n_\F\rightarrow\mathrm{H}^\infty_\F$ equivariant (see the comments before \cref{uniformconvergence}). 		If $x,y\in \mathbf{H}^n_\F$, there exists $g\in \mathrm{Isom}_\F(\mathbf{H}^n_\F)$ hyperbolic such that $gx=y$ and such that $x$ and $y$ belong to the axis of $g$. Hence, \begin{equation}\label{expansiva}d(f(x),f(y))=d\big(f(x), \rho(g)f(x)\big)\geq\ell(\rho(g))=\ell(\rho)d(x,y).\end{equation}
		
		Let $\Gamma<\mathrm{Isom}_\F(\mathbf{H}^n_\F)$ be a uniform lattice (see Chapter XIV in \cite{raghunathan1972discrete}), and  by  Selberg's lemma (see \cite{alperinselberg}),  suppose that $\Gamma$ is  torsion-free.
		The  claim is that of $\rho\vert_\Gamma$  is non-elementary.  Indeed, let   
		$\phi:\partial\mathbf{H}^n_\F\rightarrow\partial\mathbf{H}^\infty_\F$ be  the equivariant map associated to $\rho$ (see \cref{mapeoequivariante}) and 
		observe that  $\Lambda\Gamma=\partial\mathbf{H}^n_\F$ because every point in $\mathbf{H}^n_\F$ remains at a bounded distance from an orbit of $\Gamma$. Thus    $\Lambda_\rho\Gamma=\mathrm{Im}(\phi)$,   showing  that $\rho\vert_\Gamma$ is non-elementary (see Propositions 6.2.14 and 7.3.1  in \cite{das2017geometry}). 
		
		In Theorem 2.3.1 of \cite{korevaari1997global}, Korevaar $\&$ Schoen showed (in a higher generality) that there exists a $\Gamma$-equivariant and $L$-Lipschitz harmonic map $\varphi:\mathbf{H}^n_\F\rightarrow\mathbf{H}^\infty_\F$. For the definition and  further reading on harmonic maps with CAT(0) codomain see also \cite{korevaarsobolev}. 
		
		Suppose $Q\subset \mathbf{H}^n_\F$ is a compact set  such that $\Gamma\cdot Q=\mathbf{H}^n_\F$ and define  $C=\sup\{d(f(q),\varphi(q))\}_{x\in Q}$. To conclude, observe that if $g_i\in \Gamma$ and $x_i\in Q$, then 
		\begin{equation*}\begin{array}{rcl}
				d\big(f(g_1x_1),f(g_2x_2)\big)&\leq&d\big(f(g_1x_1),\varphi(g_1x_1)\big)+d\big(\varphi(g_1x_1),\varphi(g_2x_2)\big)+d\big(\varphi(g_2x_2),f(g_2x_2)\big)\\
				&\leq&Ld(g_1x_1,g_2x_2)+2C.
		\end{array}\findos\end{equation*}\renewcommand{\qedsymbol}{}
	\end{proof}
	\vspace{-\baselineskip}
	\begin{lem}\label{stronglydiscrete} If $\mathrm{Isom}_\F(\mathbf{H}^n_\F)\xrightarrow{\rho}\mathrm{Isom}_\F(\mathbf{H}^\infty_\F)$ is non-elementary and $G<\mathrm{Isom}_\F(\mathbf{H}^n_\F)$ is strongly discrete, then 
		$\rho(G)$ is strongly discrete.\end{lem}
	\begin{proof}
		The first claim  is that $\rho(G)$ is discrete for the bounded-open topology. The composition with a fixed isometry is an open map in the bounded-open topology. By contradiction, suppose that $(\rho(g_\alpha))\to \mathrm{Id}$, where $(g_\alpha)$ is a non-trivial net.  The equivariant map $f:\mathbf{H}^n_\F\rightarrow\mathbf{H}^\infty_\F$ associated to $\rho$ (see \cref{fesquasiisometrico}) is a quasi-isometric. Thus there exist $\lambda,c>0$ such that 
		for every $x\in\mathbf{H}^n_\F$,   
		\[ \lambda d(g_\alpha x,x)\leq d(\rho(g_\alpha)f(x),f(x))+c.\] As $G$ is strongly discrete,  $\{g_\alpha\}_\alpha$ is a finite set.
		
		The set $C_{\rho(G)}$ is a locally compact, convex and closed subspace of $\mathbf{H}^\infty_\F$. The restriction of $\rho(G)$ to $C_{\rho(G)}$ is discrete  for the compact-open topology, thus   $\rho(G)$ is a strongly discrete  group of isometries of $C_{\rho(G)}$ (see Proposition 5.2.7 in \cite{das2017geometry}).    	  
		Let ${\pi}:\mathbf{H}^\infty_\F\rightarrow C_{\rho(G)}$ be the projection on $C_{\rho(G)}$ (see Proposition II.2.4 in \cite{bridson2013metric}) and let  $B\subset\mathbf{H}^\infty_\F$ be a bounded set.  The map $\pi$ does not increase the distance  and is $\rho(G)$-equivariant. If $\gamma\in\rho(G)$ is such that 
		$\gamma\cdot B\cap B\neq\emptyset$, then $\pi(B)$ is a bounded set and 
		\[\pi\big(\gamma\cdot B\cap B)\subset \gamma\cdot\pi(B)\cap\pi(B)\neq\emptyset,\]   
		showing that  \[\{\gamma\in \rho(G)\mid \gamma\cdot B\cap B\neq\emptyset\}\] is a finite set. 
	\end{proof}
	\begin{teo}
		If $\mathrm{Isom}_\F(\mathbf{H}^n_\F)\xrightarrow{\rho}\mathrm{Isom}_\F(\mathbf{H}^\infty_\F)$ is non-elementary and $G<\mathrm{Isom}_\F(\mathbf{H}^n_\F)$ is convex-cobounded (resp. geometrically finite), then 
		$\rho(G)<\mathrm{Isom}_\F(\mathbf{H}^\infty_\F)$ is convex-cobounded (resp. geometrically finite). 
	\end{teo}
	\begin{proof}
		Let $f:\mathbf{H}^n_\F\rightarrow\mathbf{H}^\infty_\F$ and $\phi:\partial\mathbf{H}^n_\F\rightarrow\partial\mathbf{H}^\infty_\F$ be the equivariant maps associated to $\rho.$ As $f$ is a quasi-isometric map (see \cref{fesquasiisometrico}), by Proposition 15 in p. 89 in \cite{ghysgroupes}, there exist $\lambda,C>0$ such that  for every $x,y,z\in\mathbf{H}^n_\F$,  
		\begin{equation}\label{estambienqisometrica}\lambda^{-1}(x,y)_z-C\leq(f(x),f(y))_{f(z)}\leq\lambda(x,y)_z-C.\end{equation}
		
		Recall that by \cref{mapeoequivariante},  $\phi(\Lambda G)=\Lambda_\rho G$. After  (\ref{estambienqisometrica}), it is clear that if $(x_i)$ converges radially to $\xi$, then $(f(x_i))$ converges radially to $\phi(\xi)$, and therefore, $\phi(\Lambda_r G)\subset \Lambda_r\rho(G).$
		
		Let $\xi\in \partial\mathbf{H^n_F}$ and $S\subset \mathbf{H^n_F}\cup\partial \mathbf{H^n_F}.$ By \cref{mapeoequivariante}, \cref{fesquasiisometrico} and (\ref{estambienqisometrica}), if $S$ is $\xi$-bounded 
		then  \[\phi(S\cap\partial \mathbf{H^n_F})\cup f(S\cap \mathbf{H^n_F})\]is $\phi(\xi)$-bounded   and 	$\phi(\Lambda_{bp} G)\subset\Lambda_{bp}\rho(G).$	  	 
	\end{proof}

	\bibliographystyle{plain}{}
	\addcontentsline{toc}{section}{Bibliography}
	\bibliography{biblio}
	
\end{document}